\documentclass{amsart}

\usepackage{amssymb,amsmath,epsfig,graphics,latexsym,psfrag}
\usepackage[english]{babel}
\usepackage[all]{xy}
\newtheorem{theorem}{Theorem}[section]
\newtheorem{proposition}[theorem]{Proposition}
\newtheorem{corollary}[theorem]{Corollary}
\newtheorem{lemma}[theorem]{Lemma}
\newtheorem{claim}[theorem]{Claim}

\theoremstyle{definition}
\newtheorem{definition}[theorem]{Definition}

\theoremstyle{remark}
\newtheorem{remark}[theorem]{Remark}

\numberwithin{equation}{section}
\renewcommand{\t}{ \tilde}
\renewcommand{\b}{ \partial}
\newcommand{\Z}{\bf Z}
\newcommand{\R}{\bf R}
\newcommand{\N}{\bf N}

\newcommand{\Q}{\bf Q}
\newcommand{\Hi}{\bf H}

\newcommand{\C}{\bf C}

\renewcommand{\S}{\bf S}
\renewcommand{\l}{\langle}
\renewcommand{\r}{\rangle}

\renewcommand{\int}[1]{\stackrel{\circ}{#1}}
\renewcommand{\o}{\overline}
\newcommand{\abs}[1]{\lvert#1\rvert}
\renewcommand{\H}[1]{H_1(#1;{\Z})}

\newcommand{\co}{\colon\thinspace}
\renewcommand{\epsilon}{\varepsilon}

\usepackage{times}

\begin{document}
\sloppy

\title[]{Topological rigidity  and Gromov simplicial volume}

\author{Pierre Derbez}
\address{Laboratoire d'Analyse, Topologie et Probabilit\'es, UMR 6632, 
Centre de Math\'ematiques et d'Informatique, Universit\'e  Aix-Marseille I,
Technopole de Chateau-Gombert, 
39, rue Fr\'ed\'eric Joliot-Curie -
 13453 Marseille Cedex 13}
\email{derbez@cmi.univ-mrs.fr}



\subjclass{57M50, 51H20}
\keywords{Haken manifold, Seifert fibered space, hyperbolic 3-manifold,  Gromov simplicial volume,  non-zero degree maps, Dehn filling, subgroup separability}

\date{\today}

\begin{abstract}
A natural problem in the theory of 3-manifolds is the question of whether  two 3-manifolds are homeomorphic or not. The aim of this paper is to study this problem for the class of closed Haken manifolds using degree one maps. 

To this purpose we introduce an invariant $
\tau(N)=\left({\rm Vol}(N),\|N\|\right)$ where $\|N\|$ denotes the Gromov simplicial volume of $N$ and ${\rm Vol}(N)$ is a  2-dimensional simplicial volume which measures the  volume of the base 2-orbifolds of the Seifert pieces of $N$.

After studying the behavior of $\tau(N)$ under nonzero degree maps action, we prove  that if $M$ and $N$ are closed Haken manifolds such that $\|M\|=\abs{{\rm deg}(f)}\|N\|$ and  ${\rm Vol}(M)={\rm Vol}(N)$ then any non-zero degree  map $f\co M\to N$ is homotopic to a covering map. This extends a result of S. Wang in \cite{W1} for maps of nonzero degree from $M$ to itself.  
As a corollary we prove that if $M$ and $N$ are closed Haken manifolds such that $\tau(N)$ is sufficiently close to $\tau(M)$ then any  degree one map $f\co M\to N$ is homotopic to a homeomorphism. 
\end{abstract}
\maketitle

\vspace{-.5cm}
\tableofcontents

\section{Introduction}

\subsection{Simplicial volume of a manifold} Let $N^n$ be a  $n$-dimensional manifold. The simplicial volume of $N$ is a homotopy invariant of $N$ defined by M. Gromov in \cite{G} using the $l^1$-pseudo norm on singular homology as follows: for an element $h\in H_{\ast}(N,\b N;{\R})$, the Gromov norm is given by $$\|h\|=\inf\left\{\sum_{i=1}^{i=r}\abs{a_i}, \ {\rm when} \sum_{i=1}^{i=r}a_i\sigma_i\ {\rm  represents}\ \ h\right\}$$ The Gromov simplicial volume of $N$, denoted by $\|N\|$ is the Gromov norm of the image of a generator of $H_{n}(N,\b N;{\Z})$ under the canonical homomorphism $H_{n}(N,\b N;{\Z})\to H_{n}(N,\b N;{\R})\simeq H_{n}(N,\b N;{\Z})\otimes{\R}$.

\subsection{Simplicial volume of a Haken manifold} Let $N$ be a closed Haken manifold. Given a submanifold $K$ of $N$ we denote by $W(K)$ a regular neighborhood of $K$ in $N$.  Denote by ${\mathcal T}_N$ the JSJ-family of $N$, by ${\mathcal S}(N)$, resp. ${\mathcal H}(N)$, the Seifert, resp. hyperbolic, components of $N^{\ast}=N\setminus W\left({\mathcal T}_N\right)$ and by $\Sigma(N)=(\Sigma(N),\emptyset)$ the characteristic Seifert pair of $N$ (see \cite{JS} and \cite{J}). The Cutting off Theorem of M. Gromov (\cite{G}) combined with the fact that manifolds admitting a fixed point free ${\S}^1$-action have zero Gromov simplicial   volume (by the Mapping Theorem of M. Gromov) implies that  
$$\|N\|=\sum_{H\in{\mathcal H}(N)}\|H\|$$
In particular this means that  the Gromov simplicial  volume   of a Haken manifold only depends on its hyperbolic pieces.
In the following it will be convenient to decompose ${\mathcal S}(N)$ into two parts depending on the geometry of the components of ${\mathcal S}(N)$. We denote by ${\mathcal S}_{h}(N)$, resp. by ${\mathcal S}_{e}(N)$,  the components of ${\mathcal S}(N)$ admitting a Seifert fibration with hyperbolic, resp. Euclidean, base 2-orbifold.

\subsection{Extending the Simplicial volume}

  To get a rigidity theorem for Haken manifolds we need to  add an other invariant of $N$ which does not vanish  on ${\mathcal S}(N)$ provided ${\mathcal S}(N)$ is "non-trivial" (i.e. when ${\mathcal S}_{h}(N)\not=\emptyset$). To this purpose we define a kind of 2-dimensional simplicial volume for $N$. More precisely, let $S$ be a component of ${\mathcal S}(N)$. Fix a Seifert fibration for $S$ and denote by ${\mathcal O}_S$ the base 2-orbifold of $S$ with respect to the fixed Seifert fibration.  Then we set ${\rm Vol}(S)=\abs{\chi({\mathcal O}_S)}$ and we define the 2-dimensional volume of $N$ by setting  
$${\rm Vol}(N)=\sum_{S\in{\mathcal S}(N)}{\rm Vol}(S)$$
\begin{lemma}\label{invariance}
If $N$ is a closed Haken manifold, the 2-dimensional volume ${\rm Vol}(N)$, and thus the pair $\tau(N)=({\rm Vol}(N),\|N\|)$, is an invariant of $N$. Moreover $\tau(N)=0$ iff $N$ is a virtual torus bundle.
\end{lemma}
 It will be convenient to use the following convention: we say that $(a,b)\geq (c,d)$ if and only if $a\geq c$ and $b\geq d$, where $(a,b)$ and $(c,d)$ are in ${\R}^2$.

\subsection{Volume and nonzero degree maps}
It follows from the definition of the Gromov simplicial volume that nonzero degree maps "decreases the simplicial volume" in the following sense. Let $f\co M\to N$ be a proper nonzero degree map between orientable $n$-dimensional manifolds. Then $\|M\|\geq\abs{{\rm deg}(f)}\|N\|$.  This inequality does not hold with $\tau(N)$.  In particular the relation ${\rm Vol}(M)\geq\abs{{\rm deg}(f)}{\rm Vol}(N)$  is not true. For instance there exists a degree one map from a Euclidean 3-manifold $M$ onto  $N={\bf S}^3$ or $N={\bf S}^2\times{\bf S}^1$. Then we get ${\rm Vol}(M)=0<{\rm Vol}(N)$ and thus the \emph{decreasing problem for the volume} has to be considered only for aspherical 3-manifolds.  In the case of aspherical 3-manifolds then  consider an orientable hyperbolic surface $F$ and define for each integer $n$ a map $g_n\co M=F\times{\S}^1\to N=F\times{\S}^1$ such that $g_n(x,z)=(x,z^n)$. Then ${\rm deg}(g_n)=n$ and ${\rm Vol}(M)={\rm Vol}(N)=-\chi(F)>0$. However we have the following comparison result:
\begin{theorem}\label{decroissance}
Let $f\co M\to N$ be a nonzero degree map between closed Haken manifolds. If  $\|M\|=\abs{{\rm deg}(f)}\|N\|$ then 
${\rm Vol}(M)\geq{\rm Vol}(N)$.
\end{theorem} 
Note that the condition on the Gromov simplicial volume is necessary in Theorem \ref{decroissance}. Indeed by a construction of  \cite{BW} using nul-homotopic hyperbolic knots, we know that for any aspherical Seifert fibered space $\Sigma$ there always exists a hyperbolic 3-manifold $M$ such that there exists a degree one map $f\co M\to\Sigma$. In this case ${\rm Vol}(M)=0$ and $\Sigma$ can be choosen so that ${\rm Vol}(\Sigma)>0$.
   
In view of Theorem \ref{decroissance}  the  following  question is natural: If  $\|M\|=\abs{{\rm deg}(f)}\|N\|$ then what happens when ${\rm Vol}(M)={\rm Vol}(N)$? This is the aim of the following section.

\subsection{Volume and topological rigidity}
The purpose of this paper is to characterize those degree one (resp. nonzero degree) maps  between closed Haken manifolds which are homotopic to a homeomorphism (resp. covering).
  Then our main result states as follows:
\begin{theorem}\label{rigide} Let $f\co M\to N$ be a nonzero degree map between closed Haken manifolds such that   $\|M\|=\abs{{\rm deg}(f)}\|N\|$.  If ${\rm Vol}(M)={\rm Vol}(N)$ then $f$ is homotopic to a ${\rm deg}(f)$-fold covering.  
\end{theorem}
\begin{remark}
In Theorem \ref{rigide} we can obviously decompose $f$ into two covering maps which preserve the JSJ-decomposition. This means that after a homotopy, $f$ induces two covering maps $f|{\mathcal H}(M)\co{\mathcal H}(M)\to{\mathcal H}(N)$ and $f|{\mathcal S}(M)\co{\mathcal S}(M)\to{\mathcal S}(N)$. Since a Seifert fibered space can be seen as a \emph{generalized} ${\S}^1$-bundle over a 2-dimensional orbifold then it could be convenient to precise the behavior of the covering map $f|{\mathcal S}(M)$ with respect to this anisotropic structure. Actually, when the fibration of a Seifert manifold $S$ is unique (up to isotopy),  the action of $f|S$ can be decomposed into two transversal actions:  a \emph{vertical action} (i.e. action along the ${\S}^1$-fibers of $S$) and a \emph{horizontal action} (i.e. action along the 2-orbifolds of $S$). Then in the proof  of Theorem \ref{rigide} we will see that the hypothesis  ${\rm Vol}(M)={\rm Vol}(N)$ implies that $f|{\mathcal S}_{h}(M)$ acts only vertically and that the horizontal action is trivial.
\end{remark}

\begin{remark}
Note that in \cite{W1}, S. Wang proved  that a proper map of nonzero degree $f\co M\to M$ from a Haken manifold $M$ to itself necessarily induces an injective homomorphism at the fundamental group level. Then Theorem \ref{rigide} gives an extension of this result since when $M=N$ the conditions on the volume are satisfied. 
\end{remark}
If we consider only degree one maps then 
one can weaken the hypothesis concerning the volumes. More precisely, combining Theorem \ref{rigide} and \cite[Theorem 1.2]{D} we get the following
\begin{corollary}\label{strong rigidity}
For any closed Haken manifold $M$ there exists a constant $\eta_M\in (0,1)$ depending only on $M$ such that any 
    degree one map $f\co M\to N$ onto a closed  Haken manifold is  homotopic to a homeomorphism  iff $\tau(N)\geq\tau(M)(1-\eta_M)$.
\end{corollary} 
\subsection{Known results on Topological rigidity}
\subsubsection{Rigidity of Surface bundles}
 
    The above problem has been studied by S. Wang and M. Boileau in \cite{W} and \cite{BW} when the domain $M$ is a surface bundle over the circle and when the target $N$ is irreducible. In particular in \cite{W}, S. Wang proved that if $M$ is a virtual  torus bundle over the circle  then $f$ is homotopic to a covering map. When $M$ is a bundle over ${\S}^1$ with   a fiber of negative Euler characteristic, denote by $\alpha$ the cohomology class corresponding to the fibration of $M$. Then in \cite{BW}, M. Boileau and S. Wang proved that  if there is a rational cohomology class $\beta$ in $N$ with $f^{\ast}(\beta)=\alpha$ and such that $\|\alpha\|_{\rm Th}=\abs{{\rm deg}(f)}\|\beta\|_{\rm Th}$ then $f$  is homotopic to a covering map. Here $\|.\|_{\rm Th}$ denotes the Thurston norm.  
    
    \begin{remark} Notice that the constant ${\rm Vol}(M)$ in Theorem \ref{rigide} plays the same role in the "graph part" of $M$ as the Thurston norm of $\alpha$ in the result of M. Boileau and S. Wang in \cite[Theorem 2.1]{BW}.
    \end{remark}
    
\subsubsection{Rigidity of hyperbolic manifolds} 
 The problem above is completely solved for   hyperbolic manifolds by a result of M. Gromov and W. Thurston which states as follows:
 \begin{theorem}[M. Gromov, W. Thurston]\label{gt}
 Let $M$ and $N$ be two complete finite volume hyperbolic 3-manifolds. Then a non zero degree  map $f\co M\to N$ is homotopic to a ${\deg}(f)$-fold covering iff $\|M\|=\abs{{\deg}(f)}\|N\|$.    
 \end{theorem}
Recall that  T. Soma gave a generalization (see \cite{Smostow}) of this result for degree one maps  by proving the following result:
\begin{theorem}[T. Soma]
 For any $\epsilon>0$ there is a constant $\eta_{\epsilon}>0$ which depends only on $\epsilon$ such that,   any degree one map $f\co M\to N$ between closed hyperbolic 3-manifolds satisfying $\|M\|\leq\epsilon$ and $\|N\|\geq\|M\|(1-\eta_{\epsilon})$  is homotopic to an isometry.    
 \end{theorem}
  Notice that $\lim_{\epsilon\to +\infty} \eta_{\epsilon}=0$.
 
 Note also that this kind of result can not be extended for Haken manifolds even if the target is a closed hyperbolic manifold. This comes from the Thurston hyperbolic surgery theorem.
 
Indeed, let $Y$ be a complete finite volume orientable hyperbolic 3-manifold with $\b Y\simeq S^1\times S^1$ and let $X$ denote an orientable graph manifold with  $\b X\simeq S^1\times S^1$ in such a way that there exists a simple closed curve $l$ in $\b X$ such that the pair $(X,l)$ is \emph{pinchable}. This means that there exists a proper degree one map $\pi\co (X,\b X)\to(V,\b V)$ where $V$ is a solid torus ${\bf D}^2\times{\S}^1$ such that $\pi\co\b X\to\b V$ is a homeomorphism which sends $l$ to the meridian $m=\b{\bf D}^2\times\{\ast\}$ in $\b V$. To perform this operation it is sufficient to choose $X$ so that $l$ is nul-homologous in $\H{X}$ (for instance $X=F\times{\S}^1$ where $F$ is an orientable surface with connected boundary and $l=\b F$).

Let $\{l_n, n\in{\N}\}$ be a sequence of simple closed curves in $\b Y$ such that $\{{\rm lenght}(l_n),n\in{\N}\}$ define a strictly increasing sequence with $\lim_{n\to\infty}{\rm lenght}(l_n)=+\infty$, where ${\rm lenght}$ denotes the lenght for the Euclidean metric on $\b Y$ induced by the hyperbolic metric of ${\rm int}(Y)$.  Denote by $M_n$ the closed Haken manifold obtained by gluing $X$ and $Y$ along $\b X$ and $\b Y$ in such a way that $l$ is identified with $l_n$ and denote by $N_n$ the 3-manifold obtained from $Y$ after performing a Dehn filling along the curve $l_n$. Thus the map $\pi$ can be extended by the identity to construct a  degree one map $f_n\co M_n\to N_n$. Then $\|M_n\|=\|Y\|>0$.   By the Thurston hyperbolic surgery theorem, one sees that the $N_n$'s are closed hyperbolic manifolds for $n$ sufficiently large and $\{\|N_n\|, n\in{\N}\}$ is a strictly increasing sequence such that $\lim_{n\to\infty}\|N_n\|=\|Y\|$. Moreover the maps $f_n$ are neither homotopic to a homeomorhism. 
\subsection{Organization of the paper}
This paper is organized as follows. 

In section 2 we recall some terminology and we state  some results on finite coverings  of  Haken manifolds. This section has essentially a technical interest for our purpose.   

Section 3 is devoted to the study of characteristic maps between closed Haken manifolds  (maps that preserves the Jaco-Shalen-Johannson  decomposition). We first give a result which allow to construct by surgeries a characteristic map from a given nonzero degree map between closed Haken mamifolds. Then we describe the behavior of characteristic  maps. More precisely the  characteristic  maps gives a \emph{thick-thin} decomposition of the domain $M$ (when $f\co M\to N$) and we show that $f$ has a \emph{virtual} standard form with respect to this decomposition. 

In Section 4 we use the results stated in Section 3 to prove Theorem \ref{decroissance} for characteristic maps (see Proposition \ref{fine}). More precisely we use the \emph{standard form} for characteristic maps  to show that it suffices to check Theorem \ref{decroissance} for proper nonzero degree maps $f\co G\to\Sigma$ from a Haken graph manifold $G$ to a Seifert fibered space with geometry ${\Hi}^2\times{\R}$ or $\t{SL}_2({\R})$. If $G$ is a Seifert fibered space then the map $f$ descends to a nonzero degree map from the orbifold ${\mathcal O}_G$ of $G$ to the orbifold ${\mathcal O}_{\Sigma}$ of $\Sigma$ and thus one can compare the volume of $G$ and $\Sigma$. When $G$ is not a Seifert fibered space, we use the decomposition stated in Section 3 which cut $G$ into a thick part $G_{\rm thick}$ and a thin part $G_{\rm thin}$. 

More precisely, the results proved in Section 3 allow to assume that there exists a family of vertical tori ${\mathcal T}_v$ in $\Sigma$ such that (after a homotopy)
$f(G_{\rm thin})={\mathcal T}_v$ and,
$f|G_{\rm thick}$ is a fiber preserving map. 

When the tori of the family ${\mathcal T}_v$ are pairwise disjoint then it can be shown that the contribution of $f|G_{\rm thin}$ is negligible to produce the volume of the image of $f$ and thus we prove that ${\rm Vol}(G_{\rm thick})\geq{\rm Vol}(\Sigma)$. Note that in this case, a crucial point is that when the tori of the family ${\mathcal T}_v$ are pairwise disjoint then the volume of $\Sigma$ does not change after removing ${\mathcal T}_v$. But when the tori of the family ${\mathcal T}_v$ are not pairwise disjoint the above arguments do not hold and the contribution of the map $f|G_{\rm thin}$ may be non-trivial and hence we have to construct a convenient \emph{orbifold complex} to compare the volumes. 

More precisely, for any Seifert piece $S$ in  $G_{\rm thick}$ the map $f|S\co S\to\Sigma$ descends to a map from the base orbifold ${\mathcal O}_S$ of $S$ to ${\mathcal O}_{\Sigma}$ but when $S$ is a Seifert piece of $G_{\rm thin}$ the map $f|S\co S\to\Sigma$ does not factors throught the base 2-orbifolds in general. Hence we decompose of $G_{\rm thin}$ into $G^1_{\rm thin}\cup G^2_{\rm thin}$ where for each Seifert piece $S$ of $G^1_{\rm thin}$, resp. $G^2_{\rm thin}$, the map $f|S$ is fiber preserving, resp. is not fiber preserving. Denote by $G_{\rm main}$ the union of the Seifert pieces of  $G_{\rm thick}$ and of $G^1_{\rm thin}$.  Thus we construct an orbifold complex $\hat{\Gamma}\cup{\mathcal O}_{\rm main}$ obtained by connecting the components of  ${\mathcal O}_{\rm main}=\cup_{S\in G_{\rm main}}{\mathcal O}_S$ by a 1-dimensional graph  $\hat{\Gamma}$. We define a volume for  $\hat{\Gamma}\cup{\mathcal O}_{\rm main}$ and we prove that $\hat{\Gamma}\cup{\mathcal O}_{\rm main}$ dominates ${\mathcal O}_{\Sigma}$ with respect to the volumes. The proof of this domination depends on the geometry of $\Sigma$.

When $\Sigma$ has a ${\Hi}^2\times{\R}$-structure then $\Sigma$ admits a horizontal surface say $F_{\Sigma}$. Note that $F_{\Sigma}$ dominates ${\mathcal O}_{\Sigma}$ since they are related by a branched covering. Then consider the surface ${\mathcal F}=f^{-1}(F_{\Sigma})$ which dominates $F_{\Sigma}$ by construction. Note that ${\mathcal F}$ is naturally divided into a thick part ${\mathcal F}_{\rm thick}$ and a thin part ${\mathcal F}_{\rm thin}$ that is itself divided into the union ${\mathcal F}^1_{\rm thin}\cup{\mathcal F}^2_{\rm thin}$ and the standard form of Section 3 allows to construct a graph $\Gamma$ embedded in ${\mathcal F}$ such that $\Gamma\cup{\mathcal F}_{\rm main}$ dominates $F_{\Sigma}$ where ${\mathcal F}_{\rm main}={\mathcal F}_{\rm thick}\cup{\mathcal F}^1_{\rm thin}$. Then the key point of this case is to prove that the dominating map $f|\Gamma\cup{\mathcal F}_{\rm main}\co\Gamma\cup{\mathcal F}_{\rm main}\to F_{\Sigma}$ descends to  a dominating map $\hat{\Gamma}\cup{\mathcal O}_{\rm main}\to{\mathcal O}_{\Sigma}$.

When $\Sigma$ has a $\t{SL}_2({\R})$-structure then we have no horizontal surface and the problem is more delicate. In this case we consider the space $\Sigma'$ obtained from $\Sigma$ after removing a regular neighborhood of a regular fiber. Then using the standard form for characteristic maps one can construct a graph manifold $G'$ "over" $G$  and a corresponding nonzero degree map $f'\co G'\to\Sigma'$ such that $G'_{\rm thin}=G_{\rm thin}$. Then we can apply the arguments of the first case with a horizontal surface $F_{\Sigma'}$ in $\Sigma'$ and a surface  ${\mathcal F}=(f')^{-1}(F_{\Sigma'})$. Next we prove that $F_{\Sigma'}$ "dominates" in a convenient way ${\mathcal O}_{\Sigma}$ so that the map $\Gamma\cup{\mathcal F}_{\rm main}\to F_{\Sigma'}\to {\mathcal O}_{\Sigma}$ factors throught $\hat{\Gamma}\cup{\mathcal O}_{\rm main}$. 

On other key point is to compare ${\rm Vol}(\hat{\Gamma}\cup{\mathcal O}_{\rm main})$ and ${\rm vol}(G)$. More precisely, the purpose of Section 4 is to  prove that ${\rm Vol}(G)\geq{\rm Vol}(\Sigma)$ and in particular if $G_{\rm thin}\not=\emptyset$ then  ${\rm Vol}(G)>{\rm Vol}(\Sigma)$. Since by the step above we know that ${\rm Vol}(\hat{\Gamma}\cup{\mathcal O}_{\rm main})\geq{\rm Vol}(\Sigma)$ then it is sufficient to check that ${\rm Vol}(\hat{\Gamma}\cup{\mathcal O}_{\rm main})<{\rm Vol}(G)$. This inequality can be directly checked when the genus $g_S$ of the base 2-orbifolds of some Seifert pieces $S$ of $G^2_{\rm thin}$  is sufficiently large. In the general case,  we use some results on subgroup separability of Haken manifolds groups of E. Hamilton, and K. Gruenberg  and a system of horizontal surfaces to construct a finite covering $\t{f}\co\t{G}\to\t{\Sigma}$ of $f$ satisfying  the condition on the genus. This proves that ${\rm Vol}(\t{G})>{\rm Vol}(\t{\Sigma})$. Note that the finite covering has to be chosen carefully so that the above inequality descends to a strict inequality ${\rm Vol}(G)>{\rm Vol}(\Sigma)$. This latter point can be achieved provided we control the difference between the fiber degree of $\t{\Sigma}\to\Sigma$ and that of $\t{G}\setminus{\mathcal T}_{\t{G}}\to G\setminus{\mathcal T}_{G}$. 

 Section 5 is devoted to the proof of Theorems \ref{decroissance} and \ref{rigide} and Corollary \ref{strong rigidity}. Note that in this paper all the 3-manifolds are orientable.

\section{Preliminaries} In this section we first prove the invariance of the 2-dimensional simplicial volume. Next we   recall some well known facts on  Seifert fibered spaces and on Haken manifolds which will be used throughout this paper. These results concern finite coverings and separability in Haken manifolds groups.
Let $\Sigma$ be an orientable Seifert fibered space. Then $\Sigma$ is a ${\S}^1$-bundle over its base 2-orbifold ${\mathcal O}_{\Sigma}$ and the ${\S^1}$-action is globally well-defined since $\Sigma$ is orientable. Recall that if $\o{{\mathcal O}}_{\Sigma}$ denotes the underlying space of ${\mathcal O}_{\Sigma}$ and if $c_1,...,c_r$ denote the exceptional points of  ${\mathcal O}_{\Sigma}$ with index $\mu_1,...,\mu_r$ respectively then 
$$\chi\left({\mathcal O}_{\Sigma}\right)=\chi\left(\o{{\mathcal O}}_{\Sigma}\right)-\sum_{i=1}^{i=r}\left(1-\frac{1}{\mu_i}\right)$$
The geometry of ${\mathcal O}_{\Sigma}$ is hyperbolic, Euclidean or Spherical when $\chi\left({\mathcal O}_{\Sigma}\right)$ is respectively $<0$, $=0$ or $>0$. Hence the geometry of $S$ depends of the geometry of ${\mathcal O}(S)$ combined with the rational Euler number ${\bf e}(S)$ of the fibration. More precisely, when ${\bf e}(S)=0$ then we get respectively  a ${\Hi}^2\times{\R}$, Euclidean, ${\S}^2\times{\R}$-structure and when ${\bf e}(S)\not=0$ we get respectively  a $\t{SL}_2({\R})$, ${\rm Nil}$, Spherical structure. Note that if $N$ is a Sol-manifold then we consider it as a Haken manifold with non-empty JSJ-decomposition so that the  Seifert pieces of $N$ are Euclidean manifolds.
\subsection{Two-dimensional simplicial volume}
In this paragraph we prove Lemma \ref{invariance}. Since the JSJ-decomposition of closed Haken manifolds is unique up to isotopy then we only have to check that the volume ${\rm Vol}(N)$ does not depend of the chosen Seifert fibration on the components of ${\mathcal S}(N)$.  Let $S$ be a Seifert piece of $N$. Since $N$ is a closed Haken manifold then $S$ admits one of the following geometries: ${\Hi}^2\times{\R}$, $\t{SL}_2({\R})$, ${\rm Nil}$ or Euclidean geometry. The only aspherical Seifert fibered spaces which admit more than one non-isotopic Seifert fibrations are Euclidean manifolds. But in this case the Euler characteristic of the base orbifold of $S$ is always zero. Hence the invariance is immediate.

It remains to check the second assertion of the lemma.   Assume that $N$ admits a finite covering $\pi\co\t{N}\to N$ which is a torus bundle over the circle. Then $\t{N}$ is a geometric manifold and the structure depends on the monodromy of the bundle. Then $N$ admits a Euclidean, or a Nil or a Sol geometry. In the case of Euclidean or Nil geometry $N$ is a Seifert fibered space and the base 2-orbifold ${\mathcal O}_N$ is Euclidean and thus $\tau(N)=0$. If $N$ is a Sol-manifold then each component  of $N\setminus{\mathcal T}_N$ is a Euclidean manifold and hence  $\tau(N)=0$.   Assume that $\tau(N)=0$. If ${\mathcal T}_N=\emptyset$ then $N$ has Euclidean or Nil-geometry. In any case $N$ is a virtual torus bundle.  If ${\mathcal T}_N\not=\emptyset$ then ${\mathcal H}(N)=\emptyset$ and each Seifert piece of $N$ is a Euclidean manifold with non-empty boundary. Then by minimality of the JSJ-decomposition either 

(i) $N$ is made of two twisted $I$-bundles over the Klein bottle glued along their boundary or 

(ii) $N$ is ${\S}^1\times{\S}^1\times I/\l\varphi\r$, where $\varphi\co{\S}^1\times{\S}^1\times\{0\}\to{\S}^1\times{\S}^1\times\{1\}$ is an Anosov diffeomorphism.

In case (ii) $N$ is a torus bundle over the circle (actually a Sol-manifold) and in case (i) $N$ admits a 2-fold covering that is a torus bundle over the circle. This completes the proof of Lemma \ref{invariance}.
\subsection{Dehn fillings}
We define \emph{Seifert Dehn fillings}. Suppose $\Sigma$ is an orientable Seifert fibered space  with $\b\Sigma\not=\emptyset$ and let $T$ be a component of $\b\Sigma$. Since $\Sigma$ is orientable then  $T\simeq{\S}^1\times{\S}^1$. Let $\alpha$ be a simple closed curve in $T$. Performing a Dehn filling on $T$ along $\alpha$ means that we glue a solid torus $V={\bf D}^2\times{\S}^1$ identifying $\b{\bf D}^2\times{\S}^1$ with $T$ so that $\alpha$ is glued with the meridian $\b{\bf D}^2\times\{\ast\}$ of $V$. Denote by $\hat{\Sigma}=\Sigma(\alpha)$ the resulting manifold. When $\alpha$ is not isotopic to a generic fiber of $\Sigma$ then  the Seifert fibration of $\Sigma$ extends to a Seifert fibration of  $\hat{\Sigma}$ and we say that we have performed a Seifert Dehn fillings.

\subsection{Morphisms}
 Let $f\co\Sigma\to\Sigma'$ be a map between orientable Seifert fibered spaces. We say that $f$ is a \emph{bundle homomorphism} is there exists a Seifert fibration of $\Sigma$ and $\Sigma'$ so that $f$ is a  homomorphism for the  ${\S}^1$-bundle structures on $\Sigma$ and $\Sigma'$. According to \cite{Ro}, for bundle homomorphisms, we define the following degrees: 

The \emph{fiber degree} of $f$ is the integer $n$ given by $f_{\ast}(h)=t^n$ where $h$, resp. $t$, denotes the generic fiber of $\Sigma$, resp. of $\Sigma'$, and we denote it by $G_h(f)$.

The \emph{orbifold degree} $G_{\rm ob}(f)$ is the minimum number of regular fibers in $g^{-1}(t)$, when $g$ runs over all bundle homomorphisms properly homotopic to $f$ and transverse to $t$.  

For a bundle homomorphism $f\co\Sigma\to\Sigma'$ we have
$$\abs{{\rm deg}(f)}\leq G_h(f)G_{\rm ob}(f)$$
We say that a bundle homomorphism is  \emph{allowable} if $\abs{{\rm deg}(f)}=G_h(f)G_{\rm ob}(f)$ .
In particular,  a bundle homomorphism $f\co(\Sigma,\b\Sigma)\to(\Sigma',\b\Sigma')$ between orientable Seifert fibered spaces with non-empty boundary which is proper (i.e. $f^{-1}(\b\Sigma')=\b\Sigma$)  is \emph{allowable}.

\subsection{Finite coverings of Seifert and Haken manifolds} 
We recall the following result which can be found in \cite{JS}.
\begin{lemma}(\cite[Lemma II.6.1]{JS})\label{jacoshalen}
Let $\Sigma$ be a Seifert fibered space. Then any finite covering $\pi\co\t{\Sigma}\to\Sigma$ admits a Seifert fibration so that $\pi$ is an allowable bundle homomorphism. Moreover the Euler characteristic of the base orbifolds satisfy
$$\chi({\mathcal O}_{\t{\Sigma}})=G_{\rm ob}(f)\chi({\mathcal O}_{\Sigma})$$
\end{lemma} 
 We will use the following result whose proof follows  from the Selberg Lemma, (see \cite{Al} for a proof), that says that a finitely generated matrix group over a field of characteristic zero has a torsion-free subgroup of finite index.
 \begin{lemma}\label{covering}
 Let $N$ be a closed Haken manifold. If $\tau(N)\not=0$ then $N$ admits a finite covering $\t{N}$, which induces the trivial covering over the JSJ-family such that each component of ${\mathcal S}(\t{N})$ (if non-empty) is a circle bundle over an orientable hyperbolic surface.
 \end{lemma}
 \begin{proof}[Sketch of the proof]
 First note that we may assume that $N$ contains no Seifert piece homeomorphic to the twisted $I$-bundle over the Klein bottle ${\bf K}^2$. Indeed if a Seifert piece $S$ of ${\mathcal S}(N)$ is homeomorphic to ${\bf K}^2\t{\times}I$ then $S$ admits a 2-fold covering $\t{S}$ homeomorphic to ${\S}^1\times{\S}^1\times I$ acting trivially on the boundary. Then $N$ admits a 2-fold covering acting trivially on ${\mathcal H}(N)\cup{\mathcal T}_N$  which contains no embedded Klein bottle. 
 
 Suppose first that ${\rm Vol}(N)=0$. Since $\tau(N)\not=(0,0)$ then ${\mathcal H}(N)\not=\emptyset$. Suppose that  ${\mathcal S}(N)\not=\emptyset$. Choose a component  $S\in{\mathcal S}(N)$. Then $\b S\not=\emptyset$ and since ${\rm Vol}(S)=0$ then $S$ is homeomorphic to the twisted $I$-bundle over the Klein bottle. A contradiction. Then in this case ${\mathcal S}(N)=\emptyset$ and there is nothing to prove.
 
 Assume  ${\rm Vol}(N)\not=0$. Thus ${\mathcal S}(N)\not=\emptyset$. Choose a component  $S\in{\mathcal S}(N)$. Since $N$ contains no twisted $I$-bundle over the Klein bottle then ${\rm Vol}(S)\not=0$. Hence, $S$ admits a unique (up to isotopy) Seifert fibration over a hyperbolic 2-orbifold ${\mathcal O}_S$. Then using the Selberg Lemma we know that $S$ admits a finite covering $\t{S}$ inducing the trivial covering over $\b S$ and such that $\t{S}$ is a ${\S}^1$-bundle over an orientable hyperbolic surface. Hence we perform this construction for any component $S$ of ${\mathcal S}(N)$ and we glue together finitely many copies of the covering spaces along the boundary to get a covering $\t{N}$ of $N$ satisfying the conclusion of the lemma. Note that one can do that since the coverings $(\t{S},\b\t{S})\to(S,\b S)$ act trivially on the boundary. 
 \end{proof}
 
 We state the following result which is a very special case of Theorem \ref{decroissance}. It will be used  in Section 4.

\begin{lemma}\label{cover} Let $N$ be a Haken manifold  and let $p\co\t{N}\to N$ be a finite covering of $N$. Then ${\rm Vol}(N)\leq{\rm Vol}(\t{N})\leq\abs{{\rm deg}(p)}{\rm Vol}(N)$. Moreover if $p$ induces the trivial covering over the regular fibers of the components of ${\mathcal S}(N)$ then ${\rm Vol}(\t{N})=\abs{{\rm deg}(p)}{\rm Vol}(N)$.
\end{lemma}
\begin{proof}
First note that if ${\rm Vol}(N)=0$ then  any finite covering $\t{N}$ of $N$ satisfies  ${\rm Vol}(\t{N})=0$. Thus we may assume that ${\rm Vol}(N)\not=0$.

Let $S$ be a component of ${\mathcal S}(N)$ admitting a  ${\bf H}^2\times{\R}$ or a $\t{SL(2,{\R})}$ structure. Notice that in the latter case, necessarily $S=N$. Choose a component $\t{S}$ of $p^{-1}(S)$ in $\t{N}$. Denote by ${\mathcal O}_S$ and by ${\mathcal O}_{\t{S}}$ the base 2-orbifold of $S$ and $\t{S}$. Note that it follows from our hypothesis that ${\mathcal O}_S$ and  ${\mathcal O}_{\t{S}}$ are hyperbolic. Denote by $n$ the integer such that $p_{\ast}(\t{h})=h^n$, where $h$, resp. $\t{h}$, denotes the homotopy class of the generic fiber of $S$ and $\t{S}$ respectively. Then by Lemma \ref{jacoshalen} we know that 
$${\rm deg}(p|\t{S})=n.G_{\rm ob}(p|\t{S}) \ \ \ {\rm and} \ \ \ {\rm Vol}(\t{S})=G_{\rm ob}(p|\t{S}){\rm Vol}(S)\geq{\rm Vol}(S)$$
On the other hand, notice that
$$\abs{{\rm deg}(p)}=\sum_{\t{S}\in p^{-1}(S)}\abs{{\rm deg}(p|\t{S})}$$
   Hence $\abs{{\rm deg}(p)}{\rm Vol}(N)\geq{\rm Vol}(\t{N})\geq{\rm Vol}(N)$. It remains to prove the second part of the lemma. Let $\Sigma$ be a component of ${\mathcal S}(N)$ with ${\Hi}^2\times{\R}$ or $\t{SL(2,{\R})}$-geometry and denote by $\Sigma_1,...,\Sigma_k$ the components of $\t{\Sigma}=p^{-1}(\Sigma)$.  

Denote by $h,h_1,...,h_k$ the homotopy class of the generic fiber of $\Sigma,\Sigma_1,...,\Sigma_k$. Since by hypothesis $p_{\ast}(h_i)=h^{\pm 1}$ for $i=1,...,k$ then each covering $p_i=p|\Sigma_i$ satisfies $\abs{{\rm deg}(p_i)}=G_{\rm ob}(p_i)$ and thus ${\rm Vol}(\Sigma_i)=\abs{{\rm deg}(p_i)}{\rm Vol}(\Sigma)$ for any $i=1,...,k$ and finally, since $\abs{{\rm deg}(p)}=\abs{{\rm deg}(p|\t{\Sigma})}=\abs{{\rm deg}(p|{\Sigma_1})}+...+ \abs{{\rm deg}(p|{\Sigma_k})}$, then  ${\rm Vol}(\t{\Sigma})={\rm Vol}(\Sigma_1)+...+{\rm Vol}(\Sigma_k)=\abs{{\rm deg}(p)}{\rm Vol}(\Sigma)$. This ends the proof of the lemma. 
\end{proof}
The following result is stated only for technical reasons and will be used in Section 4.
\begin{lemma}\label{nonseparating}
Let $N$ denote a connected ${\S}^1$-bundle over an orientable hyperbolic surface $F$ with bundle projection $p\co N\to F$ and let $u$ denote a simple closed surve in $F$ that is homotopically non-trivial. Then there exists a finite abelian covering $\left(\t{N},\t{F},\t{p}\right)$ of $\left(N,F,p\right)$  acting trivially on the fiber of $N$ and on $[u]$ such that each component $\t{u}$ over $u$ 
is of infinite order in $H_1\left(\t{F};{\Z}\right)$.
\end{lemma} 
\begin{proof}
Let $u$ denote a homotopically non-trivial simple closed curve in $F$. If $u$ is not of infinite order in $\H{F}$,  since the group $\H{F}$ is torsion free, then $u$ is a separating curve in $F$. Denote by $A$ and $B$ the components of $F\setminus u$.

Assume first that both $\H{A,u}$ and $\H{B,u}$  are non-zero. Then   one can construct non-trivial finite abelian groups $L_A$, $L_B$ and epimorphisms $\rho_A\co \H{A}\to L_A$ and $\rho_B\co \H{B}\to L_B$ such that $\ker\rho_A\supset\l[u]\r$ and $\ker\rho_B\supset\l[u]\r$. Using the exact sequence $$\H{[u]}\to\H{A}\oplus\H{B}\to\H{F}\to\{0\}$$ we get an epimorphism $\rho\co\H{F}\to L_A\oplus L_B$ such that $\ker(\rho)\supset\l[u]\r$. The finite covering corresponding to the homomorphism given by
$$\pi_1N\stackrel{p_{\ast}}{\to}\pi_1F\to\H{F}\stackrel{\rho}{\to}L_A\oplus L_B$$
 satisfies the required property.

Assume that $\H{A,u}=\{0\}$ (say). This means that $\H{u}\to\H{A}$ is an epimorphism and thus $\H{A}$ is $\{0\}$ or ${\Z}$. In the first case $A$ is a disk which is impossible since $u$ is homotopically non-trivial and in the second case $A$ is an annulus. This means that $u$ is $\b$-parallel in $F$.  Moreover since $u$ is nul-homologous, then $u=\b F$ and in particular $F$ has connected boundary. Since $F$ is hyperbolic then $\H{F}\not=\{0\}$. Then there exists a non-trivial finite group $L$ and an epimorphism $\rho\co\H{F}\to L$. Then the finite covering  corresponding to the homomorphism given by $$\pi_1N\stackrel{p_{\ast}}{\to}\pi_1F\to\H{F}\stackrel{\rho}{\to}L$$
 satisfies the required property. This completes the proof of the lemma.

 \end{proof} 
 \begin{figure}[htb]
\psfrag{A}{$A$}
\psfrag{B}{$B$}
\psfrag{tu}{$\t{u}$}
\psfrag{u}{$u$}
\psfrag{k}{$u$ is boundary parallel}
\psfrag{l}{$u$ is an essential separating curve}
\centerline{\includegraphics{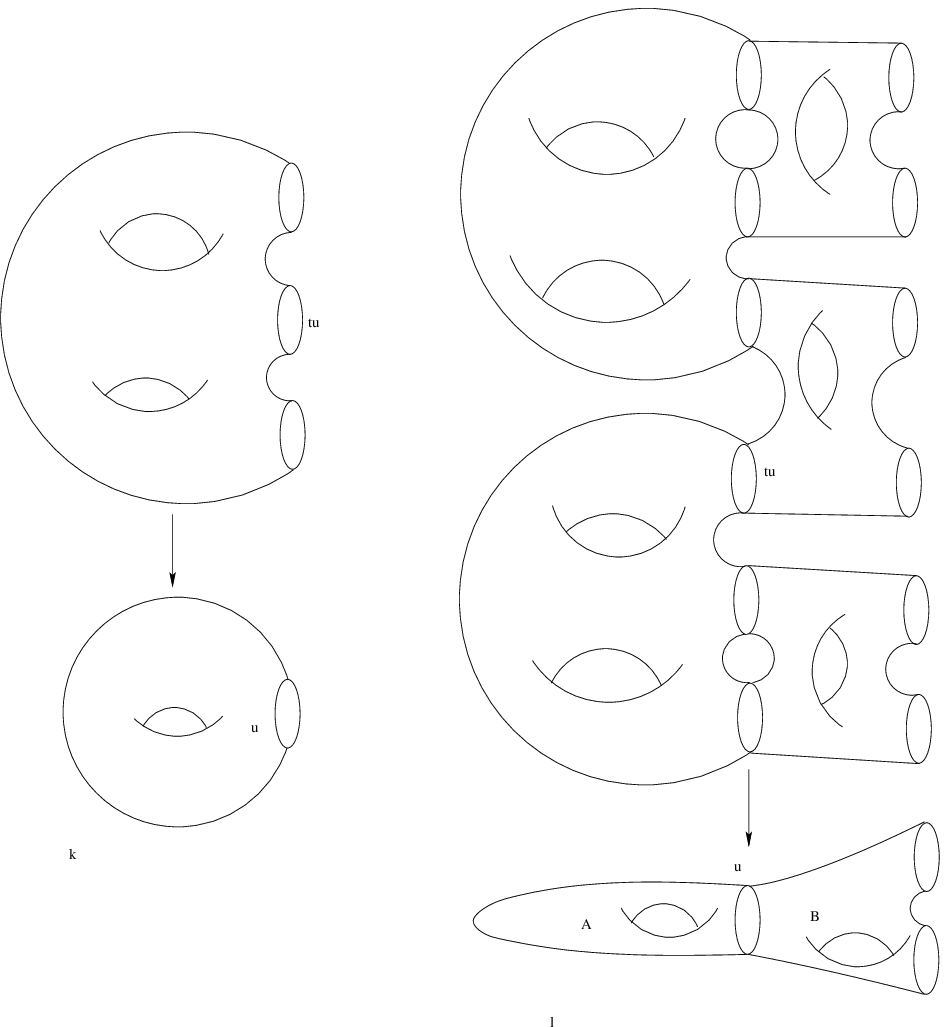}}
\end{figure}

\subsection{Separability of fundamental groups}
Throughout the proof of Theorem \ref{decroissance}, we will need the following technical result which is a direct consequence of an abelian subgroup separability Theorem of E. Hamilton combined with the residual $q$-nilpotence of free groups, for any prime $q$, proved by Gruenberg. 
\begin{lemma}[\cite{H}, \cite{Gr}]\label{residu}
Let $F$ be an orientable hyperbolic surface and let $u\in\pi_1F$ be a non-trivial element. Then for any prime $q$ there exists a finite group $H_q$ and an epimorphism $\tau\co\pi_1F\to H_q$ such that $\tau(u)\not=1$ and $q$ divides the order of $\tau(u)$. 
\end{lemma}
\begin{proof}
Consider $\pi_1F$ as a discret subgroup of  ${\rm Isom}({\Hi}^{2,+})\simeq{\rm PSL}_2({\R})$. 

Assume first that $u$ is a hyperbolic isometry (i.e. $u$ has exactly two fixed points both in $\b_{\infty}{\Hi}^{2,+}$). Then the proof of the lemma follows directly from  Propostion 5 of \cite{H} in this case. Indeed the eigen values of the matrix representing $u$ in ${\rm SL}_2({\C})$ are not root of unity. 

Assume now that $u$ is a parabolic isometry (i.e. $u$ has exactly one fixed point and it lies in $\b_{\infty}{\Hi}^{2,+}$). In this case, necessarily $\b F\not=\emptyset$ and thus $\pi_1F$ is a free group. Then it follows from \cite{Gr} that $\pi_1F$ is residually $q$-nilpotent for any prime $q$. This means that there exists a finite $q$-group $H_q$ and an epimorphism $\tau\co\pi_1F\to H_q$ such that $\tau(u)\not=1$. This completes the proof of the lemma.      
\end{proof}

We end this section with the  following result which will be used  latter in the proof of Theorem \ref{decroissance} and which follows from the residual finiteness of surface groups. 
\begin{lemma}\label{nonabgenus}
Let $f\co S\to\Sigma$ be a map from a Seifert fibered space to a ${\S}^1$-bundle over an orientable hyperbolic surface $F$ such that $f_{\ast}(\pi_1S)$ is non-abelian and  $f|\b S$ is $\pi_1$-injective (i.e.  $\pi_1$-injective on each component of $\b S$). Then for any $n\in{\N}$ there exists a finite covering $\t{f}_n\co\t{S}_n\to\t{\Sigma}_n$ satisfying the following properties:

(i) the covering $\t{\Sigma}_n\to\Sigma$ has trivial fiber degree,

(ii) each component of $\t{S}_n$ has a base 2-orbifold of genus at least $n$.  
\end{lemma}
\begin{proof} Let $T_1,...,T_p$ be the components of $\b S$.
Denote by $t$ the homotopy class of the fiber of $\Sigma$ and by $\pi\co\Sigma\to F$ the bundle projection. Denote by $d_1,...,d_{p}$ the chosen sections of $\b S$ with respect to the fixed Seifert fibration of $S$ and let $c_1,...,c_r$ denote the homotopy class of the exceptional fibers of $S$ with index $\mu_1,...,\mu_r$ respectively.

 Since $f_{\ast}(\pi_1S)$ is non-abelian it follows from \cite{JS} that for any fiber $v$ of $S$, then $f_{\ast}(v)\in\l t\r$.  Denote by ${\mathcal O}_S$ the base 2-orbifold of $S$ and by $\o{\mathcal O}_S$ the underlying space and set $g_S={\rm genus}(\o{\mathcal O}_S)$. In order to prove the lemma we will construct regular coverings of $\Sigma$ acting trivially on $t$ and which induces via $f$ some coverings of $S$ satisfying condition (ii). Consider a regular covering $(q,\t{\Sigma})$ of $\Sigma$. Consider the corresponding epimorphism $\varphi\co\pi_1\Sigma\to K$  where $K$ is a finite group. A component $(p,\t{S})$ of the induced covering over $S$ correponds to the homomorphism $\varphi\circ f_{\ast}$. This covering induces a branched covering of degree $\sigma$ between the underlying spaces of the  base 2-orbifolds of $\t{S}$ and  $S$.    Let $\beta_j$ denote the order of $\varphi f_{\ast}(c_j)$ and for each $i=1,...,p$  denote by $r_i$ the number of components of $\b\t{S}$ over $T_i$ and set $n_i=\sigma/r_i$. Then the Riemann-Hurwitz formula allows to compute the genus of ${\mathcal O}_{\t{S}}$ using the datas of ${\mathcal O}_S$ and those of $p\co\t{S}\to S$:
 $$2g_{\t{S}}=2+\sigma\left(p+2g_S+r-2-\sum_{i=1}^{i=p}\frac{1}{n_i}-\sum_{i=1}^{i=r}\frac{1}{(\mu_i,\beta_i)}\right)$$

\emph{Case 1}: Assume $g_S\geq 2$.  First note that since $f_{\ast}(\pi_1S)$ is non-abelian there exists an element  $a\in\pi_1 S$ such that $\pi_{\ast}\circ f_{\ast}(a)\not=1$ in $\pi_1F$. Since surface groups are residually finite then there exists a finite group $K$ an epimorphism $\epsilon\co\pi_1F\to K$ such that $\epsilon(\pi_{\ast}\circ f_{\ast}(a))\not=1$. Consider the homomorphism $\varphi=\epsilon\circ\pi_{\ast}$. Note that since the regular fiber of $S$ is sent via $f$ to the fiber of $\Sigma$ then necessarily $\sigma\geq 2$. Then the Riemann-Hurwitz formula gives
$$2g_{\t{S}}\geq 2+\sigma\left(2g_{S}-2\right)$$
Thus since $g_S\geq 2$ and $\sigma\geq 2$ we get $2+\sigma\left(2g_{S}-2\right)>2g_S$. This proves that $g_{\t{S}}>g_S$ and completes the proof of the lemma in this case.

\emph{Case 2}: Assume $g_S=1$. Then we claim that $p\geq 1$. Suppose the contrary. Denote by $a,b$ the standard generators of $\pi_1\o{\mathcal O}_S$, by $q_1,...,q_r$ the sections corresponding to the exeptional fibers $c_1,...,c_r$.  Since $f_{\ast}(h)\in\l t\r$ and since $\pi_1F$ is torsion free then $f_{\ast}(q_i)\in\l t\r$ for $i=1,...,r$. Hence, since $[a,b]q_1...q_r=h^b$ then $f_{\ast}([a,b])\in\l t\r$ and thus $[\pi_{\ast}f_{\ast}(a),\pi_{\ast}f_{\ast}(b)]=1$ in $\pi_1F$. Since $F$ is a hyperbolic surface then there exists $u\in\pi_1F$ such that $\l u\r=\l\pi_{\ast}f_{\ast}(a),\pi_{\ast}f_{\ast}(b)\r$. Let $g\in\pi_1\Sigma$ such that $\pi_{\ast}(g)=u$. Then $f_{\ast}(\pi_1S)\subset\l g,t\r\simeq{\Z} \times{\Z}$. A contradiction. Since $f|\b S\co\b S\to\Sigma$ is $\pi_1$-injective then $\pi_{\ast}f_{\ast}(d_i)\not=1$ in $\pi_1F$. Thus there exists an epimorphism into a finite group $K$ denoted by $\epsilon\co\pi_1F\to K$ such that $\epsilon\pi_{\ast}f_{\ast}(d_i)\not=1$ for $i=1,...,p$. Consider the homomorphism $\varphi=\epsilon\pi_{\ast}$ and the associated coverings $\t{S}$ and $\t{\Sigma}$.  Then it follows from our construction that $n_i\geq 2$ for $i=1,...,p$ and $\sigma\geq 2$. Then the Riemann-Hurwitz formula gives
$$2g_{\t{S}}\geq 2+\sigma\frac{p}{2}>2$$
Thus  $g_{\t{S}}\geq 2$ and we have a reduction to the first case. 

\emph{Case 3}: Assume  $g_S=0$. In this case the fundamental group of $S$ admits a presentation
$$\left\l d_1,...,d_p,q_1,...,q_r,h: [h,d_i]=[h,q_j]=1, q_i^{\mu_i}=h^{\gamma_i}, d_1...d_pq_1...q_r=h^b\right\r$$
where  $q_i$ denotes the chosen section correponding to the exptional fiber $c_i$. Note that when $p>0$, i.e. when $\b S\not=\emptyset$, then one can choose $b=0$.    Since $f_{\ast}(h)\in\l t\r$ and since $\pi_1F$ is torsion free then $f_{\ast}(q_i)$ and $f_{\ast}(c_i)$  are in $\l t\r$. Then we first check (using the presentation above and the fact that $f_{\ast}(\pi_1S)$ is non-abelian) that $p\geq 3$. If $p\leq 1$ then we get $f_{\ast}(\pi_1S)\subset\l t\r$. A contradiction. Assume that $p=2$. Using the presentation of $\pi_1S$ we get $\pi_{\ast}f_{\ast}(d_1).\pi_{\ast}f_{\ast}(d_2)=1$. Then this proves that $f_{\ast}(\pi_1S)\subset\l f_{\ast}(d_1), t\r\simeq{\Z}\times{\Z}$. A contradiction. From now one may we assume that $p\geq 3$. 

Note that since $f|\b S\co\b S\to\Sigma$ is $\pi_1$-injective then $\pi_{\ast}f_{\ast}(d_i)\not=1$ in $\pi_1F$. Thus there exists an epimorphism into a finite group $K$ denoted by $\epsilon\co\pi_1F\to K$ such that $\epsilon\pi_{\ast}f_{\ast}(d_i)\not=1$ for $i=1,...,p$. Consider the homomorphism $\varphi=\epsilon\pi_{\ast}$ and the associated coverings $\t{S}$ and $\t{\Sigma}$.  Then it follows from our construction that $n_i\geq 2$ for $i=1,...,p$ and $\sigma\geq 2$. Then the Riemann-Hurwitz formula gives
$$2g_{\t{S}}\geq 2+\sigma\left(\frac{p}{2}-2\right)$$

\emph{Subcase 1}: Assume $g_S=0$ and  $p\geq 4$. This implies that $g_{\t{S}}\geq 1$ and we have a reduction to the second case.

\emph{Subcase 2}: Assume $g_S=0$ and $p=3$. If  the number of connected components $\t{p}$ of  $\t{S}$ is $\geq 4$ then we have a reduction to the subcase 1. Hence assume that   $\t{p}=3$. The Riemann-Hurwitz formula gives
$$2g_{\t{S}}=2-\t{p}+\sigma\left(p+r-2-\sum_{i=1}^{i=r}\frac{1}{(\mu_i,\beta_i)}\right)\geq\sigma-1\geq 1$$
Then we get $g_{\t{S}}\geq 1$. This completes the proof of the lemma. 
\end{proof}
\section{A thick-thin decomposition for Haken manifolds with respect to nonzero degree maps}
This section is devoted to the study of a class of nonzero degree maps between closed Haken manifolds which preserve the geometric decomposition. These maps will be termed \emph{characteristic maps}. We first state two results which give sufficient conditions allowing to construct characteristic maps from a given nonzero degree map using surgeries. Next we describe the standard form of  \emph{characteristic maps}.  

\subsection{Obtaining characteristic maps} 
We first define \emph{characteristic maps}.
\begin{definition}
Let $f\co M\to N$ be a map between closed Haken manifolds. We say that $f$ is characteristic if $f({\mathcal H}(M))\subset{\rm int}{\mathcal H}(N)$, $f({\mathcal S}(M))\subset{\rm int}\Sigma(N)$. 
\end{definition} 
Note that when $N$ is not a virtual torus bundle then using arguments in \cite{Wa} one can show, that if $f\co M\to N$ is a characteristic map then for any component $T\in{\mathcal T}_N$ the space $f^{-1}(T)$ is the disjoint union of components of ${\mathcal T}_M$ possibly with some graph submanifolds of $M$.  
Next we define the \emph{Volume} and \emph{Non-degeneration} conditions.
\begin{definition}
We say that a map $f\co M\to N$ between closed Haken manifolds satisfies the volume, resp. Non-degeneration, condition if $\|M\|={\rm deg}(f)\|N\|$, resp. if $f|{\mathcal T}_{M}\co {\mathcal T}_{M}\to N$ is $\pi_1$-injective (i.e. the map is $\pi_1$-injective on each component of ${\mathcal T}_{M}$).
\end{definition} 

We state the following result which explain how a given nonzero degree map between closed Haken manifolds satisfying the Volume condition induces after surgeries a characteristic map.  
 
 \begin{lemma}\label{ducon}
 Let $f\co M\to N$ be a nonzero degree map between closed Haken manifolds with $\|M\|={\rm deg}(f)\|N\|$. Assume that $N$ is not a virtual torus bundle.   Then  there exists  a canonical submanifold $G\subset M$ which contains ${\mathcal H}(M)$ in its interior (in particular if $\b Q\cap\b G\not=\emptyset$ for $Q\in N^{\ast}$ then $Q$ is Seifert), a closed Haken manifold $\hat{M}$ obtained from $G$ after Seifert Dehn fillings along $\b G$ and an extension $\hat{f}\co\hat{M}\to N$ of $f|G\co G\to N$ with the same degree as $f$ satisfying the Volume and Non-degeneration conditions and which is homotopic to a characteristic map. 
 \end{lemma} 
\begin{proof}
The  proof of  Lemma \ref{ducon} follows from the arguments used in   \cite[Lemma 5.1]{D} without any essential change.
\end{proof}
 When a map $f\co M\to N$ is characteristic we always assume that it satisfies the following \emph{minimality condition}: over all characteristic maps homotopic to $f$ we choose a representant in such a way that the number of connected components of $G_{\Sigma}=f^{-1}(\Sigma)$ is minimal when $\Sigma\in{\mathcal S}(N)$.
As a consequence of Lemma \ref{ducon} we get the following 
\begin{corollary}\label{chienne}
Let $f\co M\to N$ be a nonzero degree maps between closed Haken manifolds. Assume that $N$ is not a virtual torus bundle. Then if $f$ satisfies  the volume and non-degeneration conditions, it is homotopic to a  characteristic map.
\end{corollary}
\subsection{Standard form for characteristic maps}
 For a characteristic map  we define a \emph{thin-thick decomposition} of the domain $M$. More precisely we set
$$M_{\rm{thick}}=\{S\in M^{\ast},\ {\rm such\ that}\
f_{\ast}(\pi_1S)\ {\rm is\ non-abelian}\}$$
$$M_{\rm{thin}}=\o{M\setminus M_{\rm{thick}}}$$

    First of all we give a convenient characterization of the components of $M_{{\rm thin}}$. This result will be used to show that the thin part of $M$ can be sent into a family of virtual tori in the target. This point will be crucial for the study of  non-zero degree maps $f\co M\to N$ satisfying the Volume and Non-degeneration conditions.  

\begin{lemma}\label{zz}
Let $f\co M\to N$ be a
map between Haken manifolds  such that $\tau(N)\not=0$  and satisfying the Non-degeneration condition.
 Then  $S\in M_{{\rm thin}}^{\ast}$ iff
$f_{\ast}(\pi_1S)\simeq{\Z}\times{\Z}$.
\end{lemma}
\begin{proof}
Assume that  $S$ is a geometric piece of $M$ such that
$f_{\ast}(\pi_1S)$ is an abelian group. Since $\pi_1N$ is torsion
free  then $f_{\ast}(\pi_1S)\simeq{\Z}^r$. Since $f|\b S$ is a
non-degenerate map, then $r\geq 2$ and since $N$ is a
three-dimensional manifold, then $r\leq 3$ since the subgroup
$f_{\ast}(\pi_1S)$ must have cohomological dimension at most 3.
Since $N$ is not a virtual torus bundle over the circle, by the
condition $\tau(N)\not=0$, then the fundamental group of $N$ can not contain a group
isomorphic  to ${\Z}\times{\Z}\times{\Z}$. Then necessarily $f_{\ast}(\pi_1S)\simeq{\Z}\times{\Z}$. This completes the proof of the lemma. 
\end{proof}
The following result gives the behavior of the thick-thin decomposition with respect to  finite coverings.
\begin{lemma}\label{abnab} Let $f\co M\to N$ be a nonzero degree
map between closed Haken manifolds such that $\tau(N)\not=0$ and satisfying the Volume and Non-degeneration conditions. Assume that we have homotoped $f$ to a characteristic map.
 Let  $\Sigma$ be a Seifert piece of $N$ that is not Euclidean and let $G_{\Sigma}$ denote $f^{-1}(\Sigma)$.
 Then the following properties hold:

(i) Let $p\co\t{G}_{\Sigma}\to G_{\Sigma}$ denote a finite
covering induced by  $f|G_{\Sigma}$ from some finite covering of
$\Sigma$. Then
$$p^{-1}((G_{\Sigma})_{\rm thin})\subset(\t{G}_{\Sigma})_{\rm thin}, \
\ (\t{G}_{\Sigma})_{\rm thick}\subset
p^{-1}((G_{\Sigma})_{\rm thick})$$ The above inclusions are
equalities when   $\Sigma$ is a circle bundle.

(ii) If $S$ is a  component of $(G_{\Sigma})_{\rm thick}$ then
 $f|S\co S\to\Sigma$ is homotopic to a bundle homomorphism,
\end{lemma}
\begin{proof}[Proof of Lemma \ref{abnab}]

 We first prove point (i). The inclusions
$p^{-1}((G_{\Sigma})_{\rm thin})\subset(\t{G}_{\Sigma})_{\rm thin}$ and $(\t{G}_{\Sigma})_{\rm thick}\subset
p^{-1}((G_{\Sigma})_{\rm thick})$
are obvious. However the
 equality is not true in general since, for instance,  there exist non-abelian groups (whose center is infinite cyclic)
  which contain finite index subgroups  isomorphic to
${\Z}\times{\Z}$.

To complete the proof of point (i) we consider the case where $\Sigma$ is  a
${\S}^1$-bundle over an orientable hyperbolic surface $F$ with bundle projection $\pi\co\Sigma\to F$ and fiber $t$. Let $S$ be a geometric component of
$({G}_{\Sigma})_{\rm thick}$.   We write the short exact sequence of the fibration
$$1\to{\Z}\stackrel{i_{\ast}}{\to}\pi_1\Sigma\stackrel{\pi_{\ast}}{\to}\pi_1F\to 1$$

We will use here a fundamental result which says that any torsion
free group which contains a finite index free  subgroup is free.
Denote by $G$ the non-abelian group  equal to $f_{\ast}(\pi_1S)\subset\pi_1\Sigma$.

Notice that $G\cap i_{\ast}({\Z})$ is non-trivial. Indeed, $\Sigma$ is a circle
bundle over an orientable hyperbolic surface then $\Sigma$
contains no embedded Klein bottles and since the centralizer of
$f_{\ast}(t_S)$, where $t_S$ denotes the regular fiber of $S$,
contains $G$ which is non abelian then by \cite[Addendum to
Theorem VI.I.6]{JS}, $f_{\ast}(t_S)$ is conjugate to a power of
the fiber of $\Sigma$. Moreover, since $f|S$ is a non-degenerate map
then this power is non trivial and thus  $G\cap i_{\ast}({\Z})\not=\{1\}$.

 On the other hand, since $G$ is
non-cyclic then $\pi_{\ast}(G)\not=\{1\}$. Let $H$ be a finite index subgroup
of $G$. If $H$ is abelian, then since  $G$ is a torsion free  and
non-cyclic subgroup of $\pi_1\Sigma$ then $H$ is necessarily  isomorphic to
${\Z}\times{\Z}$. Then  $\pi_{\ast}(H)$ is an infinite cyclic subgroup
(indeed  since $F$ is an orientable hyperbolic surface, it can not
contain a subgroup   isomorphic to ${\Z}\times{\Z}$). Since $\pi_{\ast}$ is
an epimorphism and since $H$ is a finite index subgroup of  $G$
then $\pi_{\ast}(H)$ is a finite index subgroup of  $\pi_{\ast}(G)$. Since $\pi_{\ast}(G)$ is
torsion free, this implies that $\pi_{\ast}(G)$ is itself infinite cyclic
which shows that $G$ is abelian too. This gives a contradiction. Thus
$H$ is non-abelian.

We now prove point (ii). Denote by $x$ the fiber of $\Sigma$.   Since $f_{\ast}(\pi_1S)$ is non-Abelian,
then for any fiber  $c$ of $S$, there exists an integer $n\not=0$  such that $f_{\ast}([c])=t^n$. Using the round handle decomposition we write  $S$ as the union
of three subspaces $H_0\cup H_1\cup H_2$, where $H_0$ is a regular
neighborhood of a regular fiber of $S$ union the exceptional fibers
of $S$, $H_1$ is a regular neighborhood of some vertical annuli in
$S$ whose boundaries are in $\b H_0\cup\b S$ and $H_2$ is a regular
neighborhood  of a regular fiber of $S$. One can deform  $f$ by a
homotopy on each annulus and solid torus of the decomposition of
$S=H_0\cup H_1\cup H_2$ in such a way that  $f$ is a bundle
homomorphism, using Lemmas 2.2 and 2.3 of \cite{Ro}.

\end{proof}

\begin{lemma}\label{tores}
Let $f\co G\to\Sigma$ be a proper nonzero degree map from a Haken graph manifold to an orientable ${\S}^1$-bundle over an orientable hyperbolic surface $F_{\Sigma}$. Suppose that $f$ satisfies the Non-degeneration condition.  There exists a finite covering $\t{f}\co\t{G}\to\t{\Sigma}$ of $f$ and a finite family of vertical tori ${\mathcal T}_v$ (not necessarily pairwise disjoint) in  $\t{\Sigma}$ satisfying the following properties:

(i) $\t{\Sigma}\to\Sigma$ induces the trivial covering over the fibers,

(ii) after a homotopy, $\t{f}(\t{G}_{\rm thin})={\mathcal T}_v$ and in particular for each component $W$ of $\t{G}_{\rm thin}$, $\t{f}(W)$ is contained in a single torus of ${\mathcal T}_v$.
\end{lemma}
\begin{proof}
The first step in the proof of the lemma is to construct a finite covering $\t{f}\co\t{G}\to\t{\Sigma}$, satisfying condition (i) and a finite family of vertical tori ${\mathcal T}_v$ in  $\t{\Sigma}$ such that, after a homotopy, $\t{f}(\t{G}^{\ast}_{\rm thin})={\mathcal T}_v$, in such a way that for each geometric component $\t{S}$ of $\t{G}_{\rm thin}$ then $\t{f}(\t{S})$ is contained in a single torus in ${\mathcal T}_v$.
Let $S$ be a component of $G^{\ast}_{\rm thin}$. We write the exact sequence of the ${\S}^1$-fibration of $\Sigma$ over $F_{\Sigma}$
$$({\mathcal F})\ \ \ \ \ \ \ \ \ \{1\}\to\pi_1({\S}^1)\stackrel{i_{\ast}}{\to}\pi_1\Sigma\stackrel{\pi_{\ast}}{\to}\pi_1F_{\Sigma}\to\{1\}$$
We denote by $t$ a chosen generator of $\pi_1{\S}^1$. Denote by $G$ the subgroup of $\pi_1\Sigma$ equal to $f_{\ast}(\pi_1S)$. Since $G\simeq{\Z}\times{\Z}$ and since $F_{\Sigma}$ is an orientable hyperbolic surface then there exists $a\in{\Z}^{\ast}$ and $b\in\pi_1\Sigma\setminus\ker(\pi_{\ast})$ such that $G=\l i_{\ast}(t^a), b\r$. Since by \cite{Sc}, subgroups of surface groups are almost geometric, then there exists a finite covering $\t{F}_{\Sigma}$ of $F_{\Sigma}$ such that $\pi_1\t{F}_{\Sigma}$ contains $\pi_{\ast}(b)$ and $\pi_{\ast}(b)$ is geometric in $\t{F}_{\Sigma}$.  Using the exact sequence $({\mathcal F})$, the covering $\t{F}_{\Sigma}$ induces a finite covering $\t{\Sigma}_S\to\Sigma$ which is trivial over the fibers of $\Sigma$ and satisfying the following property: $\pi_1\t{\Sigma}_S$ contains  the group $\l i_{\ast}(t), b\r$ and $\l i_{\ast}(t), b\r$ is realized by a vertical torus in $\t{\Sigma}_S$. Consider the covering $\t{f}_S\co\t{G}_S\to\t{\Sigma}_S$ of $f$ corresponding to $\t{\Sigma}_S\to\Sigma$. It follows from our construction that for each component $\t{S}$  over $S$ in $\t{G}_S$ there exists a vertical torus in $\t{\Sigma}_S$ that contains $\t{f}_S(\t{S})$. Consider the covering $\t{\Sigma}$ obtained as the fiber product of the coverings $\{\t{\Sigma}_S, S\in G^{\ast}_{\rm thin}\}$. Denote by $\t{f}\co\t{G}\to\t{\Sigma}$ the covering of $f$ corresponding to $\t{\Sigma}\to\Sigma$. Since, by point (i) of Lemma \ref{abnab}, $\t{G}_{\rm thin}=p^{-1}(G_{\rm thin})=\cup_{S\in G^{\ast}_{\rm thin}} p^{-1}(S)$, this completes the proof of the first step. 

From now on, we may assume that $\Sigma$ contains a family ${\mathcal T}_v$ of vertical tori (not necessarily pairwise disjoint) such that for each component $S$ of $G^{\ast}_{\rm thin}$, then $f(S)$ is contained in a single torus in ${\mathcal T}_v$.   It remains to prove that the same property remains true by replacing $S$ by a connected component $V$ of $G_{\rm thin}$ which contains $S$. We argue by induction on the complexity of the dual graph $\Gamma_V$ of $V$. 

Suppose first that $\Gamma_V$ is a tree. Fix a
component $S_0$ in $V^{\ast}$. We may assume,
possibly after passing to a finite covering, that there exists a
maximal vertical torus ${\bf T}$ in $\Sigma$ such that $f(S_0)={\bf
T}$ and $\pi_1{\bf T}=\l t, b\r$ where $b\in\pi_1\Sigma\setminus\ker
\pi_{\ast}$. Let $S$ be an other  component of $V^{\ast}$ adjacent to $S_0$
along a torus $T_0$. Fix a base point $x\in T_0$ and $y=f(x)\in{\bf
T}$. We set $K_0=f_{\ast}(\pi_1(T_0,x))\subset\pi_1({\bf T}, y)$ and
$H=f_{\ast}(\pi_1(S,x))$. It follows from our construction that $K_0$
and $H$ are both isomorphic to ${\Z}\times{\Z}$ and that $H\supset
K_0$. Thus $K_0$ is a finite index subgroup of $H$. Thus since
$K_0\subset\pi_1({\bf T}, y)$ for any $g\in H$ then there exists an
integer $n_g\in{\Z}$ such that $g^{n_g}\in\pi_1({\bf T}, y)=\l
t,b\r$. On the other hand there exists an integer $\beta\not=0$ and an
element $\alpha\in\pi_1\Sigma\setminus\ker\pi_{\ast}$ such that $H=\l
t^{\beta}, \alpha\r$.

Then, in particular, there exist two nonzero  integers $n,m\in{\Z}$
such that $\pi_{\ast}^n(\alpha)=\pi_{\ast}^m(b)$. It is easy to check that $\l
\pi_{\ast}(\alpha), \pi_{\ast}(b)\r$ is an infinite cyclic subgroup of $\pi_1F_{\Sigma}$.
Indeed, since $F_{\Sigma}$ is hyperbolic then we consider the elements of
$\pi_1F_{\Sigma}$ as isometries of ${\bf H}^{2,+}$. Since $\pi_{\ast}(\alpha)$ and
$\pi_{\ast}^n(\alpha)$ (resp. $\pi_{\ast}(b)$ and $\pi_{\ast}^m(b)$) commute then they are
isometries of the same type (parabolic or hyperbolic) with the same
fixed points. Thus $\pi_{\ast}(b)$ and $\pi_{\ast}(\alpha)$ have the same type with
the same fixed points. Thus the discretness of the group  $\l
\pi_{\ast}(\alpha), \pi_{\ast}(b)\r$ combined with the classification of the
isometries of ${\bf H}^{2,+}$ implies that the latter group is cyclic.
Thus since ${\bf T}$ is a maximal torus then $\pi_{\ast}(\alpha)\in\l \pi_{\ast}(b)\r$
and thus $H\subset\pi_1{\bf T}$. Hence, after a homotopy we may
assume that $f(S\cup_{T_0}S_0)\subset{\bf T}$. This completes the proof of the lemma when $\Gamma_V$ is a tree by repeating this process.

If $\Gamma_V$ is not a tree then ${\rm Rk}(H_1(\Gamma_V;{\R}))\geq
1$. Choose a characteristic non-separating torus $T$ in $V$ and consider the space $\hat{V}$ obtained by cutting $V$ along $T$. Then ${\rm Rk}(H_1(\Gamma_{\hat V};{\R}))<{\rm Rk}(H_1(\Gamma_V;{\R}))$.
Denote  by
$U_1, U_2$ the components of $\b\hat{V}$ over $T$. Consider the
map $f_1=f|\hat{V}\co\hat{V}\to\Sigma$. We know from the induction hypothesis that there
exists a vertical torus ${\bf T}$ in $\Sigma$ such that after
modifying $f_1$ by a homotopy then $f_1(\hat{V})={\bf T}$. Thus deforming slightly the map $f|\hat{V}$ on a regular  neighborhood of ${\bf T}$ identified with ${\bf
T}\times[-1,1]$ we get the following diagram
$$\xymatrix{
(\hat{V},U_1,U_2) \ar[r]^{f_1} \ar[d]_{\varphi}  & \Sigma \\
({\bf T}\times[-1,1], {\bf T}\times\{-1\}, {\bf T}\times\{1\})
\ar[ur]_{\psi}}$$
where $\psi\co{\bf T}\times[-1,1]\to\Sigma$ is an embedding. Thus the map $f|V\co V\to\Sigma$ factors
through a torus bundle denoted $N_{\Phi}$, where $\Phi\in{\rm Diff}({\bf T})$ denotes the monodromy of the fibration. Then we get the
following commutative diagram
$$\xymatrix{
{V} \ar[r]^{f|V} \ar[d]_{\pi}  & \Sigma \\
{N}_{\Phi} \ar[ur]_{f_{\Phi}}}$$

We claim that $f_{\Phi}\co {N}_{\Phi}\to\Sigma$ cannot be
$\pi_1$-injective. Indeed if $f_{\Phi}$ is $\pi_1$-injective then,
since by \cite{Sc} any finitely generated subgroup of $\pi_1\Sigma$
is separable  we may assume, passing to a
finite covering, that $\Sigma$ contains a torus bundle over the
circle   whose fiber is a vertical torus in $\Sigma$. This gives a
contradiction since $\Sigma$ is a circle bundle over an orientable
hyperbolic surface. Hence there exists a non-trivial element
$g\in\ker(f_{\Phi})_{\ast}$. Denote still by ${\bf T}$ the fiber of
$N_{\Phi}$. Recall that $\pi_1N_{\Phi}$ admits a presentation
$$\l \alpha\in\pi_1{\bf T}, h: \Phi_{\ast}(\alpha)=h\alpha h^{-1}\r$$
 where $h$ is represented by a simple
closed curve meeting each fiber exactly one time.
  Then $g$ admits a unique decomposition $g=\alpha h^n$
where $\alpha\in\pi_1{\bf T}$ and $n\not=0$ since $f_{\Phi}|{\bf T}\times\{p\}=\psi|{\bf T}\times\{p\}$ which is an embedding. Since $(f_{\Phi})_{\ast}(\alpha)=(f_{\Phi})^{-n}_{\ast}(h)$ 
and since $(f_{\Phi})_{\ast}(\alpha)\in\pi_1{\bf T}$ then $(f_{\Phi})^{-n}_{\ast}(h)\in\pi_1{\bf T}$ and thus since 
$\Sigma$ is a circle bundle over an orientable hyperbolic surface and since  ${\bf T}$ is maximal then  $(f_{\Phi})_{\ast}(h)\in\pi_1{\bf T}$. Hence we have showed   that $f_{\ast}(\pi_1V)\subset\pi_1\bf T$. This completes the proof of the lemma. 
\end{proof}
Using the same arguments as above we get the following result which is a special case of the lemma above.
\begin{corollary}\label{bouts}
Let $f\co G\to\Sigma$ be a proper nonzero degree map from a Haken graph manifold to an orientable ${\S}^1$-bundle over an orientable hyperbolic surface $F_{\Sigma}$. Suppose that $f$ satisfies the Non-degeneration condition. Assume that there exists a component   $G_{\b}$  of $G_{\rm thin}$ such that $\b G_{\b}\cap\b G$ is non-empty. Then there exists a component ${\bf T}$ of $\b\Sigma$ such that, after a homotopy,   $f(G_{\b})={\bf T}$.
\end{corollary}
We end this section with the following result
\begin{lemma}\label{adjacent}
Let $f\co M\to N$ be a nonzero degree map between closed Haken manifolds satisfying the Volume and Non-degeneration conditions. Let $\Sigma$ be a Seifert piece of $N$ which has a hyperbolic base 2-orbifold. Then for any component $G$ of $G_{\Sigma}$, the space $G_{\rm thick}$ is non-empty and $G_{\rm thin}$ is empty iff $G$ is Seifert. Moreover, if $S$ and $S'$ are distinct components of  $G_{\rm thick}$ then they cannot be adjacent in $M$. 
\end{lemma} 
\begin{proof}
Note that, since $N$ contains a Seifert piece of $\Sigma$ which has a hyperbolic base 2-orbifold, then $\tau(N)\not=0$. 
We first prove that for any component $G$ of $G_{\Sigma}$, $G_{\rm thick}\not=\emptyset$.  Suppose the contrary. Then there exists a component $G$ of $G_{\Sigma}$ such
that $G\subset(G_{\Sigma})_{\rm thin}$.

\emph{First Case.} Suppose first that $N=\Sigma$. Then $G=G_{\Sigma}=M$. It follows from Lemmas \ref{covering} and \ref{tores} (which apply since $\tau(N)\not=0$) that, after passing to a finite covering, $N$ is a circle bundle over an orientable hyperbolic surface and  there exists a maximal vertical torus ${\bf T}$ in $N$ such that, after a homotopy, $f(M)={\bf T}$ and in particular $f_{\ast}(\pi_1M)\simeq{\Z}\times{\Z}\subset\pi_1{\bf T}\subset\pi_1N$. Since $f$ has nonzero degree then this implies that $\pi_1{\bf T}$ is a finite index subgroup  of $\pi_1N$. We write the exact sequence of the fibration
$$1\to{\Z}=\l t\r\stackrel{i_{\ast}}{\to}\pi_1N\stackrel{p_{\ast}}{\to}\pi_1F\to 1$$
Then $\pi_1{\bf T}=\l i_{\ast}(t), b\r$, where $b\in\pi_1N\setminus\ker p_{\ast}$. Thus, since $p_{\ast}$ is surjective then this implies that $\l p_{\ast}(b)\r$ is a finite index subgroup of $\pi_1F$. Since $\pi_1F$ is torsion free this implies that $\pi_1F$ is infinite cyclic. This is a contradiction since $F$ is a hyperbolic surface.

\emph{Second Case.} Assume that  $N\not=\Sigma$ (in particular, ${\mathcal T}_N\not=\emptyset$ and $\b\Sigma\not=\emptyset$). Since  $G\subset\mathcal(G_{\Sigma})_{\rm thin}$ and since $\b\Sigma\not=\emptyset$ then by Corollary \ref{bouts} there exists a component ${\bf T}$ of $\b\Sigma$ such that, after deforming $f$ by a homotopy, then $f(G)={\bf T}$. This operation strictly reduces the number of components of $f^{-1}(\Sigma)$. This contradicts the minimality condition.

 Let  $G$ denotes a component of  $G_{\Sigma}$. It follows from the
above arguments that if $G$ is a Seifert fibered space then
$f_{\ast}(\pi_1G)$ is necessarily a non-abelian group and thus
$G_{\rm thin}$ is empty.

Assume that   $G$ is a graph manifold that is not a Seifert
fibered space. We know that there exists at least  one Seifert
piece  $S_1$ of $G$ such that  $f_{\ast}(\pi_1S_1)$ is non-abelian and
which is adjacent along a  characteristic torus  $T_1$ to a 
 Seifert piece, denoted by  $S_2$ in $G$. Fix a
 base point $x$ in $T_1$ and denote by $h_i$, $i=1,2$, the homotopy class
of the regular fiber in $S_i$ represented in $T_1$. Since
$f_{\ast}(\pi_1S_1)$ is non-abelian,  then $f_{\ast}(h_1)$ is a
power of a fiber in  $\Sigma$. Thus, using the same argument if  $f_{\ast}(\pi_1S_2)$ is
non-abelian we show, by the minimality of the JSJ decomposition, that  the map $f|T_1$ cannot be  $\pi_1$-injective.
This a contradiction since $f$ is non-degenerate when
restricted to any component of the JSJ-family of tori. This completes the proof of
the lemma.
\end{proof}

\section{Simplicial volume and non-degenerate maps of nonzero degree} 
We first study  the volume function ${\rm Vol}(N)$ under nonzero degree maps $f\co M\to N$ satisfying the Volume Condition  and the Non-Degeneration Condition.
The main purpose of this section is to state the following 
\begin{proposition}\label{fine}
Let $f\co M\to N$ be a nonzero degree map between closed Haken manifolds such that $\tau(N)\not=0$ and satisfying the Volume and Non-Degeneration Conditions. If $M_{\rm thin}\not=\emptyset$ then ${\rm Vol}(M)>{\rm Vol}(N)$.
\end{proposition}
\begin{remark}\label{question}
Roughly speaking, in the proof of Proposition \ref{fine} we establish the following inequality: ${\rm Vol}(N)\leq{\rm Vol}(M_{\rm thick})+\varepsilon$ where $\varepsilon\ll{\rm Vol}(M_{\rm thin})$, when $M_{\rm thin}\not=\emptyset$. This implies that ${\rm Vol}(N)<{\rm Vol}(M)$ when $M_{\rm thin}\not=\emptyset$ and this inequality is sufficient for our purpose. However the following question is natural: Is it true that ${\rm Vol}(N)\leq{\rm Vol}(M_{\rm thick})$? Note that throughout the proof of Proposition \ref{fine} we will state two special casis where inequality ${\rm Vol}(N)\leq{\rm Vol}(M_{\rm thick})$ hold (see paragraphs 4.1 and 4.3).   
\end{remark}

The strategy of the proof is (for technical reasons) to find a convenient finite covering $\t{f}\co\t{M}\to\t{N}$ of $f\co M\to N$ such that ${\rm Vol}(\t{M})>{\rm Vol}(\t{N})$. However it is not true in general that this inequality descends to a strict inequality between $M$ and $N$. This phenomenon depends on a relation between the fiber degree of $\t{N}^{\ast}\to N^{\ast}$ and that of $\t{M}^{\ast}\to M^{\ast}$. Note that when $f$ is not $\pi_1$-surjective then $\t{M}$ is not necessarily connected when $\t{f}\co\t{M}\to\t{N}$ denotes a finite covering of $f$. To avoid this kind of situation we define for a nonzero degree map the finite \emph{$f$-coverings}. More precisely, let $f\co M\to N$ be a proper nonzero degree map between orientable 3-manifolds.  Consider the finite covering $\pi\co\hat{N}\to N$ of $N$ corresponding to the finite index subgroup $f_{\ast}(\pi_1M)$ of $\pi_1N$ and denote by $\hat{f}\co M\to\hat{N}$ the lifting of $f$ so that $\pi\circ\hat{f}=f$. Let $q\co\t{N}\to\hat{N}$ denote a finite covering of $\hat{N}$. Denote by $p\co\t{M}\to M$ the finite covering corresponding to the subgroup $\hat{f}^{-1}_{\ast}(\pi_1\hat{N})$ in $\pi_1M$ and let $\t{f}\co\t{M}\to\t{N}$ be the map such that $q\circ\t{f}=\hat{f}\circ p$.
$$\xymatrix{
\t{M} \ar[r]^{\t{f}} \ar[d]_{p} & \t{N} \ar[d]^{q}\\
M \ar[r]^{\hat{f}} \ar[d]_{{\rm Id}} & \hat{N} \ar[d]^{\pi}\\
M \ar[r]^{f}  & N \\
}$$
 Then we say that $\t{f}$ is a finite $f$-covering of $f$ (i.e.  the finite $f$-coverings of $f$ are the finite coverings of $\hat{f}$).   
Using this definition  $\t{M}$ is always connected and ${\rm deg}(p)={\rm deg}(q)$.    As a special case of this remark we state the following result which will be convenient:
\begin{lemma}\label{huge}
Let $f\co M\to N$ be a nonzero degree map between closed Haken manifolds and assume that there exists a finite $f$-covering  $\t{f}\co\t{M}\to\t{N}$ of $f\co M\to N$ such that ${\rm Vol}(\t{M})>{\rm Vol}(\t{N})$. If the covering $\t{N}\to\hat{N}$ has a trivial fiber degree when restricted to the Seifert pieces of $N$ then ${\rm Vol}(M)>{\rm Vol}(N)$. 
\end{lemma}
\begin{proof}
We keep the same notations as above.  Assume that ${\rm Vol}(\t{N})<{\rm Vol}(\t{M})$. Since $q\co\t{N}\to\hat{N}$ induces the trivial covering over  the fibers, then, by Lemma \ref{cover}, ${\rm Vol}(\t{N})={\rm deg}(q){\rm Vol}(\hat{N})$. On the other hand,   ${\rm Vol}(\t{M})\leq{\rm deg}(p){\rm Vol}(M)$. Since ${\rm deg}(p)={\rm deg}(q)$ then ${\rm Vol}(M)>{\rm Vol}(\hat{N})$ and since by Lemma \ref{cover} ${\rm Vol}(\hat{N})\geq{\rm Vol}(N)$ then ${\rm Vol}(M)>{\rm Vol}(N)$. This completes the proof of the lemma.
\end{proof}
As a consequence of this construction we state the following 
\begin{claim}\label{debut}  Throughout the proof of Proposition \ref{fine}, we may assume, that
$f\co M\to N$ is a $\pi_1$-surjective nonzero degree map satisfying the Volume and Non-degeneration conditions  and  each Seifert piece of $N$ is a ${\S}^1$-bundle over an orientable hyperbolic surface and each Seifert piece of $M_{\rm thick}$ has a base 2-orbifold of genus at least 2.
\end{claim}
\begin{proof}
 Let $\pi\co\hat{N}\to N$ be the finite covering of $N$ corresponding to $f_{\ast}(\pi_1M)$ and denote by $\hat{f}$ the lifting of $f$. Since $\tau(\hat{N})\geq\tau(N)>(0,0)$ by Lemma \ref{cover}  thus  Lemma \ref{covering} implies that there exists a finite covering $q\co\t{N}\to\hat{N}$, inducing the trivial covering over the JSJ-family, such that each Seifert piece of $\t{N}$ is a circle bundle over an orientable hyperbolic surface. Hence it is sufficient to apply Lemma \ref{huge}. 
 The second part of the claim follows directly from Lemma \ref{nonabgenus} and Lemma \ref{huge}.
\end{proof}

  Let $\Sigma$ be a component of ${\mathcal S}(N)$ and let $G$ be a component of $G_{\Sigma}=f^{-1}(\Sigma)$ such that $f|G\co G\to\Sigma$ is a proper nonzero degree map. Denote by $F_{\Sigma}$ the base of $\Sigma$ and by $p\co\Sigma\to F_{\Sigma}$ the bundle projection. 
 Suppose first that $G_{\rm thin}=\emptyset$. Since $f|G$ satisfies the conditions of Lemma \ref{adjacent} then $G$ is a Seifert fibered space. Then $f|G$ is a proper allowable bundle homomorphism and, since $\pi_1F_{\Sigma}$ is torsion free, then $f$ induces a nonzero degree map $\o{f}\co\o{\mathcal O}_G\to F_{\Sigma}$, where $\o{\mathcal O}_G$ is the base surface of $G$. This proves that ${\rm Vo}(G)\geq{\rm Vol}(\Sigma)$.

 From now on, one can assume $G_{\rm thin}\not=\emptyset$. Under this additional hypothesis we have to check that
 $${\rm Vol}(G)>{\rm Vol}(\Sigma)\ \ \ \ \ \ \ \ \ \ ({\mathcal V})$$ 
 Recall that by Lemma \ref{tores}, $\Sigma$ admits a finite covering $\t{\Sigma}\to\Sigma$ with trivial fiber degree such that  there exists a finite collection of vertival tori $\t{\mathcal T}_v$ in $\t{\Sigma}$ such that $\t{f}(\t{G}_{{\rm thin}})=\t{\mathcal T}_v$ where $\t{f}|\t{G}$ is the covering of $f|G$ corresponding to  $\t{\Sigma}\to\Sigma$.
 
 Thus in order to prove inequality $({\mathcal V})$ one can assume, using  Lemma \ref{huge}, that the following property $({\mathcal P}_t)$ is satisfied:
 
  $({\mathcal P}_t)\ \ \ there\ exists\ a\ finite\ collection\ of\ vertival\ tori\ {\mathcal T}_v\ in\ \Sigma\ such\ that\ f(G_{{\rm thin}})={\mathcal T}_v\ and\ that\ f|G\ is\ \pi_1-surjective.$

\subsection{A special Case of "thick-domination"}  
In this paragraph we prove the following result. Notations and hypothesis are the same as above.
\begin{lemma}\label{special}
Assume that either ${\mathcal T}_v=\emptyset$ or the tori of the family ${\mathcal T}_v$ are pairwise disjoint. Then we get the following inequality:
$$      {\rm Vol}(G_{\rm thick})\geq{\rm Vol}(\Sigma)      $$
\end{lemma}
\begin{proof}
If ${\mathcal T}_v=\emptyset$ then necessarily $G_{\rm thin}=\emptyset$ and then $G=G_{\rm thick}$ is a Seifert fibered space by Lemma \ref{adjacent}. Hence in this case $f|G\co G\to\Sigma$ is a bundle homomorphism of non zero degree and then the inequality follows. 

From now on,  we assume $G_{\rm thin}\not=\emptyset$. 
 Changing $f$ by a homotopy which is constant when restricted on $G_{\rm thin}$ and so that $f|G_{\rm thick}$ is generic and  using standard cut and paste arguments, it is easy to see, using Lemma \ref{tores},  that $(f|G)^{-1}({\mathcal T}_v)$ is made of $G_{\rm thin}$ union a collection of two sided properly embedded  incompressible surfaces ${\mathcal U}_{\rm thick}$ in $G_{\rm thick}$. We claim that each component of ${\mathcal U}_{\rm thick}$ is a vertical surface. Indeed assume that there is a component $G_1$ of $G_{\rm thick}$ that contains a component $U_1$ of ${\mathcal U}_{\rm thick}$ that is horizontal. Denote by $h_1$ the homotopy class of the generic fiber of the Seifert fibration of $G_1$ and   denote by $A$ the subgroup of $\pi_1 G_1$ generated by $\pi_1U_1$ and $h_1$. Since $f|G_1$ is a bundle homomorphism then it follows from our construction that there exists a component $T$ of  ${\mathcal T}_v$ such that $f_{\ast}(A)\subset\pi_1T$. On the other hand, notice that $A$ is a finite index subgroup of $\pi_1G_1$. Since $\Sigma$ is a ${\S}^1$-bundle over an orientable hyperbolic surface and since $\pi_1T$ is maximal then $f_{\ast}(\pi_1G_1)\subset\pi_1T$. A contradiction.

Hence the family   ${\mathcal U}_{\rm thick}$ consists of a collection of vertical tori and properly embedded, in $G_{\rm thick}$, vertical annuli. Denote by $\Sigma'$ the space $\Sigma\setminus{\mathcal T}_v$ and by $G'_{\rm thick}$ the space $G_{\rm thick}\setminus{\mathcal U}_{\rm thick}$. Then we get a proper  nonzero degree map $f'\co G'_{\rm thick}\to\Sigma'$. Notice that since each component of $G_{\rm thick}$ is a Seifert piece of $M$ and since the components of  ${\mathcal T}_v$ and ${\mathcal U}_{\rm thick}$ are vertical annuli and tori then the spaces $G'_{\rm thick}$ and $\Sigma'$ are endowed with a Seifert fibration and $f'$ is still a bundle homomorphism. More precisely, the base surface of $\Sigma'$  denoted by  $F'_{\Sigma}$ is obtained from $F_{\Sigma}$ after cutting $F_{\Sigma}$ along  some circles obtained as the images of the components of ${\mathcal T}_v$ under the Seifert projection $p\co\Sigma\to F_{\Sigma}$. Notice that, since $\Sigma$ is an orientable circle bundle over $F_{\Sigma}$, that is also orientable, then  $\Sigma'$ is a circle bundle over $F'_{\Sigma}$ and since $\b\Sigma'\not=\emptyset$ then $\Sigma'\simeq F'_{\Sigma}\times{\S}^1$.   On the other hand $G'_{\rm thick}$ has base 2-orbifold ${\mathcal O}'_{\rm thick}$ obtained from  ${\mathcal O}_{\rm thick}$ after cutting along some circles and properly embedded arcs corresponding to the images of the component of ${\mathcal U}_{\rm thick}$ under the Seifert projection $p_{G_{\rm thick}}\co G_{\rm thick}\to{\mathcal O}_{\rm thick}$. Notice that it follows from our construction that 
$${\rm Vol}(\Sigma)={\rm Vol}(\Sigma')\ \ \ \ \ \ \ {\rm and}\ \ \ \ \ \ \ \ {\rm Vol}(G_{\rm thick})\geq{\rm Vol}(G'_{\rm thick}) \ \ \ \ \ (1)$$
Denote by ${\mathcal F}'$ the  surface $(f')^{-1}(F'_{\Sigma})$ which can be assumed to be incompressible and properly embedded. Notice that since $f'$ is a bundle homomorphism then ${\mathcal F}'$ is necessarily a horizontal surface and the map $f'|{\mathcal F}'\co{\mathcal F}'\to F_{\Sigma}'$ factors throught $\o{\mathcal O}'_{\rm thick}$, where $\o{\mathcal O}'_{\rm thick}$ denotes the underlying surface of ${\mathcal O}'_{\rm thick}$ so that  the following  diagram is consistant 
$$\xymatrix{
{\mathcal F}' \ar[r]^{f'|{\mathcal F}'} \ar[d]_{p_{G'_{\rm thick}}} & F'_{\Sigma}\\
\o{\mathcal O}'_{\rm thick} \ar[ur]_{f''}}$$
Since $f'|{\mathcal F}'$ has nonzero degree and since $p_{G'_{\rm thick}}|{\mathcal F}'$ is actually a branched finite covering then the map $f''$ has nonzero degree too. This implies that ${\rm Vol}(G'_{\rm thick})\geq{\rm Vol}(\Sigma')$.    Thus, using relations (1) we get ${\rm Vol}(G_{\rm thick})\geq{\rm Vol}(\Sigma)$.
\end{proof}
We are now ready to begin the proof of inequality $({\mathcal V})$ in the general case.
We first localize the geometric pieces of $G$ which are \emph{efficient} with respect to the domination via the simplicial volume. 
\subsection{Efficient surfaces}
First assume that the Euler number ${\bf e}(\Sigma)$ of the ${\S}^1$-fibration of $\Sigma$ is zero, which corresponds to the ${\Hi}^2\times{\R}$-geometry.  Then $\Sigma$ is homeomorphic to the product $F_{\Sigma}\times{S}^1$. Possibly after changing $f$ by a homotopy so that $f^{-1}(F_{\Sigma})$ is a collection of properly embedded incompressible orientable surfaces, consider  a component ${\mathcal F}$ of  $f^{-1}(F_{\Sigma})$ such that ${\rm deg}(f|{\mathcal F}\co{\mathcal F}\to F_{\Sigma})\not=0$ and denote by $G^{\rm eff}$ the Seifert pieces of $G$ which meet ${\mathcal F}$. 

The thick-thin decomposition of $G$ gives a decomposition of ${\mathcal F}$ into the disjoint union ${\mathcal F}_{\rm thick}\coprod{\mathcal F}_{\rm thin}$. 
Since $f|G_{\rm thick}\co G_{\rm thick}\to\Sigma$ is a bundle homomorphism, then each component of ${\mathcal F}_{\rm thick}$ is a horizontal surface and since ${\mathcal F}$ is connected then each component of 
 $${\mathcal F}^{\ast}_{\rm thin}={\mathcal F}_{\rm thin}\setminus({\mathcal T}_G\cap{\mathcal F}_{\rm thin})$$
  is either a properly embedded vertical annulus  or a horizontal surface. Let $F_{\rm hor}$ denote a horizontal component of ${\mathcal F}^{\ast}_{\rm thin}$. If $p\co\Sigma\to F_{\Sigma}$ denotes the bundle projection then $p_{\ast}f_{\ast}(\pi_1 F_{\rm hor})$ and $p_{\ast}f_{\ast}([c])$ are  infinite cyclic subgroups in $\pi_1F_{\Sigma}$, which are contained in the fundamental group  generated by simple closed curve in $F_{\Sigma}$, for any component $c$ of $\b F_{\rm hor}$. To see this, it is sufficient to check that for any component $c$ of $\b F_{\rm hor}$ then  $f_{\ast}([c])\not=1$ and this follows directly from the Non-degeneration condition.
  \begin{remark}\label{epaisse}
  Let $T$ be a canonical torus of $G$ and suppose that $T\not\subset\b G$. Then $T$ is the boundary of two Seifert pieces $S$ and $S'$ of $G$. If $S=S'$ then we denote by $\Sigma(T)$ a regular neighborhood of $T$. Then in the following it will be convenient, for technical reasons, to consider the spaces $\Sigma(T)$ as a Seifert piece of $G$. Actually $\Sigma(T)$ can be seen as a component of $G_{\rm thin}^{\ast}$. Note that each component of ${\mathcal F}\cap\Sigma(T)$ is always a properly embedded annulus. 
  \end{remark}

Denote by $G^{\rm eff}_{\rm thick}$, resp. $G^{\rm eff}_{\rm thin}$, the space which consists of the disjoint union  of the components $S$ of $G_{\rm thick}$, resp. of $G^{\ast}_{\rm thin}$, such that ${\mathcal F}\cap S\not=\emptyset$. Next we decompose $G^{\rm eff}_{\rm thin}$ into the union $G^{\rm eff}_{\rm thin,hor}\cup G^{\rm eff}_{\rm thin,ver}$ where $G^{\rm eff}_{\rm thin,hor}$, resp. $G^{\rm eff}_{\rm thin,ver}$, denotes the components $S$ of $G^{\rm eff}_{\rm thin}$ such that ${\mathcal F}\cap S$ consists of horizontal, resp. vertical, surfaces. Denote by ${\mathcal F}_{\rm thin,hor}\cup{\mathcal F}_{\rm thin,ver}$ the corresponding decomposition of  ${\mathcal F}_{\rm thin}$. To finish, decompose $G^{\rm eff}_{\rm thin,hor}$ into the union $G^{\rm eff,h}_{\rm thin,hor}\cup G^{\rm eff,\sim}_{\rm thin,hor}$ where $G^{\rm eff,h}_{\rm thin,hor}$ consists of the Seifert pieces $S$ of $G^{\rm eff}_{\rm thin,hor}$ such that $f_{\ast}(h_S)\in\l t\r$, where $h_S$ denotes the generic fiber of $S$ and $t$ is the fiber of $\Sigma$. This gives a decomposition of ${\mathcal F}$ into ${\mathcal F}_{\rm thick}\cup{\mathcal F}_{\rm thin, hor}^{h}\cup{\mathcal F}_{\rm thin, hor}^{\sim}\cup{\mathcal F}_{\rm thin, ver}$. Note that two components of ${\mathcal F}_{\rm thick}\cup{\mathcal F}_{\rm thin, hor}^{h}$, resp. of ${\mathcal F}_{\rm thin, ver}$, cannot be adjacent. This follows directly from the minimality of the JSJ-decomposition combined with the non-degeneration of the map $f|{\mathcal T}_G\co{\mathcal T}_G\to N$.

On the other hand, given a component $F$ of ${\mathcal F}^{\ast}={\mathcal F}\setminus{\mathcal F}\cap{\mathcal T}_G$, we decompose $\b F$ into $\b_{\rm ext} F\cup\b_{\rm int} F$ where $\b_{\rm ext}F=\b F\cap\b G$. Notice that since $f|{\mathcal F}\co{\mathcal F}\to F_{\Sigma}$ has nonzero degree, then as in the proof of Lemma \ref{adjacent}, ${\mathcal F}_{\rm thick}\not=\emptyset$.

\begin{figure}[htb]
\psfrag{ext}{$\b_{\rm ext}$}
\psfrag{thick}{$G^{\rm eff}_{\rm thick}$}
\psfrag{anneau}{${\S}^1\times I$}
\psfrag{vertical}{$G^{\rm eff}_{\rm thin,ver}$}
\psfrag{tilde}{$G^{\rm eff,\sim}_{\rm thin,hor}$}
\psfrag{h}{$G^{\rm eff,h}_{\rm thin,hor}$}
\psfrag{W}{$W(T)\simeq{\S}^1\times{\S}^1\times I$}
\centerline{\includegraphics{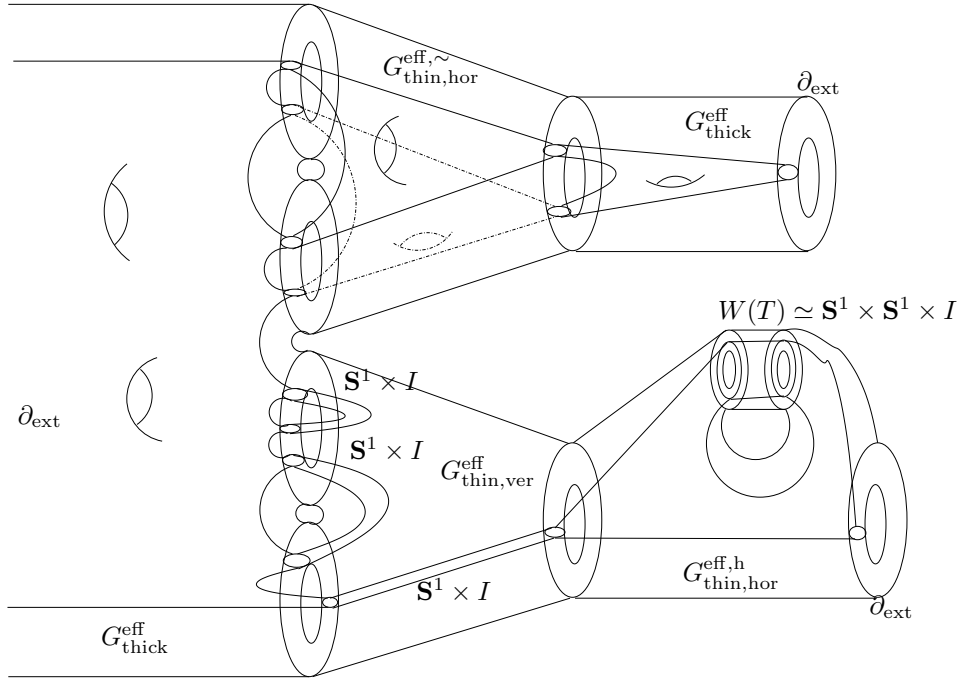}}
\caption{Efficient surfaces and thin-thick decomposition}
\end{figure}

Assume now that ${\bf e}(\Sigma)\not=0$. In particular this means that $\Sigma$ is closed. Since $f|G_{\rm thick}\co G_{\rm thick}\to\Sigma$ is a bundle homomorphism then one can choose a fiber $t$ in $\Sigma\setminus W({\mathcal T}_v)$ such that $f^{-1}(t)$ is a finite unions of fibers $h_1,...,h_l$ in ${\rm int}(G_{\rm thick})$. Then denote by $\Sigma'$, resp. $G'$, the space $\o{\Sigma\setminus W(t)}$, resp. $\o{G\setminus\cup_{i=1,...,l}W(h_i)}$, and by $f'\co G'\to \Sigma'$ the induced proper nonzero degree map. Now $\Sigma'$ is a circle bundle over a surface $F_{\Sigma'}$   with zero Euler number. Note that $F_{\Sigma'}$ is  $F_{\Sigma}$  minus the interior of a slight 2-disk.  As in the first case consider  a component ${\mathcal F}$ of  $f'^{-1}(F_{\Sigma'})$ such that ${\rm deg}(f'|{\mathcal F}\co{\mathcal F}\to F_{\Sigma'})\not=0$. Denote by $\pi_G$, resp. $\pi_{\Sigma}$, the natural quotient map $\pi_G\co G'\to G$, resp. $\pi_{\Sigma}\co{\Sigma}'\to \Sigma$. Then one can perform a decomposition of ${\mathcal F}$ and of $G'$ into $G'^{\rm eff}_{\rm thick}\cup G'^{\rm eff,h}_{\rm thin, hor}\cup G'^{\rm eff,\sim}_{\rm thin, hor}\cup G'^{\rm eff}_{\rm thin, ver}$ as in the first case and the projection $\pi_G$ gives the decomposition  $G^{\rm eff}_{\rm thick}\cup G^{\rm eff,h}_{\rm thin, hor}\cup G^{\rm eff,\sim}_{\rm thin, hor}\cup G^{\rm eff}_{\rm thin, ver}$ of $G^{\rm eff}=\pi_G(G'^{\rm eff})$. 

\begin{remark}\label{coconne}
Note that it follows from our construction that $\pi_G|G'_{\rm thin}\co G'_{\rm thin}\to G_{\rm thin}$ is the identity map and then ${\mathcal F}_{\rm thin}$ can be regarded as a surface in $G$ properly embedded in $G_{\rm thin}$. However ${\mathcal F}_{\rm thick}$ cannot be properly embedded in  $G_{\rm thick}$. This comes from the observation that each components $c_i$ of  ${\mathcal F}_{\rm thick}\cap\b W(h_i)$  is related to the meridian $m_i$ of $W(h_i)$ by the equation $c_i=m_i^{a_i}h^{n_i}$ where $a_i\not=0$  and $(a_i,n_i)=1$. This $n_i$ is generally non-zero. 
\end{remark}
Since $\Sigma'$ is the trivial orientable circle bundle over $F_{\Sigma'}$ then we denote by $i'\co F_{\Sigma'}\to N'$ the canonical inclusion and by $p'\co \Sigma'\to F_{\Sigma'}$ the bundle projection. 
In the case where ${\bf e}(\Sigma)\not=0$ the following map
$$I:\ \ \ \ {\mathcal F}\stackrel{f'|{\mathcal F}}{\to}F_{\Sigma'}\stackrel{i'}{\to}\Sigma'\stackrel{\pi_{\Sigma}}{\to}\Sigma\stackrel{p}{\to}F_{\Sigma}$$
 plays a crucial role in the proof of inequality ${\mathcal V}$ since it measures the contribution of the thin and thick parts to produce the simplicial volume of the target. Note that it is easy to check that the map $p\circ\pi_{\Sigma}\circ  i'$ is surjective at the $\pi_1$-level and thus the induced  homomorphism $I_{\ast}$ as a finite index image in $\pi_1F_{\Sigma}$.
 
 Since $f'|G'_{\rm thick}$ and $f'|G'^{\rm eff, h}_{\rm thin, hor}$ is a bundle homomorphism then $I|{\mathcal F}_{\rm thick}\cup{\mathcal F}^{h}_{\rm thin, hor}$ factors throught $\o{\mathcal O}'^{\rm eff}_{\rm thick}\cup\o{\mathcal O}^{\rm eff,h}_{\rm thin, hor}$, which denotes the union of the base surfaces of the Seifert pieces of $G'^{\rm eff}_{\rm thick}\cup G'^{\rm eff, h}_{\rm thin, hor}$ (recall that $G'_{\rm thin}=G_{\rm thin}$),  in such a way that there is a map 
$$I'\co\o{\mathcal O}'^{\rm eff}_{\rm thick}\cup\o{\mathcal O}^{\rm eff,h}_{\rm thin, hor}\to F_{\Sigma}$$ 
such that 
$$I|{\mathcal F}_{\rm thick}\cup{\mathcal F}^h_{\rm thin,hor}\simeq I'\circ\left(p_{G'^{\rm eff}_{\rm thick}}\cup p_{G'^{\rm eff, h}_{\rm thin, hor}}\right)$$
 where $p_{G'_{\rm thick}}\cup p_{G'^{\rm eff, h}_{\rm thin, hor}}$ denotes the Seifert projections. Denote by $c$ the boundary of $F_{\Sigma'}$ in $\b\Sigma'$ and denote by $m$ the simple closed curve in $\b\Sigma'$ which is identified with the meridian of a solid torus to obtain $\Sigma$ from $\Sigma'$. Since ${\bf e}(\Sigma)\not=0$ and since $\Sigma$ has no exceptional fiber then there exists a nonzero integer $n$ such that $m=ct^n$, where $t$ denotes the homotopy class of the fiber represented in $\b\Sigma'$.  In the same way denote by $c_i$ a boundary of ${\mathcal F}_{\rm thick}$ in $\b W(h_i)$ and by $m_i$ the meridian of $W(h_i)$. There exists coprime  integers $(a_i,n_i)$ with $a_i\not=0$ such that $c_i=m_i^{a_i}h_i^{n_i}$. Denote by $\delta_i$ the boundary component of $\o{\mathcal O}'^{\rm eff}_{\rm thick}$ corresponding to $\b W(h_i)$. Then $(p_{G'^{\rm eff}_{\rm thick}})_{\ast}(c_i)=[\delta_i]^{\pm a_i}$. On the other hand it follows from our construction that there exists $\alpha_i\in{\Z}$ such that $f'_{\ast}(c_i)=c^{\alpha_i}$ and thus $I_{\ast}(c_i)=p_{\ast}\circ(\pi_{\Sigma})_{\ast}(m^{\alpha_i}t^{-n\alpha_i})=p_{\ast}(t^{-n\alpha_i})=1$. Hence $(I')_{\ast}\circ(p_{G'^{\rm eff}_{\rm thick}})_{\ast}(c_i)=(
I')_{\ast}([\delta_i]^{\pm a_i})=1$. Since $\pi_1F_{\Sigma}$ is torsion free then $(I')_{\ast}([\delta_i]^{\pm 1})=1$ and thus the map $I'$ factors throught $\o{\mathcal O}^{\rm eff}_{\rm thick}$. Thus $I$ induces a map $I''\co\o{\mathcal O}^{\rm eff}_{\rm thick}\cup\o{\mathcal O}^{\rm eff,h}_{\rm thin, hor}\to F_{\Sigma}$ such that $I''\circ\pi_{\mathcal O}=I'$ where $\pi_{\mathcal O}\co\o{\mathcal O}'^{\rm eff}_{\rm thick}\cup\o{\mathcal O}^{\rm eff,h}_{\rm thin, hor}\to\o{\mathcal O}^{\rm eff}_{\rm thick}\cup\o{\mathcal O}^{\rm eff,h}_{\rm thin, hor}$ denotes the natural quotient map.
 
\subsection{An intermediate result when the image of $G_{\rm thin}$ is large at the homological level}  For technical reasons,  the proof of inequality $({\mathcal V})$ depends of the image of $G_{\rm thin}$ at the homological level. For instance the more easy case is when the induced homomorphism is surjective. More precisely we prove the following
 \begin{lemma}\label{plein} If the homomorphism $(f|G_{\rm thin})_{\sharp}\co H_1(G_{\rm thin};{\Q})\to H_1(\Sigma;{\Q})$ is surjective then ${\rm Vol}(G_{\rm thick})\geq{\rm Vol}(\Sigma)$.
 \end{lemma} 
 \begin{proof}
 We first check that the hypothesis implies that $(f|\b G_{\rm thick})_{\sharp}\co H_1(\b G_{\rm thick};{\Q})\to H_1(\Sigma;{\Q})$ is surjective. Since by construction $\b_{\rm int}G_{\rm thick}=\b_{\rm int}G_{\rm thin}$ it is sufficient to check that $(f|\b_{\rm int}G_{\rm thin})_{\sharp}\co H_1(\b_{\rm int}G_{\rm thin};{\Q})\to H_1(\Sigma;{\Q})$ is surjective. Let $L$ be a component of $G_{\rm thin}$. Choose a component $D_L$ of $\b_{\rm int}L\subset\b_{\rm int}G_{\rm thin}$. Then $D_L\subset\b_{\rm int}G_{\rm thick}$. Then we claim that 
$$(A)\ \ \ \ \ \ \ \oplus_{L\subset G_{\rm thin}}H_1(D_L;{\Q})\to H_1(\Sigma;{\Q})$$
is an epimorphism.
 Recall that there exists a maximal vertical torus $T$ in $\Sigma$ such that $f_{\ast}(\pi_1L)\subset\pi_1T$ and $f_{\ast}(\pi_1D)\simeq{\Z}\times{\Z}$ where $D$ denotes a component of $\b_{\rm int}L$. Indeed this follows from property ${\mathcal P}_t$ combined with the Non-degeneration condition. Then we have the following commutative diagramm
 $$\xymatrix{
 \pi_1D_L \ar[r]^{i_{\ast}} \ar[d] & \pi_1L \ar[r]^{f_{\ast}} \ar[d] & \pi_1T\subset\pi_1\Sigma \ar[d]\\
 H_1(D_L;{\Z}) \ar[r]^{i_{\sharp}} & \H{L} \ar[r]^{f_{\sharp}} \ar[ur] & \H{F_{\Sigma}}}
 $$
Then ${\rm Rk}(f_{\sharp}\circ i_{\sharp})={\rm Rk}(f_{\sharp})$. Since the components of $G_{\rm thin}$ cannot be adjacent by construction then 
$$H_1(G_{\rm thin})=\oplus_{L\in G_{\rm thin}}H_1(L)$$ 
which proves that $(A)$ is surjective (with coefficient ${\Q}$) and thus so is $(f|\b_{\rm int}G_{\rm thick})_{\sharp}$. Since $f|G_{\rm thick}$ is a fiber preserving map then it descends to a map $f'\co\o{\mathcal O}_{\rm thick}\to F_{\Sigma}$ such that $f'|\b\o{\mathcal O}_{\rm thick}\co\b\o{\mathcal O}_{\rm thick}\to F_{\Sigma}$ induces   an epimorphism at the $H_1$-level (with coefficient ${\Q}$ and where $\o{\mathcal O}_{\rm thick}$ denotes the base surfaces of the components of $G_{\rm thick}$).  
 
 We are now ready to check that ${\rm Vol}(G_{\rm thick})\geq{\rm Vol}(\Sigma)$. On one hand we know that ${\rm Vol}(\Sigma)=\beta_1(F_{\Sigma})-\epsilon$, with $\epsilon=2$ or 1 depending on whether $\Sigma$ is closed or not, and $${\rm Vol}(G_{\rm thick})\geq\sum_{S\in G_{\rm thick}}\left(2g_S+p_S-2\right)$$
 where $g_S$ denotes the genus of $\o{\mathcal O}_S$ and $p_S$ denotes the number of components of $\b S$. On the other hand we know from the paragraph above that 
  $$\sum_{S\in G_{\rm thick}}p_S\geq\beta_1F_{\Sigma}$$ and by Claim \ref{debut} we know that $g_S\geq 2$ when $S\in G_{\rm thick}$. Thus we get
 $${\rm Vol}(G_{\rm thick})\geq\beta_1(F_{\Sigma})+2\sum_{S\in G_{\rm thick}}\left(g_S-1\right)>{\rm Vol}(\Sigma)$$
 This completes the proof of the claim. 
 \end{proof}
\begin{remark}\label{remarqueplein}
Recall that the image of $G_{\rm thin}$ is a family of vertical tori ${\mathcal T}_v$. Denote by ${\mathcal C}_v$ the family of circles in $F_{\Sigma}$ corresponding to ${\mathcal T}_v$.  Then $H_1({\mathcal T}_v;{\Q})\to H_1(\Sigma;{\Q})$ is surjective iff $H_1({\mathcal C}_v;{\Q})\to H_1(F_{\Sigma};{\Q})$ is surjective and in this case we have ${\rm Vol}(G_{\rm thick})\geq{\rm Vol}(\Sigma)$.
\end{remark}
 
 Thus from now on, when $N$ is a circle bundle with non-zero Euler number we can assume that the following condition  is checked:  
 $$({\mathcal C})\ {\rm the\ homomorphism}\ (f|G_{\rm thin})_{\sharp}\co H_1(G_{\rm thin};{\Q})\to H_1(\Sigma;{\Q})\ {\rm is\ not\ surjective.}$$

\subsection{Parametrization of non-degenerate maps}
We define a set of parameters which characterize the map $f\co G\to\Sigma$. First of all, given a Seifert fibered space $S$ with non-empty boundary, endowed with a fixed fibration, a  generic fiber $h_S$, exeptional fibers $c_1,...,c_{r_S}$ denote by $T_1(S),...,T_{p_S}(S)$ its boundary components and for each $i=1,...,p_S$ denote by $d_i(S)$ a section of $T_i(S)$ so that $d_1(S)+...+d_{p_S}(S)=q_1+...+q_{r_S}$ in $\H{S}$, where each $q_i$ is a chosen section corresponding to the exceptional fiber $c_i$. 

If $S$ denotes  a Seifert piece of $G^{\rm eff}_{\rm thick}\cup G^{\rm eff, h}_{\rm thin, hor}$ then let $q_S$ be the nonzero integer satisfying $f_{\ast}(h_S)=t^{q_S}$. 
If $S$ denotes a  Seifert piece of $G^{\rm eff,\sim}_{\rm thin, hor}$ then recall that there exists a vertical torus  $T_S$ in ${\mathcal T}_v$ such that $f(S)=T_S$. Denote by $u_S$ the simple closed curve of $F_{\Sigma}$ such that $\l[u_S]\r=p_{\ast}(\pi_1T_S)$. Moreover when  $S$ is in $G^{\rm eff,\sim}_{\rm thin, hor}$ then we choose a lifting of $u_S$ denoted by $\o{u}_S$ in $\Sigma$ such that $\l\o{u}_S,t\r=\pi_1T_S$ in the following way:  consider a component $c$ of $\b S\cap {\mathcal F}$ and choose $\o{u}_S$ so that $f_{\ast}([c])\in\l[\o{u}_S]\r$. Let $(\beta_S,\alpha_S)$ be the integers such that  $f_{\ast}(h_S)=\o{u}_S^{\beta_S}t^{\alpha_S}$. 
Note that by definition of $G^{\rm eff,\sim}_{\rm thin, hor}$ then ${\mathcal F}\cap S$ is a horizontal surface and then for each $i=1,...,p_S$ there exists $\gamma^i_S\not=0$ and coprime integers  $(a^i_S,n^i_S)$ with $a^i_S\not=0$  such that $f_{\ast}(d^{a_S^i}_i(S)h_S^{-n^i_S})=\o{u}_S^{\gamma^i_S}$.

If $S$ is a Seifert piece of $G^{\rm eff}_{\rm thin, ver}$ then we denote by $\nu_S$ the non-zero integer such that $p_{\ast}f_{\ast}(h_S)=u_S^{\nu_S}$.

 Then we define the \emph{parameters space} of the maps $f$ by setting
$${\mathcal M}(f,{\mathcal F}):=\left\{
\begin{array}{lll}
q_S & {\rm when}\ S\in G^{\rm eff}_{\rm thick}\cup G^{\rm eff, h}_{\rm thin, hor}\\
(\alpha_S,\beta_S), \gamma^i_S, n_S^i & i=1,...,p_S\ {\rm when}\ S\in G^{\rm eff,\sim}_{\rm thin, hor}\\
\nu_S & {\rm when}\ S\in G^{\rm eff}_{\rm thin, ver}
\end{array}
\right\}$$
\begin{claim}\label{combi}
For any Seifert piece $S$ of $G^{\rm eff,\sim}_{\rm thin,hor}$ the couple $(\beta_S,\alpha_S)$ defined above always satisfies the condition $\beta_S\alpha_S\not=0$. 
\end{claim}
\begin{proof}[Proof of the Claim]
The fact that $\beta_S\not=0$ follows from the definition of $G^{\rm eff,\sim}_{\rm thin,hor}$.

Denote by $F_S$ the subsurface ${\mathcal F}\cap S$. We know from the definition of $G^{\rm eff,\sim}_{\rm thin,hor}$ that $F_S$ is a horizontal surface of $S$. Hence if $\alpha_S=0$ then $f_{\ast}(h_S)=\o{u}_S^{\beta_S}$ and thus, since $f|T$ is $\pi_1$-injective $h_S$ is a component of $F_S\cap T$. This implies that $F_S$ is a properly embedded horizontal annulus in $S$. A contradiction. This proves the claim.
\end{proof}
\subsection{Increasing the genus of the base 2-orbifolds of  the efficient thin part}
In this part we state a technical result which allows to construct suitable coverings which increase the genus of the base of the thin part of $G$ assuming some technical  conditions. More precisely, in this paragraph we state the following
\begin{lemma}\label{genus} Assume that either ${\bf e}(\Sigma)=0$ or if ${\bf e}(\Sigma)\not=0$ then assume condition $({\mathcal C})$.
For any $n\in{\N}^{\ast}$ there exists a finite regular covering $f_n\co G_n\to \Sigma_n$ of $f\co G\to\Sigma$ satisfying the following properties:

(i) Any component of $(G_n)^{\ast}$ over a geometric piece of $G^{\rm eff,\sim}_{\rm thin,hor}$ admits a fibration over a  2-orbifold of genus at least $n$,

(ii) the covering $\pi_n\co\Sigma_n\to\Sigma$ has a fiber degree $\leq$ to the fiber degree of $S_n\to S$ for any  geometric piece $S$ of $G^{\rm eff}$ and for any $S_n$ in $G_n^{\ast}$ over $S$.
\end{lemma}
To prove the lemma when ${\bf e}(\Sigma)\not=0$, we need some technical refinement. More precisely in this case  it will  be convenient to assume that the "base curves" $u_S$ that define the  tori $T_S$ satisfy 
$$({\mathcal C}')\ \ \ {\rm Im}\left(H_1\left(u_S;{\Z}\right)\to H_1\left(F_{\Sigma};{\Z}\right)\right)\not=\{0\} \ \ {\rm when}\ S\ {\rm is\ a\ Seifert\ piece\ of}\ G^{\rm eff}_{\rm thin}$$
This is possible using Lemma \ref{nonseparating} combined with Lemma \ref{huge}.
 On the other hand, when ${\bf e}(\Sigma)\not=0$, then
recall that the group $\pi_1\Sigma$ has a presentation
$$({\mathcal P}_{\bf e})\ \ \ \l t,a_1,b_1,...,a_g,b_g: a^{-1}_ita_i=t, b^{-1}_jtb_j=t, [a_1,b_1]...[a_g,b_g]=t^n\r$$ The integer $n$ has the following interpretation: the group $\pi_1\Sigma$ is obtained as a central extension of $\l t\r={\Z}$ by $\pi_1F_{\Sigma}$ using the exact sequence of the fibration
$$\{1\}\to\l t\r\simeq{\Z}\stackrel{i_{\ast}}{\to}\pi_1\Sigma\stackrel{p_{\ast}}{\to}\pi_1F_{\Sigma}\to\{1\}$$
Recall that central extensions of ${\Z}$ by $\pi_1F_{\Sigma}$ correpond to elements of $H^2(\pi_1F_{\Sigma},{\Z})$ and the integer $n$ is the element of ${\Z}\simeq H^2(\pi_1F_{\Sigma},{\Z})$  corresponding to $\pi_1\Sigma$.  
The following result will be convenient because it allows to increase the integer $n$ without modifying the parameter space of ${\mathcal M}(f,{\mathcal F})$. More precisely:
\begin{lemma}\label{euler}
Assume that ${\bf e}(\Sigma)\not=0$ and that condition ${\mathcal C}$  is satisfied for $\Sigma$. Then for any prime $q$ there exists a finite abelian covering $\pi\co\t{\Sigma}\to\Sigma$ acting trivially on ${\mathcal T}_v$, hence in particular on  $t$ and on $G_{\rm thin}$ via $f|G_{\rm thin}$ by remark \ref{remarqueplein}, such that $\t{n}\in q{\Z}$, where $\t{n}$ is the element of $H^2(\pi_1F_{\t{\Sigma}},{\Z})$ corresponding to $\pi_1\t{\Sigma}$. 
\end{lemma}
\begin{proof}
Let $q$ be a prime number. By condition ${\mathcal C}$ combinned with remark \ref{remarqueplein}, there exists an epimorphism $\epsilon\co\H{F_{\Sigma}}\to{\Z}_q$ such that $\ker\epsilon\supset\H{{\mathcal C}_v}$.  Consider the finite covering $\pi\co\t{\Sigma}\to\Sigma$ induced by $\epsilon$ via $p\co\Sigma\to F_{\Sigma}$. It follows from the construction that $\pi$ acts trivially on  ${\mathcal T}_v$. On the other hand $\t{\Sigma}$ is a ${\S}^1$-bundle over a surface ${F}_{\t{\Sigma}}$ that is the covering of $F_{\Sigma}$ corresponding to $\epsilon$. Note that the inclusion $\pi_1{F}_{\t{\Sigma}}\to\pi_1F_{\Sigma}$ gives a map 
$$H^2(\pi_1F_{\Sigma},{\Z})\simeq{\Z}\ni 1\mapsto q\times 1\in{\Z}\simeq H^2(\pi_1\t{F}_{\Sigma},{\Z})$$
and thus the integer $\t{n}$ corresponding to the fibration of $\t{\Sigma}$ satisfies the equation $\t{n}=qn$. This completes the proof. 
\end{proof}
\begin{remark}
Note that since the covering acts trivially on $t$ and on $G^{\rm eff}_{\rm thin}$ via $f|G^{\rm eff}_{\rm thin}$ then the covering $\t{f}\co\t{G}\to\t{\Sigma}$ does not affect the parameters space of $\t{G}^{\rm eff}$. More precisely for any Seifert piece $S$ in $G^{\rm eff,\sim}_{\rm thin,hor}$, resp. $G^{\rm eff}_{\rm thick}\cup G^{\rm eff,h}_{\rm thin,hor}$, resp. $G^{\rm eff}_{\rm thin,ver}$ and any component $\t{S}$ over $S$ in $\t{G}^{\ast}$ then $\alpha_{\t{S}}=\alpha_S$, $\beta_{\t{S}}=\beta_S$,  $\gamma^{j_i}_{\t{S}}=\gamma^i_S$, resp. $q_{\t{S}}=q_S$, resp. $\nu_{\t{S}}=\nu_S$.   
\end{remark}
Thus from now on one may assume that the following condition is checked for $f\co G\to\Sigma$: There exists a prime $q$ such that
$$({\mathcal C}'')\left\{
\begin{array}{ll}
q>{\rm l.c.m}\left\{
\begin{array}{ll}
q_S,\ S\in G^{\rm eff}_{\rm thick}\cup G^{\rm eff,h}_{\rm thin,hor}\\
\gamma^i_S, \alpha_S,\ S\in G^{\rm eff,\sim}_{\rm thin,hor},\ i=1,...,p_S\\
\nu_S,\ S\in G^{\rm eff}_{\rm thin,ver} \end{array}\right\}\\
n\in q{\Z}
\end{array}
\right\}
$$
where ${\rm l.c.m}$ denotes the lowest common multiple and $n$ is defined in ${\mathcal P}_{\bf e}$. The following result is the key step for the proof of Lemma \ref{genus}.
\begin{lemma}\label{separation}
Assume that either ${\bf e}(\Sigma)=0$ or if ${\bf e}(\Sigma)\not=0$ then assume  conditions $({\mathcal C}), ({\mathcal C}')$ and $({\mathcal C}'')$ are satisfied.
Let $S$ be a geometric piece of $G_{\rm thin,hor}^{\rm eff,\sim}$. Let  $g$ be an element of $\pi_1S$ which denotes either the homotopy class of an exceptional fiber or the homotopy class of a section of a boundary component of $S$. Then there exists a finite group $H$ and an epimorphism $\varphi\co\pi_1\Sigma\to H$ such that:

(i) Separation: $\varphi f_{\ast}(g)\not\in\l\varphi f_{\ast}(h_{S})\r$

(ii) Action on the fibers: Let $p\co\t{\Sigma}\to\Sigma$ denote the covering of $\Sigma$ corresponding to $\varphi$ and for any geometric piece $S$ of  $G^{\rm eff}$ denote by $\pi_S\co\t{S}\to S$ the finite covering of $S$ induced by $p$ via $f|S$. Then $G_h(\pi_S)\geq G_h(p)$. 
\end{lemma}
\begin{proof}
Let $S$ be a geometric piece of $G_{\rm thin,hor}^{\rm eff,\sim}$ and let $g$ be an element of $\pi_1S$ satisfying the hypothesis of the lemma. 

First assume that $g$ is the homotopy class of an exceptional fiber $c$ of $S$ and denote by $\mu>1$ the index of this fiber. Let $(\beta,\alpha)\in{\Z}^2$ such that $f_{\ast}(g)=\o{u}_S^{\beta}t^{\alpha}$.  In particular we have $\beta\mu=\beta_{S}\not=0$.

 Let $p$ be a prime number such that $p|\mu$. According to Lemma \ref{residu} there exists a finite group $H_p$ and an epimorphism $\tau\co\pi_1F_{\Sigma}\to H_p$ such that $\tau(u_S^{\beta})\not=1$ and $p$ divides the order of $\tau(u_S^{\beta})$.  Consider the homomorphism  $\varphi$ given by  $$\pi_1\Sigma\stackrel{p_{\ast}}{\to}\pi_1F_{\Sigma}\stackrel{\tau}{\to}H_p$$
This completes the proof when $g=c$. Indeed suppose that there exists $n\in{\Z}$ such that $\varphi f_{\ast}(g)=\varphi f_{\ast}(h^n_S)$.  Then $\tau(u_S^{\beta})=\tau(u_S^{n\beta\mu})$. Then $p$ divides $1-n\mu$. A contradiction since $p|\mu$. On the other hand, the second point of the lemma is satisfied since the covering on the target corresponding to $\varphi$ acts trivially on the fiber. 

Assume now that $g$ denotes the homotopy class of a section $d$ of a component of $\b S$ According to the notation of paragraph 4.4 we know that there exists $i\in\{1,...,p_S\}$ such that   $d=d_i(S)$. In particular we have $f_{\ast}(d^{a^i_S})=\o{u}^{\gamma^i_S+ n^i_{S}\beta_{S}}_{S}t^{n^i_{S}\alpha_{S}}$ where $a^i_S\not=0$. 

Assume that ${\bf e}(\Sigma)\not=0$. 
From the presentation $({\mathcal P}_{\bf e})$ of $\pi_1\Sigma$ and by condition ${\mathcal C}''$ one sees that   $H_1(\Sigma;{\Z})\simeq {\Z}_{n}\oplus \H{F_{\Sigma}}$ where $n\in q{\Z}$. Since $n\in q{\Z}$ then there exists an epimorphism $\lambda_q\co{\Z}_{n}\to {\Z}_{q}$. On the other hand, it follows from condition ${\mathcal C}'$ that the $u_{S}$'s are  non-trivial elements of $H_1(F_{\Sigma};{\Q})$ (when $S$ runs over the Seifert pieces of  $G_{\rm thin}$). Then there exists a $q$-group $F_q$ and an epimorphism  $\tau_q\co\H{F_{\Sigma}}\to F_q$ such that $\tau_q(u_{S})\not=0$.  

Note that if ${\bf e}(\Sigma)=0$ then $H_1(\Sigma;{\Z})\simeq {\Z}\oplus\H{F_{\Sigma}}$ and thus the above construction still hold.
Consider now the homomorphism $\varphi$ defined by
$$\pi_1\Sigma\to\H{\Sigma}\simeq{\Z}_{n}\oplus \H{F_{\Sigma}}\stackrel{\lambda_q\times\tau_q}{\to}{\Z}_{q}\times F_q$$
Using condition ${\mathcal C}''$ we claim that $\varphi$ satisfies the conclusion of the lemma. First we check point (i).  To see this it is sufficient to check that $\varphi f_{\ast}(d^{a^i_S})\not\in\l\varphi f_{\ast}(h)\r$. Assume  that there exists $n\in{\Z}$ such that $\varphi f_{\ast}(d^{a^i_S})=\varphi f_{\ast}(h^n_S)$. Then using our notations this means that 
$$\tau_q\left(u^{\gamma^i_S+ n^i_{S}\beta_{S}}_{S}\right)=\tau_q(u_S^{n\beta_S})\ \  {\rm and}\ \  \lambda_q\left(t^{n^i_S\alpha_S}\right)=\lambda_q\left(t^{n\alpha_S}\right)$$
 Then $q$ divides $\gamma^i_S+\beta_S(n^i_S-n)$ and $(n^i_S-n)\alpha_S$. Since $(\alpha_S,q)=1$ then $q$ divides $n^i_S-n$ and thus $q$ divides $\gamma^i_S$. A contradiction. It remains to check the second point of the lemma. First it follows from the construction of $\varphi$ that $G_h(p)=q$. On the other hand for any Seifert piece $S$ of $G^{\rm eff}$ then it follows from our construction and from conditions ${\mathcal C}'$ and ${\mathcal C}''$ that $f_{\ast}(h_S)$ has order $q^{r_S}$ with $r_S\geq 1$, since $f_{\ast}(h_S)=t^{q_S}$ and $(q,q_S)=1$ or $f_{\ast}(h_S)=\o{u}_S^{\beta_S}t^{\alpha_S}$ and $(\alpha_S,q)=1$ or $p_{\ast}f_{\ast}(h_S)=u_S^{\nu_S}$ with $(\nu_S,q)=1$.  This completes the proof of the lemma.
\end{proof}
\begin{proof}[Proof of Lemma \ref{genus}]
We follow here the same kind of arguments as in \cite{PS}[Lemma 4.2.1, paragraph 4.2.14] using Lemma \ref{separation}.
Let $S$ be a Seifert piece of $G_{\rm thin,hor}^{\rm eff,\sim}$ and assume that the genus $g_S$ of the base 2-orbifold ${\mathcal O}_S$ of $S$ satisfies $g_S\geq 1$.  Denote by $d_1,...,d_{p_S}$ the chosen section of $\b S$ (with respect to the fixed Seifert fibration of $S$) and let $c_1,...,c_r$ denote the homotopy class of the exceptional fibers of $S$ with index $\mu_1,...,\mu_r$.  Using Lemma \ref{separation}  and Lemma \ref{euler} we know that there exists a homomorphism $\varphi\co\pi_1\Sigma\to K$ onto a finite group such that 

(i) $\varphi f_{\ast}(d_i)\not\in\l\varphi f_{\ast}(h_S)\r$, for $i=1,...,p_S$ and  $\varphi f_{\ast}(c_j)\not\in\l\varphi f_{\ast}(h_S)\r$ for $j=1,...,r$. 

Denote by $p\co\t{S}\to S$ a component of the covering over $S$ induced by $\varphi$ via $f|S$. This covering induces a branched covering of degree denoted by $\sigma$ between the underlying space of the base 2-orbifolds of $S$ and $\t{S}$. Let $\beta_j$ the order of $\varphi f_{\ast}(c_j)$ in $K$ and for  each $i=1,...,p_S$ denote by $r_i$ the number of component of $\b\t{S}$ over $T_i$ and set $n_i=\sigma/r_i$.
  Then the Riemann-Hurwitz formula allows the compute the genus of the base 2-orbifold of $\t{S}$ in the following way:
$$2g_{\t{S}}=2+\sigma\left(p_S+2g_S+r-2-\sum_{i=1}^{i=p_S}\frac{1}{n_i}- \sum_{i=1}^{i=r}\frac{1}{(\mu_i,\beta_i)}\right)$$
By condition (i) one can check that $\sigma\geq 2$, $n_i\geq 2$ for $i=1,...,p_S$ and $(\mu_i,\beta_i)\geq 2$ for $i=1,...,r$. Then since moreover $p_S\geq 1$ it is easy to check that $g_{\t{S}}>g_S$ when $g_S\geq 1$. Note that condition (ii) of Lemma \ref{genus} is garanteed by condition (ii) of Lemma \ref{separation}. 

Assume now that $g_S=0$. We follow here the same construction as in the case $g_S\geq 1$ using Lemma \ref{separation}. The Riemann-Hurwitz formula gives:
$$2g_{\t{S}}\geq 2+\sigma\left(\frac{p_S}{2}-2\right)$$
Hence if $p_S\geq 4$ then $g_{\t{S}}\geq 1$ and we have a reduction to the first case. Assume that $p_S\leq 3$ and perform the same construction as above. Denote by $\t{S}$ the finite covering of $S$ corresponding to $\varphi\circ(f|S)_{\ast}$ and denote by $p_{\t{S}}$ the number of boundary components of $\t{S}$. Then the Riemann-Hurwitz forlmula gives
$$2g_{\t{S}}=2-p_{\t{S}}+\sigma\left(p_S+r-2- \sum_{i=1}^{i=r}\frac{1}{(\mu_i,\beta_i)}\right)$$ 
Assume $p_S=3$. If $p_{\t{S}}\geq 4$ then we have a reduction to the case above. If $p_{\t{S}}=3$ the Riemann-Hurwitz forlmula gives, since $\sigma\geq 2$ then  $2g_{\t{S}}\geq -1+\sigma\geq 1$ and thus $g_{\t{S}}\geq 1$.
Assume $p_S=2$. Applying the same argument we get a reduction to the case  $p_S=3$ or $g_{\t{S}}\geq 1$. Note that the case $p_S=1$ is impossible since $f_{\ast}(\pi_1S)\simeq{\Z}\times{\Z}$ and $f_{\ast}(\pi_1\b S)\simeq{\Z}\times{\Z}$ (Indeed if $\b S$ is connected then ${\rm Rk}(H_1(\b S)\to H_1(S))=1$).  
Next we perform this construction for each component of $G^{\rm eff,\sim}_{\rm thin,hor}$.  

  To complete the proof of the lemma it remains to check that one can find a regular covering. More precisely assume that there exists a finite covering $f_n\co G_n\to\Sigma_n$ satisfying the conclusion of the lemma. Denote by $\pi_n\co\Sigma_n\to\Sigma$ the associated covering of $\Sigma$, by $H_n$ the finite index subgroup of $\pi_1\Sigma$ corresponding to this covering and denote by $p_n\co G_n\to G$ the corresponding finite covering induced by $\pi_n$ via $f\co G\to \Sigma$. Denote by $\epsilon_n\co\hat{\Sigma}_n\to\Sigma_n$ the finite covering so that $\pi_n\circ\epsilon_n$ is the regular covering of $\Sigma$ corresponding to the normal subgroup $$K_n=\bigcap_{g\in\pi_1\Sigma}gH_ng^{-1}\lhd\pi_1\Sigma$$ 
Then consider the induced  regular covering $\hat{f}_n\co\hat{G}_n\to\hat{\Sigma}_n$. Since $G_n$ satisfies point (i) of Lemma \ref{genus} and since  $\hat{G}_n$ is a finite covering of $G_n$ then point (i) is still true for $\hat{G}_n$. On the other hand since the fiber of $\Sigma$ is central in $\pi_1\Sigma$ then it follows from the construction the fiber degree of $\pi_n\circ\epsilon_n$ is equal to the fiber degree of $\pi_n$. Hence the covering $\hat{f}_n\co\hat{G}_n\to\hat{\Sigma}_n$ satisfies point (ii).
 This completes the proof of the lemma.
\end{proof}
\subsection{Comparing the volume}

We use here the efficient surfaces constructed in paragraph 4.2 and we keep the same notations.
 We construct a dual graph $\Gamma$ for ${\mathcal F}$ in the following way: the vertex space $V(\Gamma)$ is the connected components of ${\mathcal F}_{\rm thick}$ and ${\mathcal F}^{\ast}_{\rm thin}$ and the edge space $E(\Gamma)$ consists of the components of $\b_{\rm int}{\mathcal F}_{\rm thick}\cup\b_{\rm int}{\mathcal F}^{\ast}_{\rm thin}$. We assume that $\Gamma$ is embedded in $G$  with the canonical inclusion.   On the other hand, for each edge $e\in E(\Gamma)$ then $e\cap{\mathcal T}_G$ consists of a single point $v_e(T)$, where $T$ denotes the component of ${\mathcal T}_G$ such that $e\cap T\not=\emptyset$. Then the set $\{v_e(T), e\in E(\Gamma), T\in{\mathcal T}_G\}=\Gamma\cap{\mathcal T}_G$ will be termed the \emph{middle space} of $\Gamma$ and we denote it by $M(\Gamma)$.

Denote by $G_{\rm main}=G^{\rm eff}_{\rm thick}\cup G^{\rm eff,h}_{\rm thin,hor}$, and by ${\mathcal F}_{\rm main}={\mathcal F}\cap G_{\rm main}={\mathcal F}_{\rm thick}\cup{\mathcal F}^h_{\rm thin,hor}$.  It will be convenient to assume that the following conditions are satisfied for the graph $\Gamma$:

(i) \emph{Vertex condition}: Given a Seifert piece $S$ of $G^{\rm eff}$, resp. a component $T$ of ${\mathcal T}_{G^{\rm eff}}$,  then we assume that there is a point $x_S$, resp. $x_T$, in $\Sigma$ such that $f(S\cap V(\Gamma))=x_S$, resp. $f(T\cap M(\Gamma))=x_T$.

This is possible after performing a homotopy on $f$ moving only a small regular neighborhood of $V(\Gamma)\cup M(\Gamma)$.

(ii) \emph{Equivariance}: Let $S$ be a Seifert piece of $G_{\rm main}$ and let $a_1,a_2$ denote two components of $S\cap E(\Gamma)$. For each $i=1,2$, then $a_i=[v_i,v_{e_i}(T_i)]$ for some $v_i\in V(\Gamma)$ and $v_{e_i}(T_i)\in M(\Gamma)$. If $T_1=T_2$ then we assume that $\pi_S|a_i\co a_i\to\pi_S(a_i)$ is a homeomorphism and that $\pi_S(a_1)=\pi_S(a_2)$, where $\pi_S\co S\to\o{\mathcal O}_S$ denotes the Seifert projection.

Note that the equivariance condition is possible, after readjusting $\Gamma$, since the Seifert pieces of $G_{\rm main}$ are not adjacent. 
\begin{figure}[htb]
\psfrag{ext}{$\b_{\rm ext}$}
\psfrag{thick}{$G^{\rm eff}_{\rm thick}$}
\psfrag{anneau}{${\S}^1\times I$}
\psfrag{vertical}{$G^{\rm eff}_{\rm thin,ver}$}
\psfrag{tilde}{$G^{\rm eff,\sim}_{\rm thin,hor}$}
\psfrag{h}{$G^{\rm eff,h}_{\rm thin,hor}$}
\psfrag{W}{$W(T)\simeq{\S}^1\times{\S}^1\times I$}
\centerline{\includegraphics{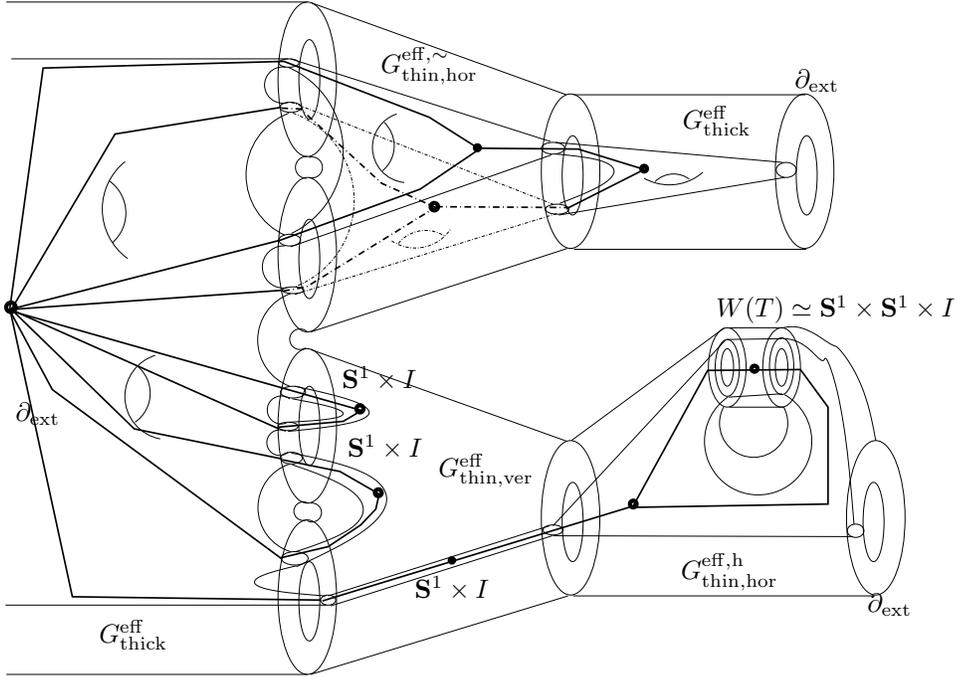}}
\caption{Dual Graph}
\end{figure} 

 \begin{claim}\label{epi} The induced homomorphisms
 $$f_{\sharp}\co H_1\left(\Gamma\cup{\mathcal F}_{\rm main};{\Q}\right)\to H_1(F_{\Sigma},{\Q}),$$
 $${\rm resp.}\ \ I_{\sharp}\co H_1\left(\Gamma\cup{\mathcal F}_{\rm main};{\Q}\right)\to H_1(F_{\Sigma},{\Q})\ \ {\rm when}\ {\bf e}(\Sigma)\not=0$$
 are surjective.
 \end{claim}
 \begin{proof}
 Recall that $f_{\ast}\co\pi_1\left({\mathcal F}\right)\to\pi_1F_{\Sigma}$, resp. $I_{\ast}\co\pi_1\left({\mathcal F}\right)\to\pi_1F_{\Sigma}$, has a finite index image in $\pi_1F_{\Sigma}$. Given a group $G$ denote by ${\mathcal G}_G$ a set of generators of $G$. Then 
 $$\bigcup_{U\in{\mathcal F}_{\rm main}}{\mathcal G}_{\pi_1U}\bigcup_{U\in\left({\mathcal F}\setminus{\mathcal F}_{\rm main}\right)} {\mathcal G}_{\pi_1U}\bigcup{\mathcal G}_{\pi_1\Gamma}$$
  is a set of generator for $\pi_1{\mathcal F}$. On the other hand it follows from our construction that  any component $U$ of ${\mathcal F}\setminus{\mathcal F}_{\rm main}$ has at least one boundary component, say $c_U$ adjacent to a component of ${\mathcal F}_{\rm main}$. On the other hand we know that  $f_{\ast}([c_U])$, resp. $I_{\ast}([c_U])$ is a finite index subgroup of $f_{\ast}(\pi_1U)$, resp. $I_{\ast}(\pi_1U)$. This completes the proof of the claim. 
 \end{proof}

 Denote by $\o{\mathcal O}^{\rm eff}$ the disjoint union of the bases of the Seifert pieces of $G^{\rm eff}$ decomposed as the union $\o{\mathcal O}^{\rm eff}_{\rm thick}\coprod\o{\mathcal O}^{\rm eff,h}_{\rm thin,hor}\coprod\o{\mathcal O}^{\rm eff,\sim}_{\rm thin,hor}\coprod\o{\mathcal O}^{\rm eff}_{\rm thin,ver}$ and by $\o{\mathcal O}_{\rm main}$ the union $\o{\mathcal O}^{\rm eff}_{\rm thick}\coprod\o{\mathcal O}^{\rm eff,h}_{\rm thin,hor}$.

 \subsubsection{Connecting the main surfaces}  
The vertex and equivariance conditions (i) and (ii) give rise  to an equivalence relation on $\Gamma$ denoted by $\sim$. Denote by $\hat{\Gamma}$ the quotient space $\Gamma/\sim$ and by $q\co\Gamma\to\hat{\Gamma}$ the projection. Note that $V(\hat{\Gamma})=q(M(\Gamma)\cup V(\Gamma))$.  Then the vertex and equivariance conditions imply that the map $f|\Gamma$ factors through $\hat{\Gamma}$. More precisely we get the following commutative diagram
$$\xymatrix{
 \Gamma \ar[r]^{f|\Gamma} \ar[d]_{q}  & F_{\Sigma}\\
 \hat{\Gamma} \ar[ur]_{h}
 }$$

Assume that ${\bf e}(\Sigma)=0$. Let $S$ be a component of $G_{\rm main}$. Since $f|S\co S\to\Sigma$ is a bundle homomorphism (by definition) and since $\pi_1F_{\Sigma}$ is torsion free then there is a map  $f'\co\o{\mathcal O}_{\rm main}\to F_{\Sigma}$ such that the following diagram is consistant.
$$\xymatrix{
 {\mathcal F}_{\rm main} \ar[r]^{f} \ar[d]_{\Pi=p_{G_{\rm main}}} & F_{\Sigma}\\
 \o{\mathcal O}_{\rm main} \ar[ur]_{f'}
 }$$ 
 On the other hand when ${\bf e}(\Sigma)\not=0$ then we know from paragraph 4.2 that the  map $I\co{\mathcal F}\to F_{\Sigma}$ induces by restriction a map such that the following diagram is consistant.
$$\xymatrix{
{\mathcal F}_{\rm main} \ar[r]^{I} \ar[d]_{\Pi=\pi_{\mathcal O}\circ p_{G'_{\rm main}}} & F_{\Sigma}\\
\o{\mathcal O}_{\rm main}  \ar[ur]_{I''}
 }$$ 
We now define a "quotient space"  of $\Gamma\cup{\mathcal F}_{\rm main}$  in the following way. First note that $\Pi(\Gamma\cap {\mathcal F}_{\rm main})$ gives a subgraph of $\hat{\Gamma}$. Then define the space $\o{\mathcal O}_{\rm main}\cup\hat{\Gamma}$ as the attachement of $\hat{\Gamma}$ to $\o{\mathcal O}_{\rm main}$ along $\Pi(\Gamma\cap{\mathcal F}_{\rm main})=\hat{\Gamma}\cap\o{\mathcal O}_{\rm main}$. Then the vertex and equivariance conditions allow to extend the above diagram to the following one $({\mathcal D})$:
 $$\xymatrix{
 \Gamma\cup{\mathcal F}_{\rm main} \ar[r]^{\star} \ar[d]_{\Pi\cup q} & F_{\Sigma}\\
 \hat{\Gamma}\cup\o{\mathcal O}_{\rm main} \ar[ur]_{\star'}
 }$$
where $\star$ , resp. $\star'$ represents $f$, resp. $f'\cup h$ or $I$, resp.  $I''\cup h$ depending on whether ${\bf e}(\Sigma)=0$ or not.  Notice that it follows from our construction and from Claim \ref{epi} that  $\star'$ induces an epimorphism at the $H_1$-level (with coefficient ${\Q}$). 
\begin{figure}[htb]
\psfrag{q}{$\Pi\cup q$}
\psfrag{A}{$\o{\mathcal O}_{\rm main}$}
\centerline{\includegraphics{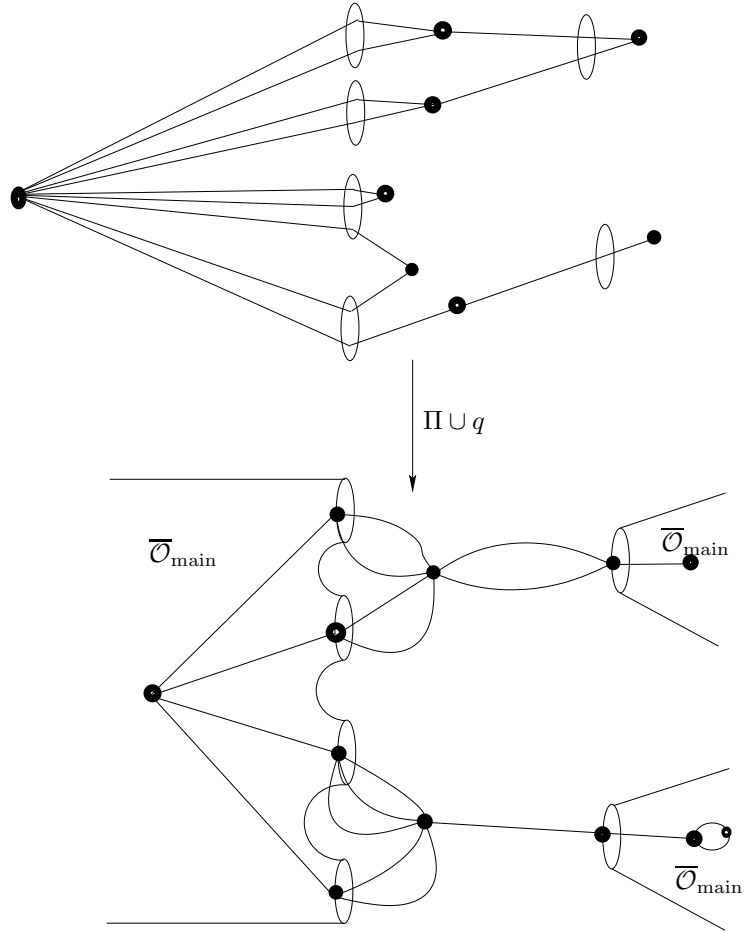}}
\caption{Quotient map}
\end{figure} 
\subsubsection{Vertical  identifications}    Let $\alpha$ and $\beta$ be two components of  $\b\o{\mathcal O}_{\rm main}$.  Then there exist  Seifert pieces $S$ in $G_{\rm main}$ whose base orbifold  ${\mathcal O}_S$ has boundary components $c^1_S,...,c^{p_S}_S$ and $S'$ in $G_{\rm main}$ whose base orbifold  ${\mathcal O}_{S'}$ has boundary components $c^1_{S'},...,c^{p_{S'}}_{S'}$ such that  $\alpha=c^i_S$ and  $\beta=c^j_{S'}$ for some $i\in\{1,...,p_S\}$ and $j\in\{1,...,p_{S'}\}$.  
      We say that $\alpha\equiv\beta$  if there is a simple closed curve $\t{c}^i_S$ in ${\mathcal F}\cap T^i_S$ and a simple closed curve $\t{c}^j_{S'}$ in ${\mathcal F}\cap T^j_{S'}$ and a sequence $A_1,...,A_n$ of vertical annuli of $G_{\rm thin}^{\ast}$ which connects $\t{c}^i_S$ with $\t{c}^j_{S'}$, where $T^i_S$, resp. $T^j_{S'}$, denotes the component of $\b S$, resp. $\b S'$, over $c^i_S$, resp. $c^j_{S'}$. Denote by $\o{\alpha}$ the set which consists of  the components $\beta$ such that $\beta\equiv\alpha$.    Then there exists an element $l$ in $\pi_1F_{\Sigma}$ such that for each $\beta$ in $\o{\alpha}$ there exists an integer $a_{\beta}$  such that $\star'_{\ast}(\beta)=l^{a_{\beta}}$. Let ${\S}^1_{\o{\alpha}}$ denote a circle. Then we glue each component $\beta$ of $\o{\alpha}$  along ${\S}^1_{\o{\alpha}}$  with the attaching maps $z\mapsto z^{a_{\beta}}$. We do that for each equivalence class of $\b\o{\mathcal O}_{\rm main}$. The resulting space is a \emph{branched surface} denoted   by $\o{\mathcal O}^{\rm new}_{\rm main}$  and we call $r$ the quotient map. Note that these identifications do not change the volume of ${\mathcal O}_{\rm main}$.
\begin{figure}[htb]
\psfrag{A}{$\o{\mathcal O}_{\rm main}$}
\psfrag{B}{$\o{\mathcal O}^{\rm new}_{\rm main}$}
\psfrag{r}{$r$}
\centerline{\includegraphics{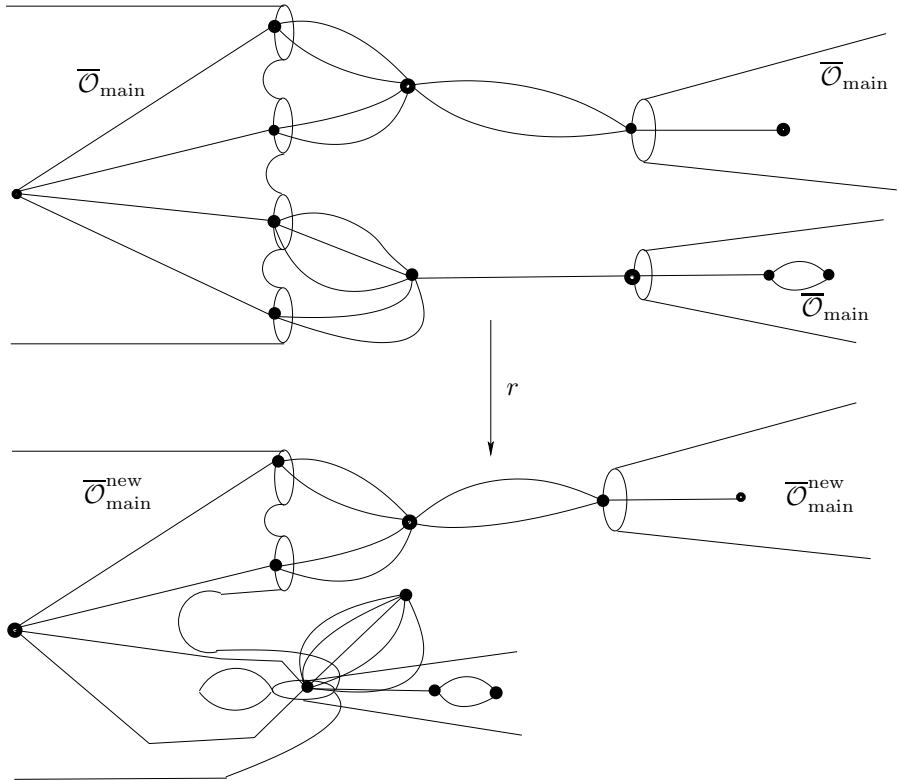}}
\caption{Vertical identifications}
\end{figure} 

Let $x_1$ and $x_2$ denote two vertices of $\hat{\Gamma}$ corresponding to two equivalent boundary components of $\o{\mathcal O}_{\rm main}$. Then we may assume that $r(x_1)=r(x_2)$. Extending $r$ trivially over $\hat{\Gamma}\setminus(\o{\mathcal O}_{\rm main}\cap\hat{\Gamma})$ we get a graph  $r(\hat{\Gamma})$ such that the following diagram is commutative
$$\xymatrix{
\hat{\Gamma}\cup\o{\mathcal O}_{\rm main} \ar[r]_{\star'} \ar[d]_{r} & F_{\Sigma}\\
 r(\hat{\Gamma})\cup\o{\mathcal O}^{\rm new}_{\rm main} \ar[ur]_{\star''}
 }$$

Consider the subgraph $\hat{\Gamma}_0$ of $r(\hat{\Gamma})$ defined as follows: $\hat{\Gamma}_0$ consists of the set, denoted by $r(\hat{\Gamma})_{\rm ext}$, of all edges of $r(\hat{\Gamma})$ which are not contained in $\o{\mathcal O}^{\rm new}_{\rm main}$ union the edges of $r(\hat{\Gamma})$ adjacent to an edge of $r(\hat{\Gamma})_{\rm ext}$. Of course we have $r(\hat{\Gamma})\cup\o{\mathcal O}^{\rm new}_{\rm main}=\hat{\Gamma}_0\cup\o{\mathcal O}^{\rm new}_{\rm main}$.

\subsubsection{Elimination Lemma} In the paragraph above we have constructed a surface complex $\hat{\Gamma}_0\cup\o{\mathcal O}^{\rm new}_{\rm main}$ which dominates $F_{\Sigma}$ at the $H_1$-level.  The purpose of this paragraph is to reduce the graph $\hat{\Gamma}_0$  to eliminate "redundant" generators of the image of $\star''$ at the $H_1$-level. This reduction is crucial to compare the volume of $G$ with $\Sigma$. In the following, given a vertex $x$ of a graph, we denote by $v(x)$ the valence of $x$. The main result of this section is the following

\begin{lemma}\label{comparing} 
There exists a graph $\hat{\Gamma}_1$ which is a subset  of $\hat{\Gamma}_{0}$ whose vertex space satisfies the following conditions:

(i) $V(\hat{\Gamma}_1)\subset r\circ q(V(\Gamma))$,
$${\rm Card}(V(\hat{\Gamma}_1))\geq{\rm Card}(\pi_0(\o{\mathcal O}^{\rm new}_{\rm main}))+{\rm Card}(\pi_0(G^{\rm eff,\sim}_{\rm thin, hor}))+{\rm Card}(\pi_0(\hat{\Gamma}_1))-1$$
On the other hand, the edge space satisfies the following conditions: 

(ii) the set $E(\Gamma_1)$ can be identified with a subset of $\b G^{\rm eff,\sim}_{\rm thin, hor}\cup\b G^{\rm eff}_{\rm thin, ver}$. Hence denote by $E_1$ the edges of $\hat{\Gamma}_1$ which corresponds to a component of $\b G^{\rm eff,\sim}_{\rm thin, hor}$ and by $E_2$ the edges of $\hat{\Gamma}_1$ which corresponds to a component of $\b G^{\rm eff}_{\rm thin, ver}$ which is not adjacent to a component of $G^{\rm eff,\sim}_{\rm thin, hor}$ (i.e. $E_1\cap E_2=\emptyset$). Then the following relations are satisfied: 
  $${\rm Card}(E_1)\leq\sum_{S\in G^{\rm eff,\sim}_{\rm thin, hor}}{\rm Card}(\pi_0(\b S))$$
  $${\rm Card}(E_2)\leq\sum_{S\in G^{\rm eff}_{\rm thin, ver}}\left({\rm Card}(\pi_0(\b S))-2\right)$$
(iii) the space $\hat{\Gamma}_1\cup\o{\mathcal O}^{\rm new}_{\rm main}$ is still connected and the map $\star''|\hat{\Gamma}_1\cup\o{\mathcal O}^{\rm new}_{\rm main}\co\hat{\Gamma}_1\cup\o{\mathcal O}^{\rm new}_{\rm main}\to F_{\Sigma}$ induces an epimorphism at the $H_1$-level (with coefficient ${\Q}$). 
\end{lemma}
In order to prove the lemma we first check the following
\begin{claim}\label{valence}
There exists a subgraph $\hat{\Gamma}'_0$ of $\hat{\Gamma}_0$ satisfying the following conditions:

(i) ${\rm Card}(V(\hat{\Gamma}'_0))\geq{\rm Card}(\pi_0(\o{\mathcal O}^{\rm new}_{\rm main}))+{\rm Card}(\pi_0(G^{\rm eff,\sim}_{\rm thin, hor}))+{\rm Card}(\pi_0(\hat{\Gamma}_1))-1$,

(ii) for any $x\in r\circ q(M(\Gamma))\cap V(\hat{\Gamma}'_0)$,  then $v(x)=2$,

 (iii) the space ${\Gamma}'_0\cup\o{\mathcal O}^{\rm new}_{\rm main}$ is still connected and the map $\star''|\hat{\Gamma}'_0\cup\o{\mathcal O}^{\rm new}_{\rm main}\co\hat{\Gamma}'_0\cup\o{\mathcal O}^{\rm new}_{\rm main}\to F_{\Sigma}$ induces an epimorphism at the $H_1$-level (with coefficient ${\Q}$).
\end{claim}

\begin{proof}
First notice that $\hat{\Gamma}_0$ satisfies points (i) and (iii) by construction. On the other hand it follows from our costruction and from the definition of $M(\Gamma)$ that for any $x\in r\circ q(M(\Gamma))\cap V(\hat{\Gamma}_0)$ then $v(x)\geq 2$. 

Then assume that there exists a point $x\in r\circ q(M(\Gamma))\cap V(\hat{\Gamma}_0)$ such that $v(x)\geq 3$. Then  there exists at least three edges, say $e_1$, $e_2$ and $e_3$ of $\hat{\Gamma}_0$ such that $x$ is an end of $e_i$, for $i=1,...,3$.  For each $i=1,...3$ denote by $y_i$ the end of $e_i$ such that $\b e_i=\{x,y_i\}$.
Note that each $y_i$ is a point of $r\circ q(V(\Gamma))\cap\hat{\Gamma}_0$ and thus each $y_i$ correponds to a unique Seifert piece of $G^{\rm eff}$. 

\emph{First Case:} Assume that there exist at least two elements $i,j$ in $\{1,2,3\}$ such that $y_i$ and $y_j$ are in  $\o{\mathcal O}^{\rm new}_{\rm main}$. Then remove the edge $e_i$ (say) from the graph $\hat{\Gamma}_0$ to get a new graph $\hat{\Gamma}_{0,1}\subset\hat{\Gamma}_0$. By construction $y_i$ and $y_j$ are in the same component of  $\o{\mathcal O}^{\rm new}_{\rm main}$ thus point (i) is still satisfied. 
On the other hand the valence of $x$ in $\hat{\Gamma}_{0,1}$ is strictly less than the valence of $x$ in $\hat{\Gamma}_0$. Point (iii) is easily checked with the space  $\hat{\Gamma}_{0,1}\cup\o{\mathcal O}^{\rm new}_{\rm main}$ since $e_i\subset\o{\mathcal O}^{\rm new}_{\rm main}$.

\emph{Second Case:} Up to re-indexing, assume that $y_1,y_2$ correpond to Seifert pieces of $G^{\rm eff}_{\rm thin}$ and that $y_3\in\o{\mathcal O}^{\rm new}_{\rm main}$.  

 Thus it follows from our construction that there exists edges $\t{e}_1,\t{e}_2,\t{e}_3$ of $\Gamma$ with end points $\{\t{x}_i,\t{y}_i\}$ $i=1,2,3$, Seifert pieces $S_3$ of $G_{\rm main}$, $S_1,S_2$ of $G^{\rm eff,\sim}_{\rm thin,hor}\cup G^{\rm eff}_{\rm thin,ver}$ and boundary $T_3\subset \b S_3$ and $T_i\subset\b S_i$, $i=1,2$ such that:

(1) $r\circ q(\t{e}_i)=e_i$ and $r\circ q(\t{x}_i)=x$,  $r\circ q(\t{y}_i)=y_i$ for $i=1,2,3$,

(2) $\t{x}_3\in T_3$, $\t{x}_1\in T_1$, $\t{x}_2\in T_2$, $\t{y}_3\in{\rm int}(S_3)$, $\t{y}_1\in{\rm int}(S_1)$ and $\t{y}_2\in{\rm int}(S_2)$.

Since $r\circ q(\t{x}_i)=x$ for $i=1,2,3$ then $T_1,T_2$ and $T_3$ are also  components of $G_{\rm main}$ which are connected  by a finite sequence of vertical annuli of $G^{\ast}_{\rm thin}$. Thus necessarily $S_1=S_2=S$ and in particular $y_1=y_2$.

Denote by $c$ the curve defined by $e_1\cup e_2$. It follows from our construction that there exists a curve $\t{c}$ in $S$ such  that $p_{\ast}f_{\ast}([\t{c}])=\star''_{\ast}([c])$ in $\pi_1F_{\Sigma}$, where $p\co\Sigma\to F_{\Sigma}$ denotes the bundle projection.  Denote by $A$ a the connected component of $G_{\rm thin}$ which contains $S$. Then, since ${\mathcal F}$ is connected it follows from our construction that $A$ contains at least one boundary component, say $T_A$, which is adjacent to a component $B$ of $G^{\rm eff}_{\rm thick}$. Note that $B$ is necessarily a Seifert manifold. It follows from our construction that $p_{\ast}f_{\ast}(\pi_1A)\simeq{\Z}$ and that if $\t{s}$ denotes a section of $T_A$ with respect to the Seifert fibration of $B$ then $p_{\ast}f_{\ast}(\l[\t{s}]\r)\simeq{\Z}$. This point comes from the non-degeneration of $f|T_A$ and from the fact that $f|B\co B\to\Sigma$  is a fiber preserving map.  Denote by $s$ the component of $\b\o{\mathcal O}_{\rm main}$ such that $s=\Pi(\t{s})$. Then ${\rm Im}(\star''_{\sharp}([r(s)]))={\rm Im}(\star''_{\sharp}([c]))$ at the $H_1$-level with coefficient ${\Q}$, where $r$ denote the projection $\o{\mathcal O}_{\rm main}\to\o{\mathcal O}^{\rm new}_{\rm main}$.

Consider the graph $\hat{\Gamma}_{0,1}$ obtained from $\hat{\Gamma}_0$ after removing ${\rm int}(e_1)$. Then $\hat{\Gamma}_{0,1}$ satisfies points (i) and (iii) and the valence of $x$  in $\hat{\Gamma}_{0,1}$ is strictly less than the valence of $x$ in $\hat{\Gamma}_0$.

\emph{Third Case:} Assume that $y_1,y_2$ and $y_3$ correspond to Seifert pieces of $G^{\rm eff,\sim}_{\rm thin,hor}\cup G^{\rm eff}_{\rm thin,ver}$. Then we can apply the same construction as in the second case.
This completes the proof of the claim.
  
\end{proof}
  
\begin{figure}[htb]
\psfrag{B}{$\o{\mathcal O}^{\rm new}_{\rm main}$}

\centerline{\includegraphics{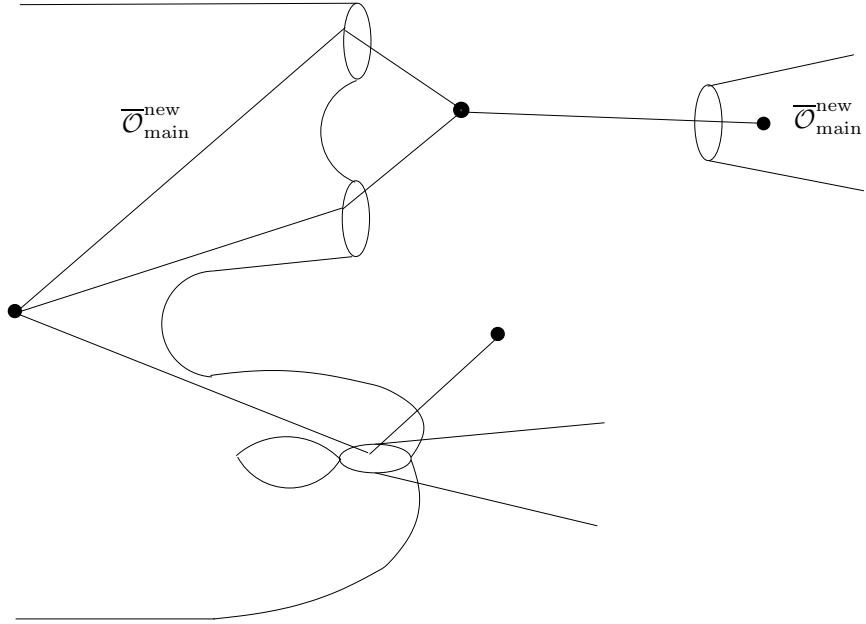}}
\caption{Resulting grah after "Elimination"}
\end{figure}

\begin{proof}[Proof of Lemma \ref{comparing}]

Let $x$ be an element of $r\circ q(M(\Gamma))\cap\hat{\Gamma}'_0$. We know by Claim \ref{valence} that $v(x)=2$. Then  there exists excatly two edges $e_1, e_2$ whose $x$ is an end point. Hence one can replace the edges $e_1,e_2$ by a single edge $e_1\cup_x e_2$. By performing this operation for all points of $r\circ q(M(\Gamma))\cap\hat{\Gamma}'_0$ we get a new graph $\hat{\Gamma}"_0$ satisfying the following properties:

(1) $V(\hat{\Gamma}"_0)\subset r\circ q(V(\Gamma))$,
$${\rm Card}(V(\hat{\Gamma}"_0))\geq{\rm Card}(\pi_0(\o{\mathcal O}^{\rm new}_{\rm main}))+{\rm Card}(\pi_0(G^{\rm eff,\sim}_{\rm thin, hor}))+{\rm Card}(\pi_0(\hat{\Gamma}"_0))-1$$
(2)The edge space $E(\hat{\Gamma}"_0)$ can be identified with a subset of $\b G^{\rm eff,\sim}_{\rm thin, hor}\cup\b G^{\rm eff}_{\rm thin, ver}$ and 
 
 (3) the space $\hat{\Gamma}"_0\cup\o{\mathcal O}^{\rm new}_{\rm main}$ is still connected and the map $\star''|\hat{\Gamma}"_0\cup\o{\mathcal O}^{\rm new}_{\rm main}\co\hat{\Gamma}"_0\cup\o{\mathcal O}^{\rm new}_{\rm main}\to F_{\Sigma}$ induces an epimorphism at the $H_1$-level (with coefficient ${\Q}$), since $\hat{\Gamma}"_0$ and $\hat{\Gamma}'_0$ are homeomorphic.
 
 Then one can consider the disjoint decomposition of $E(\hat{\Gamma}"_0)$ into $E_1\cup E_2$.
 To complete the the proof of the Lemma it remains to check the inequalities (ii).
 
 First note that for any Seifert piece $S$ of $G_{\rm thin}$ then  ${\rm Card}(\pi_0\b S)\geq 2$.
 Indeed it follows from the construction  of  $G_{\rm thin}$ that $f_{\ast}(\pi_1S)\simeq{\Z}\times{\Z}$ and from the non-degeneration condition that $f_{\ast}(\pi_1 T)\simeq{\Z}\times{\Z}$ for any component $T$ of $\b S$. Thus we get the following commutative diagram
 $$\xymatrix{
 \pi_1T \ar[r] \ar[d] & \pi_1S \ar[r] \ar[d] & {\Z}\times{\Z}\subset\pi_1\Sigma\\
 \H{T} \ar[r] & \H{S} \ar[ur]
 }$$
This implies that ${\rm Rk}(\H{T}\to\H{S})=2$. If $\b S$ is connected then it follows from the  exact sequence corresponding to the pair $(S,\b S)$ that ${\rm Rk}(\H{\b S}\to\H{S})=1$ then  $\b S$ can not be connected.

Let $e$ be an element of $E(\hat{\Gamma}"_0)$ with end points $x_1$ and $x_2$. Notice that $x_i\in r\circ q(V(\Gamma))$ for $i=1,2$.  For the points $x_1,x_2$ the following  possibilities hold: 

\emph{Case 1:} Assume that  $x_1$ and $x_2$ correpond to Seifert pieces of $G^{\rm eff,\sim}_{\rm thin,hor}\cup G^{\rm eff}_{\rm thin,ver}$.
Let $T$ denote the characteristic torus corresponding to $e$. Let $S_1$ and $S_2$ denote the Seifert pieces of $G\setminus{\mathcal T}_G$ adjacent to $T$ ($S_i\not\simeq{\S}^1\times{\S}^1\times I$ for $i=1,2$). 

\emph{Subcase 1.1:} Suppose first that $S_1=S_2=S$ then it follows from our construction (after re-indexing) that $x_1$ corresponds to $S_1$ and $x_2$ correpsonds to a piece homeomorphic to ${\S}^1\times{\S}^1\times I$ (see Remark \ref{epaisse}). But in this case, $v(x_2)=2$. Let $e'$ the edge of $\hat{\Gamma}"_0$ adjacent to $e$ along $x_2$ with end points $\{x_2,x_3\}$. Necessarily $x_3$ corresponds to $S$. Then one can replace the two edges $e,e'$ by a single edge $e\cup_{x_2}e'$. Note that by construction $x_1=x_3$ and thus $e\cup_{x_2}e'$ is a simple closed curve. Then one can remove the curve $e\cup_{x_2}e'$ so that properties (1), (2) and (3) remain true, using the same arguments as in the proof of Claim \ref{valence}.       

\emph{Subcase 1.2:} Suppose that $S_1\not=S_2$. Then $e\in E_1$.
This comes from the fact that to Seifert pieces of $G^{\rm eff}_{\rm thin, ver}$ can not be adjacent along the characteristic torus corresponding to $e$ by minimality of the JSJ-decomposition. 

\emph{Case 2:} If $x_1\in\o{\mathcal O}^{\rm new}_{\rm main}$ and $x_2$ corresponds to a Seifert piece of $G^{\rm eff,\sim}_{\rm thin, hor}$ then $e\in E_1$.

\emph{Case 3:} To finish, suppose that  $x_1\in\o{\mathcal O}^{\rm new}_{\rm main}$ and $x_2$ corresponds to a Seifert piece $S$ of $G^{\rm eff}_{\rm thin, ver}$. If ${\rm Card}(\pi_0\b S)\geq 3$ then $e\in E_2$. Assume that ${\rm Card}(\pi_0\b S)=2$. Then $v(x_2)\leq 2$. If $v(x_2)=1$ then one can remove the edge $e$ from the graph $\hat{\Gamma}"_0$ so that properties (1), (2) and (3) remain true. If $v(x)=2$ then denote by $e'$ the edge of $\hat{\Gamma}"_0$ such that $e\cap e'=\{x_2\}$ and denote by $x_3$ the vertex of $e_3$ so that $\{x_2,x_3\}=\b e'$. Notice that it follows from our construction that $x_3$ necessarily correponds to a Seifert piece of $G^{\rm eff,\sim}_{\rm thin, hor}$. Then one can replace  the edges $e,e'$ by a single edge $e\cup_{x_2} e'$ so that $e\cup_{x_2} e'$ can be seen as an edge of $E_1$ in the new graph. We perform this operation for any edge in Case 3 and we denote by $\hat{\Gamma}_1$ the resulting graph. Notice that $\hat{\Gamma}_1$ satisfies conditions (1), (2) and (3).   Denote still by $E_1\cup E_2$ the decomposition of $E(\hat{\Gamma}_1)$. By condition (2) the inequality 
$${\rm Card}(E_1)\leq\sum_{S\in G^{\rm eff,\sim}_{\rm thin, hor}}{\rm Card}(\pi_0(\b S))$$
is clearly true. Thus it remains to check the second inequality. Let $S$ be a component of $G^{\rm eff}_{\rm thin,ver}$ and let $\{e_1,...,e_k\}$ denote the edges of $E_2$ corresponding to $\pi_0(\b S)$. We have to check that $k\leq{\rm Card}(\pi_0(\b S))-2$. If $k=0$ the result is obvious. Thus assume that $k\geq 1$. Then it follows from our construction that ${\rm Card}(\pi_0\b S)\geq 3$.  Note that the edges $\{e_1,...,e_k\}$ correspond to canonical tori which can be seen as boundary  components of some components of $G_{\rm main}$. But these boundary components are related by a sequence of vertical annuli in $S$ and so $k=1$. This completes the proof of the lemma.

\end{proof}

\subsubsection{Estimating the volume.} In this paragraph we check the inequality ${\rm Vol}(G)>{\rm Vol}(\Sigma)$ when $G_{\rm thin}\not=\emptyset$. We distinguish two casis.

\emph{First Case:} First of all assume that $G^{\rm eff,\sim}_{\rm thin,hor}=\emptyset$. Then in this case, it follows directly fom the proof of Lemma \ref{comparing} that $\star''\co\o{\mathcal O}^{\rm new}_{\rm main}\to F_{\Sigma}$ is an epimorphism at the $H_1$-level with coefficient ${\Q}$. This implies that ${\rm Vol}(G_{\rm main})\geq{\rm Vol}(\Sigma)$. Thus  ${\rm Vol}(G)>{\rm Vol}(\Sigma)$ (since the Seifert pieces of $G_{\rm main}$ can not be adjacent in $G$). 

\emph{Second Case:} We assume that $G^{\rm eff,\sim}_{\rm thin,hor}\not=\emptyset$. First suppose that ${\rm genus}({\mathcal O}_S)\geq 1$ for any $S\in G^{\rm eff,\sim}_{\rm thin,hor}$.
Consider the Mayer-Vietoris exact sequence corresponding to  the decomposition of $\hat{\Gamma}_1\cup\o{\mathcal O}^{\rm new}_{\rm main}$ given by  $\left(\hat{\Gamma}_1,  \o{\mathcal O}^{\rm new}_{\rm main}, {\mathcal C}\right)$ where ${\mathcal C}=\hat{\Gamma}_1\cap\o{\mathcal O}^{\rm new}_{\rm main}$. Denote by $S_1,...,S_k$ the components of $\o{\mathcal O}^{\rm new}_{\rm main}$ and by $\Sigma_1,...,\Sigma_l$ the components of  $G^{\rm eff,\sim}_{\rm thin,hor}$. Then we get
$$\{0\}\to H_1\left(\o{\mathcal O}^{\rm new}_{\rm main}\right)\oplus H_1\left(\hat{\Gamma}_1\right)\to H_1\left(\hat{\Gamma}_1\cup\o{\mathcal O}^{\rm new}_{\rm main}\right)\to H_0\left({\mathcal C}\right)\to ...$$
$$...\to H_0\left(\o{\mathcal O}^{\rm new}_{\rm main}\right)\oplus H_0\left(\hat{\Gamma}_1\right)\to H_0\left(\hat{\Gamma}_1\cup\o{\mathcal O}^{\rm new}_{\rm main}\right)\to\{0\}   $$
Thus we get the following relation
$$\beta_1\left(\o{\mathcal O}^{\rm new}_{\rm main}\right)=\beta_1\left(\hat{\Gamma}_1\cup\o{\mathcal O}^{\rm new}_{\rm main} \right)-\beta_1\left(\hat{\Gamma}_1\right)$$
On one hand we know that ${\rm Vol}\left(\o{\mathcal O}^{\rm new}_{\rm main}\cup\o{\mathcal O}^{\rm eff,\sim}_{\rm thin, hor}\right)=\beta_1\left(\o{\mathcal O}^{\rm new}_{\rm main}\right)+\beta_1\left(\o{\mathcal O}^{\rm eff, \sim}_{\rm thin,hor}\right)-k-l$  Thus we get
$${\rm Vol}\left(\o{\mathcal O}^{\rm new}_{\rm main}\cup\o{\mathcal O}^{\rm eff,\sim}_{\rm thin, hor}\right)=\beta_1\left(\hat{\Gamma}_1\cup\o{\mathcal O}^{\rm new}_{\rm main} \right)-\beta_1(\hat{\Gamma}_1)+\beta_1\left(\o{\mathcal O}^{\rm eff,\sim}_{\rm thin,hor}\right)
-k-l$$
 Moreover we know that $\beta_1\left(\hat{\Gamma}_1\right)={\rm Card}\left(E\left(\hat{\Gamma}_1\right)\right)-{\rm Card}\left(V\left(\hat{\Gamma}_1\right)\right)+{\rm Card}(\pi_0(\hat{\Gamma}_1))$.

By point (i) of Lemma \ref{comparing} we know that ${\rm Card}\left(V(\hat{\Gamma}_1)\right)\geq k+l+{\rm Card}(\pi_0(\hat{\Gamma}_1))-1$. 
On the other hand,  using point (iii) of Lemma \ref{comparing}, then  $\beta_1\left(\hat{\Gamma}_1\cup\o{\mathcal O}^{\rm new}_{\rm main}\right)\geq\beta_1\left(F_{\Sigma}\right)$ and ${\rm Vol}(\Sigma)=\beta_1\left(F_{\Sigma}\right)-\epsilon$ where $\epsilon=2$ or $1$ depending on whether $F_{\Sigma}$ is closed or not. This implies that 
$${\rm Vol}\left(\o{\mathcal O}^{\rm new}_{\rm main}\cup\o{\mathcal O}^{\rm eff,\sim}_{\rm thin, hor}\right)\geq{\rm Vol}(\Sigma)+(\epsilon-1)+\beta_1\left(\o{\mathcal O}^{\rm eff,\sim}_{\rm thin,hor}\right)-{\rm Card}\left(E\left(\hat{\Gamma}_1\right)\right)$$
By point (i) of Lemma \ref{comparing} one can decompose  $E\left(\hat{\Gamma}_1\right)$ into the union $E_1\cup E_2$ where $E_1$ consists of edges corresponding to components of $\b G^{\rm eff,\sim}_{\rm thin,hor}$ and $E_2$ consists of edges corresponding to components of $\b G^{\rm eff}_{\rm thin,ver}$ which are not adjacent to coponents of $G^{\rm eff,\sim}_{\rm thin,hor}$. Then we get
$${\rm Vol}\left(G^{\rm eff}\right)\geq{\rm Vol}(\Sigma)+(\epsilon-1)+\beta_1\left(\o{\mathcal O}^{\rm eff,\sim}_{\rm thin,hor}\right)-{\rm Card}\left(E_1\right)+{\rm Vol}\left(G^{\rm eff}_{\rm thin,ver}\right)-{\rm Card}\left(E_2\right)$$
First note that it follows from Lemma \ref{comparing} 
$${\rm Card}\left(E_2\right)\leq\sum_{S\in G^{\rm eff}_{\rm thin,ver}}\left({\rm Card}(\pi_0(\b S))-2\right)$$
and  ${\rm Vol}\left(G^{\rm eff}_{\rm thin,ver}\right)\geq\sum_{S\in G^{\rm eff}_{\rm thin,ver}}\left({\rm Card}(\pi_0(\b S))-2\right)$ then ${\rm Vol}\left(G^{\rm eff}_{\rm thin,ver}\right)\geq {\rm Card}\left(E_2\right)$ and thus we get the following inequality
$${\rm Vol}\left(G^{\rm eff}\right)\geq{\rm Vol}(\Sigma)+(\epsilon-1)+\beta_1\left(\o{\mathcal O}^{\rm eff,\sim}_{\rm thin,hor}\right)-{\rm Card}\left(E_1\right)
$$
Note that $\beta_1\left(\o{\mathcal O}^{\rm eff,\sim}_{\rm thin,hor}\right)=\sum_{i=1}^{l}\left(2g_i+r_i-1\right)$ where $g_i$, resp. $r_i$, denotes the genus, resp. the number boundary compoments,  of $\Sigma_i$, $i=1,...,l$. Then
$${\rm Vol}\left(G^{\rm eff}\right)\geq{\rm Vol}(\Sigma)+2\sum_{i=1}^{l}g_i-l+\sum_{i=1}^{l}r_i-{\rm Card}\left(E_1\right)$$
Using  point (ii) of Lemma \ref{comparing} we know that $\sum_{i=1}^{l}r_i-{\rm Card}\left(E_1\right)\geq 0$ and since, $g_i\geq 1$ for $i=1,...,l$ then we get 
$${\rm Vol}\left(G^{\rm eff}\right)\geq{\rm Vol}(\Sigma)+2\sum_{i=1}^{l}g_i-l>{\rm Vol}(\Sigma)$$ 
This proves that ${\rm Vol}(G^{\rm eff})>{\rm Vol}(\Sigma)$ since $l\geq 1$ by hypothesis. Hence this
 completes the proof in this case. If the condition on the genus of ${\mathcal O}_{\rm thin,hor}^{\rm eff,\sim}$ is not satisfied then, since condition $({\mathcal C})$ is satisfied,  we know from Lemma \ref{genus} that there exists a finite regular covering $f_1\co G_1\to\Sigma_1$ of $f\co G\to \Sigma$  satisfying the following properties:
let $\pi\co\Sigma_1\to\Sigma$ and $p\co G_1\to G$ denote the finite regular  coverings corresponding to $f_1$ then

(i) any geometric component of $p^{-1}(G^{\rm eff,\sim}_{\rm thin,hor})$ admits a Seifert fibration over a 2-orbifold of genus at least 1,

(ii)  for any geometric piece $S$ of  $G^{\rm eff}$ and for any component $S_1$ of $p^{-1}(S)$ then $G_h(p|S_1)\geq G_h(\pi)$. 

One can apply the above arguments to the map $f_1\co G_1\to\Sigma_1$.
First note that it follows from our construction that $(G_1)^{\rm eff}\subset p^{-1}(G^{\rm eff})$. Using point (i), it follows from the paragraph above that we have ${\rm Vol}((G_1)^{\rm eff})>{\rm Vol}(\Sigma_1)$. Thus we get
$${\rm Vol}\left(\Sigma_1\right)={\rm Vol}\left(\Sigma\right)\frac{{\rm deg}(\pi)}{G_h(\pi)}<{\rm Vol}\left((G_1)^{\rm eff}\right)\leq{\rm Vol}\left(p^{-1}\left(G^{\rm eff}\right)\right)$$
Denote by $Q_1,...,Q_l$ the geometric components of $G^{\rm eff}$ and by $p_i$ the induced covering $p|p^{-1}(Q_i)\co p^{-1}(Q_i)\to Q_i$. Then since $p$ is a regular covering we have
$${\rm Vol}\left(G^{\rm eff}\right)=\frac{1}{{\rm deg}(p)}\sum_{i=1}^{i=l}{\rm Vol}\left(p^{-1}\left(Q_i\right)\right)G_h(p_i)$$ 
Since $f\co G\to\Sigma$
   is $\pi_1$-surjective then ${\rm deg}(\pi)={\rm deg}(p)$ and by condition (ii) we get
$${\rm Vol}\left(G^{\rm eff}\right)\geq\frac{G_h(\pi)}{{\rm deg}(\pi)}{\rm Vol}\left(p^{-1}\left(G^{\rm eff}\right)\right)$$
By combining this latter inequality with the first one we get ${\rm Vol}(\Sigma)<{\rm Vol}(G^{\rm eff})$. The proof of inequality $({\mathcal V})$ is now complete. This ends the proof of Proposition \ref{fine}.

\section{Proof of the Theorems and Corollary} 
\subsection{Nonzero degree maps decreases the volume}
In order to prove Theorem \ref{decroissance} we state the following result which will be used in the proof of Theorem \ref{rigide} (see section \ref{enfin}).
\begin{lemma}\label{decroissance'}
Let $f\co M\to N$ be a nonzero degree map between closed Haken manifolds satisfying $\|M\|={\rm deg}(f)\|N\|$. Then ${\rm Vol}(M)\geq{\rm Vol}(N)$ and if there exists a canonical torus $T\in{\mathcal T}_M$ such that $f|T\co T \to N$  is not $\pi_1$-injective then ${\rm Vol}(M)>{\rm Vol}(N)$. 
\end{lemma}
\begin{proof}
First assume that $\tau(N)=0$. If $\tau(M)=0$ then $M$ is a virtual torus bundle and then $f$ is homotopic to a finite covering by \cite{W}, in particular $f_{\ast}\co\pi_1M\to\pi_1N$ is injective. In the other cases $\tau(M)\not=0$ and thus ${\rm Vol}(M)>0$. Thus from now one one can assume $\tau(N)\not=0$.
 
Suppose that $f|{\mathcal T}_M\co {\mathcal T}_M\to N$ is $\pi_1$-injective.  If $M_{\rm thin}\not=\emptyset$ then one can applies Proposotion \ref{fine}.  If $M_{\rm thin}=\emptyset$  then for any Seifert piece $\Sigma$  of $N$ that is not homeomorphic to ${\bf K}^2\t{\times}I$ each component of $f^{-1}(\Sigma)$ is a Seifert piece of $M$. This follows from Lemma \ref{adjacent} since $\Sigma$ has a hyperbolic 2-orbifold base. Choose a component $G$ of $f^{-1}(\Sigma)$ so that $f|G\co G\to\Sigma$ has nonzero degree. Thus we get ${\rm Vol}(G)\geq{\rm Vol}(\Sigma)$. This proves the lemma if $f$ is non-degenerate when restricted to ${\mathcal T}_M$.

 Suppose that $f|{\mathcal T}_M\co {\mathcal T}_M\to N$ is not $\pi_1$-injective.
Passing to a finite covering, we may assume, by Lemma \ref{huge} and Claim \ref{debut} that $N$ contains no embedded Klein bottles and thus using   Lemma \ref{ducon} we know that there exists a closed Haken manifold $\hat{M}_1$ and a map $f_1\co\hat{M}_1\to N$  satisfying the following conditions:
 
 (i)  $\hat{M}_1$ is obtained from a canonical submanifold $M_1$ of $M$ after Seifert Dehn fillings. This means in particular that ${\mathcal H}(\hat{M}_1)={\mathcal H}(M)$ and that each Seifert piece $\hat{S}$ of $\hat{M}_1$ is obtained as an extension from a unique Seifert piece  $S$ of $M$ after Seifert Dehn fillings,
 
 (ii) $\hat{M}_1$ and $M$ have the same Gromov simplicial volume and  ${\rm deg}(f_1)={\rm deg}(f)$,
 
 (iii) the map $f|{\mathcal T}_{\hat{M}_1}\co{\mathcal T}_{\hat{M}_1}\to N$ is non-degenerate. 
 
Since $f|{\mathcal T}_M\co{\mathcal T}_M\to N$ is degenerate, then there exists at leat one Seifert piece $\hat{S}$ in $\hat{M}_1$ obtained from $S$ after non-trivial (i.e. with slope $\not=\infty$) Seifert Dehn fillings. Assume that a Seifert fibration of $S$ is fixed. The base 2-orbifold ${\mathcal O}_{\hat{S}}$ of $\hat{S}$ is ${\mathcal O}_S$ after gluing some cone points along some components of $\b{\mathcal O}_{S}$.    Note that $S$ necessarily supports a ${\Hi}^2\times{\R}$-geometry.  

Indeed if not then $S$ is the twisted $I$-bundle over the Klein bottle and thus $\hat{M}_1=\hat{S}$ is a closed Seifert fibered space whose base is a 2-sphere with cone points $(2,2,n)$. Then $\hat{M}_1$ is a Seifert fibered space whose base 2-orbifold admits a spherical geometry. This contradicts the fact that $\hat{M}_1$ is a Haken manifold.  

Then we get $\chi({\mathcal O}_S)<\chi(\mathcal O_{\hat{S}})\leq 0$.
 This proves that ${\rm Vol}(\hat{M}_1)<{\rm Vol}(M)$. On the other hand since $f_1\co\hat{M}_1\to N$ has nonzero degree and since $\|\hat{M}_1\|={\rm deg}(f_1)\|N\|$ then ${\rm Vol}(\hat{M}_1)\geq{\rm Vol}(N)$ by the first case. This completes the proof of the lemma and of Theorem \ref{decroissance}. 
 \end{proof}

\subsection{Proof of the rigidity theorem}\label{enfin}
In this paragraph we prove Theorem \ref{rigide}.
Let $f\co M\to N$ be a nonzero degree map between closed Haken manifolds  satisfying the Volume Condition $\|M\|=\abs{{\rm deg}(f)}\|N\|$ and ${\rm Vol}(M)={\rm Vol}(N)$. Then it follows from Lemma \ref{decroissance'} that $f|{\mathcal T}_M$ is $\pi_1$-injective.  
\subsubsection{Assume that $N$ admits a geometry $E^3$, $Nil$ or $Sol$} This means that  $\tau(M)=\tau(N)=(0,0)$. Then $M$ is a virtual torus bundle (in particular $M$ is geometric) and since $N$ is irreducible then  $f$ is homotopic to a ${\rm deg}(f)$-fold covering by a result of \cite{W}.

\subsubsection{Assume that $N$ admits a geometry ${\Hi}^2\times{\R}$ or  $\t{SL}(2,{\R})$} Then we claim that $M$ and $N$ are both Seifert fibered spaces. 

Indeed, if not then ${\mathcal T}_M\not=\emptyset$. By Lemma \ref{adjacent} we know that $M_{\rm thick}\not=\emptyset$ and since ${\mathcal T}_M\not=\emptyset$ and since $f|{\mathcal T}_M$ is $\pi_1$-injective then $M_{\rm thin}\not=\emptyset$. This implies using Proposition \ref{fine} that ${\rm Vol}(M)>{\rm Vol}(N)$. A contradiction. Thus we may assume that $M$ is Seifert.

 Note that $f$ is homotopic to a fiber preserving map.  Let $q\co\hat{N}\to N$ be the finite covering of $N$ correpsonding to $f_{\ast}(\pi_1M)$ and let $\hat{f}\co M\to\hat{N}$ denote the lifting of $f$. There exists  a finite covering $\t{f}\co\t{M}\to\t{N}$ of $\hat{f}$ such that $\t{M}\to M$ and $\t{N}\to\hat{N}$ have fiber degree $\pm 1$ and such that $\t{N}$ is a ${\S}^1$-bundle over a closed orientable hyperbolic surface $\t{F}$. Note that it follows from our construction that ${\rm Vol}(\t{M})={\rm Vol}(\t{N})$. 
  Then the map $\t{f}$ descends to a nonzero degree map $\pi\co\o{\mathcal O}_{\t{M}}\to\t{F}$, where $\o{\mathcal O}_{\t{M}}$ denotes the base surface of $\t{M}$. Note that $-\chi({\mathcal O}_{\t{M}})\geq-\chi(\o{\mathcal O}_{\t{M}})\geq{\rm deg}(\pi)(-\chi(\t{F}))>0$ and since ${\rm Vol}(\t{M})={\rm Vol}(\t{N})$ then $\chi({\mathcal O}_{\t{M}})=\chi(\o{\mathcal O}_{\t{M}})=\chi(\t{F})<0$ and thus $\t{M}$ is a ${\S}^1$-bundle over a closed orientable hyperbolic surface $\t{K}=\o{\mathcal O}_{\t{M}}={\mathcal O}_{\t{M}}$ and ${\rm deg}(\pi)=1$ which implies that $\pi\co\t{K}\to\t{F}$ is homotopic to a homeomorphism. Denote by $h$ (resp. $t$) the homotopy class of the fiber in $\t{M}$ (in $\t{N}$ resp.) and let $n$ denote the nonzero integer such that $f_{\ast}(h)=t^n$. Using the exact sequence
  $$\xymatrix{
  \{1\} \ar[r] \ar[d] & {\Z} \ar[r] \ar[d]_{\times n} & \pi_1(\t{M}) \ar[r] \ar[d]_{\t{f}_{\ast}} & \pi_1(K) \ar[r] \ar[d]_{\pi_{\ast}} & \{1\} \ar[d]\\
  \{1\} \ar[r]  & {\Z} \ar[r]  & \pi_1(\t{N}) \ar[r] \ & \pi_1(F) \ar[r]  & \{1\}
  }$$
   we check that $\t{f}_{\ast}$ is an isomorphism. Thus so is $\hat{f}$, by \cite{Wa}, and finally $f$ is a covering map. Moreover we claim  that  $G_h(f)={\rm deg}(f)$ and $G_{\rm ob}(f)=1$.  Indeed the map $f$ induces a map $f'\co{\mathcal O}_M\to{\mathcal O}_N$ of degree $G_{\rm ob}(f)$. This implies that $\abs{\chi({\mathcal O}_M)}\geq G_{\rm ob}(f)\abs{\chi({\mathcal O}_N)}>0$ and since ${\rm Vol}(M)=\abs{\chi({\mathcal O}_M)}={\rm Vol}(N)=\abs{\chi({\mathcal O}_N)}$ then $G_{\rm ob}(f)=1$ and since ${\rm deg}(f)=G_h(f)\times G_{\rm ob}(f)$ our claim is checked.
   
\subsubsection{Assume that $N$ is hyperbolic} In this case the condition on the volume implies that $M$ is still a hyperbolic manifold and  $f$ is homotopic to a covering map by a rigidity result of T. Soma in \cite[Theorem 1]{Srigide}.   

\subsubsection{Assume  that $N$ is a non-geometric Haken manifold} This means in particular that $\tau(N)\not=0$. Let $q\co\hat{N}\to N$ be the finite covering of $N$ corresponding to $f_{\ast}(\pi_1M)$ and let $\hat{f}\co M\to\hat{N}$ denote the lifting of $f$. There exists  a finite covering $\t{f}\co\t{M}\to\t{N}$ of $\hat{f}$ acting trivially on ${\mathcal T}_N\cup{\mathcal H}(N)$ (resp. ${\mathcal T}_M\cup{\mathcal H}(M)$)  (in particular they have fiber degree 1) and such that $\t{N}$ contains no embedded Klein bottles.
 After adjusting $\t{f}\co\t{M}\to\t{N}$ in standard form, using Corollary \ref{chienne} we fix a Seifert piece $\Sigma$ in $\t{N}$ and consider a component $G$ of $\t{f}^{-1}(\Sigma)$ so that ${\rm deg}(\t{f}|G\co G\to\Sigma)\not=0$. We know from Proposition \ref{fine} that ${\rm Vol}(G)\geq{\rm Vol}(\Sigma)>0$. This implies, since ${\rm Vol}(\t{M})={\rm Vol}(\t{N})$ that $\t{f}^{-1}(\Sigma)$ is actually connected and equal to $G$. On the other hand, we know that if $G_{\rm thin}\not=\emptyset$ then ${\rm Vol}(G)>{\rm Vol}(\Sigma)$. Then $G_{\rm thin}=\emptyset$ and thus $G$ is a Seifert piece of $\t{M}$ with $\chi({\mathcal O}_G)=\chi({\mathcal O}_{\Sigma})$. Hence using paragraph 5.2.2 we know that $\t{f}|G$ is a covering map such that  $G_h(\t{f}|G)={\rm deg}(\t{f}|G)$ and $G_{\rm ob}(f|G)=1$. This proves that $\t{f}|{\mathcal S}(\t{M})\co{\mathcal S}(\t{M})\to{\mathcal S}(\t{N})$ is a covering map. On the other hand  $\t{f}|{\mathcal H}(\t{M})\co{\mathcal H}(\t{M})\to{\mathcal H}(\t{N})$ is a covering map by a result of T. Soma in \cite{Srigide}. But since $\t{f}$ is $\pi_1$-surjective then $\t{f}$  and thus $\hat{f}$ is actually a homeomorphism, using \cite{Wa}, and hence $f$ is a covering map. Note that the induced proper map $f|{\mathcal S}(M)\co{\mathcal S}(M)\to{\mathcal S}(N)$ is a covering map  such that $G_h(f|{\mathcal S}_h(M))={\rm deg}(f)$ and $G_{\rm ob}(f|{\mathcal S}_h(M))=1$. This completes the proof of Theorem \ref{rigide}.

\subsection{Proof of Corollary \ref{strong rigidity}}
We consider here  degree one maps  between closed Haken manifolds. In view of Theorem \ref{rigide}, to prove Corollary \ref{strong rigidity} we have to check the following
\begin{claim}\label{equality} For any closed Haken manifold $M$ 
 there exists a constant $\eta_M\in (0,1)$, which depends only on $M$, such that if $N$ is a closed Haken manifold 1-dominated by $M$ and satisfying  $\tau(N)\geq\tau(M)(1-\eta_M)$ then $\tau(M)=\tau(N)$.
 \end{claim}
 \begin{proof} 
Suppose the contrary. Then there is a closed Haken manifold $M_0$ and a sequence of closed Haken manifolds $N_n$  such that there are degree one maps $f_n\co M_0\to N_n$ satisfying $\tau(N_n)\geq\tau(M_0)(1-1/n)$ and $\tau(N_n)\not=\tau(M_0)$ for any $n\in{\N}$.  This implies in particular that $\|M_0\|\geq\|N_n\|\geq\|M_0\|(1-1/n)$. Then $\lim_{n\to\infty}\|N_n\|=\|M_0\|$. Hence by \cite{D} this implies that the sequence $\{N_n\}_{n\in{\N}}$ is finite up to homeomorphism. This contradicts the inequalities 
$$\|M_0\|\left(1-\frac{1}{n}\right)\leq\tau\left(N_n\right)<\tau\left(M_0\right)$$
This completes the proof of the claim.
\end{proof}
Thus one can apply Theorem \ref{rigide} with the hypothesis ${\rm deg}(f)=1$. This completes the proof of the corollary.


\begin{thebibliography}{ZZZZZZ}

\bibitem[Al]{Al} {\sc Alperin, Roger C.}, {\it 
An elementary account of Selberg's lemma},
Enseign. Math. (2) 33 (1987), no. 3-4, 269--273.

\bibitem[BW]{BW} {\sc M. Boileau, S. Wang}, {\it Non-zero degree maps and surface bundles over ${\S}^1$}, J. Differential Geom.  43  (1996),  no. 4, 789--806.

\bibitem[D]{D} {\sc P. Derbez}, {\it Non-zero degree maps between closed orientable three-manifolds}, 
    To appear in  Trans. Amer. Math. Soc (2005).

\bibitem[G]{G} {\sc M. Gromov}, {\it Volume and bounded cohomology},  Inst. Hautes \'{E}tudes Sci. Publ. Math. No. 56, (1982), 5--99.

\bibitem[Gr]{Gr} {\sc  K.  Gruenberg}, {\it Residual properties of infinite soluble groups},  Proc. London Math. Soc. (3)  7  (1957), 29--62.

\bibitem[H]{H} {\sc E. Hamilton},  {\it Abelian subgroup separability of Haken 3-manifolds and closed hyperbolic $n$-orbifolds},   Proc. London Math. Soc. (3)  83  (2001),  no. 3, 626--646.


 \bibitem[JS]{JS} {\sc W. Jaco, P.B. Shalen},  {\it Seifert fibered space in
  3-manifolds},  Mem. Amer. Math. Soc.  21  (1979).

\bibitem[J]{J}  {\sc K. Johannson}, {\it Homotopy equivalences of 3-manifolds with
 boundaries},  Lecture Notes in Mathematics, 761. Springer, Berlin, 1979.
 
 
 \bibitem[PS]{PS} {\sc B. Perron, P. Shalen}, {\it  Homeomorphic graph manifolds: a contribution to the $µ$ constant problem},   Topology Appl.  99  (1999),  no. 1, 1--39. 
 
 \bibitem[Ro]{Ro} {\sc Y.  Rong}, {\it Maps between Seifert fibered spaces of infinite $\pi\sb 1$},  Pacific J. Math.  160  (1993),  no. 1, 143--154.
 
 \bibitem[Sc]{Sc} {\sc P. Scott}, {\it Subgroups of surface groups are almost geometric},  J. London Math. Soc. (2)  17  (1978), no. 3, 555--565

\bibitem[S1]{Srigide} {\sc T. Soma},  {\it A rigidity Theorem for Haken manifolds},  Math. Proc. Cambridge Philos. Soc.  118  (1995),  no. 1, 141--160.

\bibitem[S2]{Smostow} {\sc T. Soma},  {\it Degree-one maps between hyperbolic 3-manifolds with the same volume limit},  Trans. Amer. Math. Soc.  353  (2001),  no. 7, 2753--2772.

\bibitem[S3]{Sdom} {\sc T. Soma}, {\it Sequences of degree-one maps between geometric
  3-manifolds},  Math. Ann.  316  (2000),  no. 4, 733--742.

 \bibitem[Th]{Th} {\sc W. Thurston}, {\it The geometry and topology of 3-manifolds}, Lectures Notes, Princeton Univ., 1979.
 
\bibitem[Wa]{Wa} {\sc F. Waldhausen}, {\it On irreducible 3-manifolds which are sufficiently large}, Ann. of Math. 87 (1968), 56-88.

\bibitem[W]{W} {\sc S. Wang}, {\it  The existence of maps of nonzero degree between aspherical $3$-manifolds}, 
Math. Z. 208 (1991), no. 1, 147--160.

\bibitem[W1]{W1} {\sc S. Wang}, {\it  The $\pi\sb 1$-injectivity of self-maps of nonzero degree on $3$-manifolds},  Math. Ann.  297  (1993),  no. 1, 171--189.


 
\end{thebibliography}
\end{document}